\documentclass[12pt] {article}
\usepackage{amsmath,amsthm}
\usepackage{amssymb,latexsym}

\usepackage{fullpage}
\usepackage{amscd}
\usepackage{epic,eepic}

\usepackage[dvips]{graphicx}
\usepackage[all]{xy}
\usepackage{amsfonts}
\usepackage{amsmath}
\usepackage{amssymb}
\usepackage{amscd}
\usepackage{color}
\usepackage{amsthm}
\usepackage{indentfirst}
\usepackage[hmargin=3cm,vmargin=3cm]{geometry}
\usepackage[small,nohug,heads=littlevee]{diagrams}

\title{A uniform estimate for general quaternionic Calabi problem (with appendix by Daniel Barlet).}
\date{\today}

\author{Semyon Alesker\footnote{Partially supported by ISF grant 1447/12.}\\
{ \normalsize Department of Mathematics}\\
{ \normalsize Tel Aviv University, Ramat Aviv}\\
{ \normalsize 69978 Tel Aviv, Israel}\\
{ \normalsize e-mail: semyon@post.tau.ac.il}
\and
 Egor Shelukhin\\
{ \normalsize  CRM}\\
{ \normalsize  University of Montreal}\\
{ \normalsize e-mail: egorshel@gmail.com}}

\def\RR{\mathbb{R}}
\def\CC{\mathbb{C}}

\def\NN{\mathbb{N}}
\def\ZZ{\mathbb{Z}}
\def\HH{\mathbb{H}}

\def\PP{\mathbb{P}}

\def\eps{\varepsilon}
\def\alp{\alpha}
\def\Ome{\Omega}
\def\ome{\omega}
\def\lam{\lambda}
\def\Lam{\Lambda}

\def\to{\longrightarrow}

\newcommand{\del}{\partial}

\newcommand{\sm}[1]{C^\infty(#1)}
\newcommand{\av}[1]{\langle#1\rangle}
\newcommand{\delbar}{\overline{\partial}}
\newcommand{\Sum}{\Sigma}

\newcommand{\normL}[3]{||#1||_{L^{#2}(#3)}}
\newcommand{\normC}[3]{||#1||_{C^{#2}(#3)}}

\newcommand{\R}{{\mathbb{R}}}

\newcommand{\C}{{\mathbb{C}}}

\newcommand{\om}{\omega}
\newcommand{\Om}{\Omega}

\newcommand{\ep}{\epsilon}

\def\qed { Q.E.D. }

\swapnumbers
\newtheorem{theorem}{Theorem}[subsection]
\newtheorem{corollary}[theorem]{Corollary}
\newtheorem{lemma}[theorem]{Lemma}
\newtheorem{proposition}[theorem]{Proposition}
\newtheorem{claim}[theorem]{Claim}
\theoremstyle{definition}

\newtheorem{definition}[theorem]{Definition}
\newtheorem{remark}[theorem]{Remark}
\def\ca{{\cal A}} \def\cb{{\cal B}} 
  \def\cf{{\cal F}}
\def\cg{{\cal G}}  \def\ci{{\cal I}}
\def\cj{{\cal J}} \def\ck{{\cal K}} \def\cl{{\cal L}}
 \def\cn{{\cal N}} \def\co{{\cal O}}
 
\def\car{{\cal R}}
\def\cs{{\cal S}} \def\ct{{\cal T}} \def\cu{{\cal U}}
\def\cv{{\cal V}}  
 
\def\pt{\partial}
\def\inj{\hookrightarrow}
\def\BC{\mathbb{B}\mathbb{C}^\bullet}

\makeatletter

\def\diagram{\m@th\leftwidth=\z@ \rightwidth=\z@ \topheight=\z@
\botheight=\z@ \setbox\@picbox\hbox\bgroup}

\def\enddiagram{\egroup\wd\@picbox\rightwidth\unitlength
\ht\@picbox\topheight\unitlength \dp\@picbox\botheight\unitlength
\hskip\leftwidth\unitlength\box\@picbox}

\def\bfig{\begin{diagram}}
\def\efig{\end{diagram}}
\newcount\wideness \newcount\leftwidth \newcount\rightwidth
\newcount\highness \newcount\topheight \newcount\botheight

\def\ratchet#1#2{\ifnum#1<#2 \global #1=#2 \fi}

\def\putbox(#1,#2)#3{%
\horsize{\wideness}{#3} \divide\wideness by 2
{\advance\wideness by #1 \ratchet{\rightwidth}{\wideness}}
{\advance\wideness by -#1 \ratchet{\leftwidth}{\wideness}}
\vertsize{\highness}{#3} \divide\highness by 2
{\advance\highness by #2 \ratchet{\topheight}{\highness}}
{\advance\highness by -#2 \ratchet{\botheight}{\highness}}
\put(#1,#2){\makebox(0,0){$#3$}}}

\def\putlbox(#1,#2)#3{%
\horsize{\wideness}{#3}
{\advance\wideness by #1 \ratchet{\rightwidth}{\wideness}}
{\ratchet{\leftwidth}{-#1}}
\vertsize{\highness}{#3} \divide\highness by 2
{\advance\highness by #2 \ratchet{\topheight}{\highness}}
{\advance\highness by -#2 \ratchet{\botheight}{\highness}}
\put(#1,#2){\makebox(0,0)[l]{$#3$}}}

\def\putrbox(#1,#2)#3{%
\horsize{\wideness}{#3}
{\ratchet{\rightwidth}{#1}}
{\advance\wideness by -#1 \ratchet{\leftwidth}{\wideness}}
\vertsize{\highness}{#3} \divide\highness by 2
{\advance\highness by #2 \ratchet{\topheight}{\highness}}
{\advance\highness by -#2 \ratchet{\botheight}{\highness}}
\put(#1,#2){\makebox(0,0)[r]{$#3$}}}

\def\adjust[#1]{} 

\newcount \coefa
\newcount \coefb
\newcount \coefc
\newcount\tempcounta
\newcount\tempcountb
\newcount\tempcountc
\newcount\tempcountd
\newcount\xext
\newcount\yext
\newcount\xoff
\newcount\yoff
\newcount\gap%
\newcount\arrowtypea
\newcount\arrowtypeb
\newcount\arrowtypec
\newcount\arrowtyped
\newcount\arrowtypee
\newcount\height
\newcount\width
\newcount\xpos
\newcount\ypos
\newcount\run
\newcount\rise
\newcount\arrowlength
\newcount\halflength
\newcount\arrowtype
\newdimen\tempdimen
\newdimen\xlen
\newdimen\ylen
\newsavebox{\tempboxa}%
\newsavebox{\tempboxb}%
\newsavebox{\tempboxc}%

\newdimen\w@dth

\def\setw@dth#1#2{\setbox\z@\hbox{\m@th$#1$}\w@dth=\wd\z@
\setbox\@ne\hbox{\m@th$#2$}\ifnum\w@dth<\wd\@ne \w@dth=\wd\@ne \fi
\advance\w@dth by 1.2em}

\def\t@^#1_#2{\allowbreak\def\n@one{#1}\def\n@two{#2}\mathrel
{\setw@dth{#1}{#2}
\mathop{\hbox to \w@dth{\rightarrowfill}}\limits
\ifx\n@one\empty\else ^{\box\z@}\fi
\ifx\n@two\empty\else _{\box\@ne}\fi}}
\def\t@@^#1{\@ifnextchar_{\t@^{#1}}{\t@^{#1}_{}}}
\def\to{\@ifnextchar^{\t@@}{\t@@^{}}}

\def\t@left^#1_#2{\def\n@one{#1}\def\n@two{#2}\mathrel{\setw@dth{#1}{#2}
\mathop{\hbox to \w@dth{\leftarrowfill}}\limits
\ifx\n@one\empty\else ^{\box\z@}\fi
\ifx\n@two\empty\else _{\box\@ne}\fi}}
\def\t@@left^#1{\@ifnextchar_{\t@left^{#1}}{\t@left^{#1}_{}}}
\def\toleft{\@ifnextchar^{\t@@left}{\t@@left^{}}}

\def\two@^#1_#2{\allowbreak
\def\n@one{#1}\def\n@two{#2}\mathrel{\setw@dth{#1}{#2}
\mathop{\vcenter{\lineskip\z@\baselineskip\z@
                 \hbox to \w@dth{\rightarrowfill}%
                 \hbox to \w@dth{\rightarrowfill}}%
       }\limits
\ifx\n@one\empty\else ^{\box\z@}\fi
\ifx\n@two\empty\else _{\box\@ne}\fi}}
\def\tw@@^#1{\@ifnextchar _{\two@^{#1}}{\two@^{#1}_{}}}
\def\two{\@ifnextchar ^{\tw@@}{\tw@@^{}}}

\def\tofr@^#1_#2{\def\n@one{#1}\def\n@two{#2}\mathrel{\setw@dth{#1}{#2}
\mathop{\vcenter{\hbox to \w@dth{\rightarrowfill}\kern-1.7ex
                 \hbox to \w@dth{\leftarrowfill}}%
       }\limits
\ifx\n@one\empty\else ^{\box\z@}\fi
\ifx\n@two\empty\else _{\box\@ne}\fi}}
\def\t@fr@^#1{\@ifnextchar_ {\tofr@^{#1}}{\tofr@^{#1}_{}}}
\def\tofro{\@ifnextchar^ {\t@fr@}{\t@fr@^{}}}

\def\mon{\mathop{\m@th\hbox to
      14.6\P@{\lasyb\char'51\hskip-2.1\P@$\arrext$\hss
$\mathord\rightarrow$}}\limits} 
\def\leftmono{\mathrel{\m@th\hbox to
14.6\P@{$\mathord\leftarrow$\hss$\arrext$\hskip-2.1\P@\lasyb\char'50%
}}\limits} 
\mathchardef\arrext="0200       

\def\settypes(#1,#2,#3){\arrowtypea#1 \arrowtypeb#2 \arrowtypec#3}
\def\settoheight#1#2{\setbox\@tempboxa\hbox{#2}#1\ht\@tempboxa\relax}%
\def\settodepth#1#2{\setbox\@tempboxa\hbox{#2}#1\dp\@tempboxa\relax}%
\def\settokens`#1`#2`#3`#4`{%
     \def\tokena{#1}\def\tokenb{#2}\def\tokenc{#3}\def\tokend{#4}}
\def\setsqparms[#1`#2`#3`#4;#5`#6]{%
\arrowtypea #1
\arrowtypeb #2
\arrowtypec #3
\arrowtyped #4
\width #5
\height #6
}
\def\setpos(#1,#2){\xpos=#1 \ypos#2}

\def\settriparms[#1`#2`#3;#4]{\settripairparms[#1`#2`#3`1`1;#4]}%

\def\settripairparms[#1`#2`#3`#4`#5;#6]{%
\arrowtypea #1
\arrowtypeb #2
\arrowtypec #3
\arrowtyped #4
\arrowtypee #5
\width #6
\height #6
}

\def\resetparms{\settripairparms[1`1`1`1`1;500]\width 500}

\resetparms

\def\mvector(#1,#2)#3{
\put(0,0){\vector(#1,#2){#3}}%
\put(0,0){\vector(#1,#2){26}}%
}
\def\evector(#1,#2)#3{{
\arrowlength #3
\put(0,0){\vector(#1,#2){\arrowlength}}%
\advance \arrowlength by-30
\put(0,0){\vector(#1,#2){\arrowlength}}%
}}

\def\horsize#1#2{%
\settowidth{\tempdimen}{$#2$}%
#1=\tempdimen
\divide #1 by\unitlength
}

\def\vertsize#1#2{%
\settoheight{\tempdimen}{$#2$}%
#1=\tempdimen
\settodepth{\tempdimen}{$#2$}%
\advance #1 by\tempdimen
\divide #1 by\unitlength
}

\def\putvector(#1,#2)(#3,#4)#5#6{{%
\ifnum3<\arrowtype
\putdashvector(#1,#2)(#3,#4)#5\arrowtype
\else
\ifnum\arrowtype<-3
\putdashvector(#1,#2)(#3,#4)#5\arrowtype
\else
\xpos=#1
\ypos=#2
\run=#3
\rise=#4
\arrowlength=#5
\ifnum \arrowtype<0
    \ifnum \run=0
        \advance \ypos by-\arrowlength
    \else
        \tempcounta \arrowlength
        \multiply \tempcounta by\rise
        \divide \tempcounta by\run
        \ifnum\run>0
            \advance \xpos by\arrowlength
            \advance \ypos by\tempcounta
        \else
            \advance \xpos by-\arrowlength
            \advance \ypos by-\tempcounta
        \fi
    \fi
    \multiply \arrowtype by-1
    \multiply \rise by-1
    \multiply \run by-1
\fi
\ifcase \arrowtype
\or \put(\xpos,\ypos){\vector(\run,\rise){\arrowlength}}%
\or \put(\xpos,\ypos){\mvector(\run,\rise)\arrowlength}%
\or \put(\xpos,\ypos){\evector(\run,\rise){\arrowlength}}%
\fi\fi\fi
}}

\def\putsplitvector(#1,#2)#3#4{
\xpos #1
\ypos #2
\arrowtype #4
\halflength #3
\arrowlength #3
\gap 140
\advance \halflength by-\gap
\divide \halflength by2
\ifnum\arrowtype>0
   \ifcase \arrowtype
   \or \put(\xpos,\ypos){\line(0,-1){\halflength}}%
       \advance\ypos by-\halflength
       \advance\ypos by-\gap
       \put(\xpos,\ypos){\vector(0,-1){\halflength}}%
   \or \put(\xpos,\ypos){\line(0,-1)\halflength}%
       \put(\xpos,\ypos){\vector(0,-1)3}%
       \advance\ypos by-\halflength
       \advance\ypos by-\gap
       \put(\xpos,\ypos){\vector(0,-1){\halflength}}%
   \or \put(\xpos,\ypos){\line(0,-1)\halflength}%
       \advance\ypos by-\halflength
       \advance\ypos by-\gap
       \put(\xpos,\ypos){\evector(0,-1){\halflength}}%
   \fi
\else \arrowtype=-\arrowtype
   \ifcase\arrowtype
   \or \advance \ypos by-\arrowlength
       \put(\xpos,\ypos){\line(0,1){\halflength}}%
       \advance\ypos by\halflength
       \advance\ypos by\gap
       \put(\xpos,\ypos){\vector(0,1){\halflength}}%
   \or \advance \ypos by-\arrowlength
       \put(\xpos,\ypos){\line(0,1)\halflength}%
       \put(\xpos,\ypos){\vector(0,1)3}%
       \advance\ypos by\halflength
       \advance\ypos by\gap
       \put(\xpos,\ypos){\vector(0,1){\halflength}}%
   \or \advance \ypos by-\arrowlength
       \put(\xpos,\ypos){\line(0,1)\halflength}%
       \advance\ypos by\halflength
       \advance\ypos by\gap
       \put(\xpos,\ypos){\evector(0,1){\halflength}}%
   \fi
\fi
}

\def\putmorphism(#1)(#2,#3)[#4`#5`#6]#7#8#9{{%
\run #2
\rise #3
\ifnum\rise=0
  \puthmorphism(#1)[#4`#5`#6]{#7}{#8}#9%
\else\ifnum\run=0
  \putvmorphism(#1)[#4`#5`#6]{#7}{#8}#9%
\else
\setpos(#1)%
\arrowlength #7
\arrowtype #8
\ifnum\run=0
\else\ifnum\rise=0
\else
\ifnum\run>0
    \coefa=1
\else
   \coefa=-1
\fi
\ifnum\arrowtype>0
   \coefb=0
   \coefc=-1
\else
   \coefb=\coefa
   \coefc=1
   \arrowtype=-\arrowtype
\fi
\width=2
\multiply \width by\run
\divide \width by\rise
\ifnum \width<0  \width=-\width\fi
\advance\width by60
\if l#9 \width=-\width\fi
\putbox(\xpos,\ypos){#4}
{\multiply \coefa by\arrowlength
\advance\xpos by\coefa
\multiply \coefa by\rise
\divide \coefa by\run
\advance \ypos by\coefa
\putbox(\xpos,\ypos){#5} }%
{\multiply \coefa by\arrowlength
\divide \coefa by2
\advance \xpos by\coefa
\advance \xpos by\width
\multiply \coefa by\rise
\divide \coefa by\run
\advance \ypos by\coefa
\if l#9%
   \putrbox(\xpos,\ypos){#6}%
\else\if r#9%
   \putlbox(\xpos,\ypos){#6}%
\fi\fi }%
{\multiply \rise by-\coefc
\multiply \run by-\coefc
\multiply \coefb by\arrowlength
\advance \xpos by\coefb
\multiply \coefb by\rise
\divide \coefb by\run
\advance \ypos by\coefb
\multiply \coefc by70
\advance \ypos by\coefc
\multiply \coefc by\run
\divide \coefc by\rise
\advance \xpos by\coefc
\multiply \coefa by140
\multiply \coefa by\run
\divide \coefa by\rise
\advance \arrowlength by\coefa
\ifcase\arrowtype
\or \put(\xpos,\ypos){\vector(\run,\rise){\arrowlength}}%
\or \put(\xpos,\ypos){\mvector(\run,\rise){\arrowlength}}%
\or \put(\xpos,\ypos){\evector(\run,\rise){\arrowlength}}%
\fi}\fi\fi\fi\fi}}

\newcount\numbdashes \newcount\lengthdash \newcount\increment

\def\howmanydashes{
\numbdashes=\arrowlength \lengthdash=40
\divide\numbdashes by \lengthdash
\lengthdash=\arrowlength
\divide\lengthdash by \numbdashes
\increment=\lengthdash
\multiply\lengthdash by 3
\divide\lengthdash by 5
}

\def\putdashvector(#1)(#2,#3)#4#5{%
\ifnum#3=0 \putdashhvector(#1){#4}#5
\else
\ifnum#2=0
\putdashvvector(#1){#4}#5\fi\fi}

\def\putdashhvector(#1,#2)#3#4{{%
\arrowlength=#3 \howmanydashes
\multiput(#1,#2)(\increment,0){\numbdashes}%
{\vrule height .4pt width \lengthdash\unitlength}
\arrowtype=#4 \xpos=#1
\ifnum\arrowtype<0 \advance\arrowtype by 7 \fi
\ifcase\arrowtype
\or \advance\xpos by 10
    \put(\xpos,#2){\vector(-1,0){\lengthdash}}
    \advance\xpos by 40
    \put(\xpos,#2){\vector(-1,0){\lengthdash}}
\or \advance \xpos by 10
    \put(\xpos,#2){\vector(-1,0){\lengthdash}}
    \advance\xpos by  \arrowlength
    \advance\xpos by  -50
    \put(\xpos,#2){\vector(-1,0){\lengthdash}}
\or \advance\xpos by 10
    \put(\xpos,#2){\vector(-1,0){\lengthdash}}
\or \advance\xpos by \arrowlength
    \advance\xpos by -\lengthdash
    \put(\xpos,#2){\vector(1,0){\lengthdash}}
\or {\advance\xpos by 10
    \put(\xpos,#2){\vector(1,0){\lengthdash}}}
    \advance\xpos by \arrowlength
    \advance\xpos by -\lengthdash
    \put(\xpos,#2){\vector(1,0){\lengthdash}}
\or \advance\xpos by \arrowlength
    \advance\xpos by -\lengthdash
    \put(\xpos,#2){\vector(1,0){\lengthdash}}
    \advance\xpos by -40
    \put(\xpos,#2){\vector(1,0){\lengthdash}}
   \fi
}}

\def\putdashvvector(#1,#2)#3#4{{%
\arrowlength=#3 \howmanydashes
\ypos=#2 \advance\ypos by -\arrowlength
\multiput(#1,#2)(0,\increment){\numbdashes}%
    {\vrule width .4pt height \lengthdash\unitlength}
\arrowtype=#4 \ypos=#2
\ifnum\arrowtype<0 \advance\arrowtype by 7 \fi
\ifcase\arrowtype
\or \advance\ypos by \arrowlength \advance\ypos by -40
    \put(#1,\ypos){\vector(0,1){\lengthdash}}
    \advance\ypos by -40
    \put(#1,\ypos){\vector(0,1){\lengthdash}}
\or \advance\ypos by 10
    \put(#1,\ypos){\vector(0,1){\lengthdash}}
    \advance\ypos by \arrowlength \advance\ypos by -40
    \put(#1,\ypos){\vector(0,1){\lengthdash}}
\or \advance\ypos by \arrowlength \advance\ypos by -40
    \put(#1,\ypos){\vector(0,1){\lengthdash}}
\or \advance\ypos by 10
    \put(#1,\ypos){\vector(0,-1){\lengthdash}}
\or \advance\ypos by 10
    \put(#1,\ypos){\vector(0,-1){\lengthdash}}
    \advance\ypos by \arrowlength \advance\ypos by -40
    \put(#1,\ypos){\vector(0,-1){\lengthdash}}
\or \advance\ypos by 10
    \put(#1,\ypos){\vector(0,-1){\lengthdash}}
    \advance\ypos by 40
    \put(#1,\ypos){\vector(0,-1){\lengthdash}}
\fi
}}

\def\puthmorphism(#1,#2)[#3`#4`#5]#6#7#8{{%
\xpos #1
\ypos #2
\width #6
\arrowlength #6
\arrowtype=#7
\putbox(\xpos,\ypos){#3\vphantom{#4}}%
{\advance \xpos by\arrowlength
\putbox(\xpos,\ypos){\vphantom{#3}#4}}%
\horsize{\tempcounta}{#3}%
\horsize{\tempcountb}{#4}%
\divide \tempcounta by2
\divide \tempcountb by2
\advance \tempcounta by30
\advance \tempcountb by30
\advance \xpos by\tempcounta
\advance \arrowlength by-\tempcounta
\advance \arrowlength by-\tempcountb
\putvector(\xpos,\ypos)(1,0)\arrowlength\arrowtype
\divide \arrowlength by2
\advance \xpos by\arrowlength
\vertsize{\tempcounta}{#5}%
\divide\tempcounta by2
\advance \tempcounta by20
\if a#8 %
   \advance \ypos by\tempcounta
   \putbox(\xpos,\ypos){#5}%
\else
   \advance \ypos by-\tempcounta
   \putbox(\xpos,\ypos){#5}%
\fi}}

\def\putvmorphism(#1,#2)[#3`#4`#5]#6#7#8{{%
\xpos #1
\ypos #2
\arrowlength #6
\arrowtype #7
\settowidth{\xlen}{$#5$}%
\putbox(\xpos,\ypos){#3}%
{\advance \ypos by-\arrowlength
\putbox(\xpos,\ypos){#4}}%
{\advance\arrowlength by-140
\advance \ypos by-70
\ifdim\xlen>0pt
   \if m#8%
      \putsplitvector(\xpos,\ypos)\arrowlength\arrowtype
   \else
   \putvector(\xpos,\ypos)(0,-1)\arrowlength\arrowtype
   \fi
\else
   \putvector(\xpos,\ypos)(0,-1)\arrowlength\arrowtype
\fi}%
\ifdim\xlen>0pt
   \divide \arrowlength by2
   \advance\ypos by-\arrowlength
   \if l#8%
      \advance \xpos by-40
      \putrbox(\xpos,\ypos){#5}%
   \else\if r#8%
      \advance \xpos by40
      \putlbox(\xpos,\ypos){#5}%
   \else
      \putbox(\xpos,\ypos){#5}%
   \fi\fi
\fi
}}

\def\putsquarep<#1>(#2)[#3;#4`#5`#6`#7]{{%
\setsqparms[#1]%
\setpos(#2)%
\settokens`#3`%
\puthmorphism(\xpos,\ypos)[\tokenc`\tokend`{#7}]{\width}{\arrowtyped}b%
\advance\ypos by \height
\puthmorphism(\xpos,\ypos)[\tokena`\tokenb`{#4}]{\width}{\arrowtypea}a%
\putvmorphism(\xpos,\ypos)[``{#5}]{\height}{\arrowtypeb}l%
\advance\xpos by \width
\putvmorphism(\xpos,\ypos)[``{#6}]{\height}{\arrowtypec}r%
}}

\def\putsquare{\@ifnextchar <{\putsquarep}{\putsquarep%
   <\arrowtypea`\arrowtypeb`\arrowtypec`\arrowtyped;\width`\height>}}
\def\square{\@ifnextchar< {\squarep}{\squarep
   <\arrowtypea`\arrowtypeb`\arrowtypec`\arrowtyped;\width`\height>}}
\def\squarep<#1>[#2`#3`#4`#5;#6`#7`#8`#9]{{
\setsqparms[#1]
\diagram
\putsquarep<\arrowtypea`\arrowtypeb`\arrowtypec`
\arrowtyped;\width`\height>
(0,0)[#2`#3`#4`{#5};#6`#7`#8`{#9}]
\enddiagram
}}                                                 
\def\putptrianglep<#1>(#2,#3)[#4`#5`#6;#7`#8`#9]{{%
\settriparms[#1]%
\xpos=#2 \ypos=#3
\advance\ypos by \height
\puthmorphism(\xpos,\ypos)[#4`#5`{#7}]{\height}{\arrowtypea}a%
\putvmorphism(\xpos,\ypos)[`#6`{#8}]{\height}{\arrowtypeb}l%
\advance\xpos by\height
\putmorphism(\xpos,\ypos)(-1,-1)[``{#9}]{\height}{\arrowtypec}r%
}}

\def\putptriangle{\@ifnextchar <{\putptrianglep}{\putptrianglep
   <\arrowtypea`\arrowtypeb`\arrowtypec;\height>}}
\def\ptriangle{\@ifnextchar <{\ptrianglep}{\ptrianglep
   <\arrowtypea`\arrowtypeb`\arrowtypec;\height>}}
\def\ptrianglep<#1>[#2`#3`#4;#5`#6`#7]{{
\settriparms[#1]
\diagram
\putptrianglep<\arrowtypea`\arrowtypeb`
\arrowtypec;\height>
(0,0)[#2`#3`#4;#5`#6`{#7}]
\enddiagram
}}                                            

\def\putqtrianglep<#1>(#2,#3)[#4`#5`#6;#7`#8`#9]{{%
\settriparms[#1]%
\xpos=#2 \ypos=#3
\advance\ypos by\height
\puthmorphism(\xpos,\ypos)[#4`#5`{#7}]{\height}{\arrowtypea}a%
\putmorphism(\xpos,\ypos)(1,-1)[``{#8}]{\height}{\arrowtypeb}l%
\advance\xpos by\height
\putvmorphism(\xpos,\ypos)[`#6`{#9}]{\height}{\arrowtypec}r%
}}

\def\putqtriangle{\@ifnextchar <{\putqtrianglep}{\putqtrianglep
   <\arrowtypea`\arrowtypeb`\arrowtypec;\height>}}
\def\qtriangle{\@ifnextchar <{\qtrianglep}{\qtrianglep
   <\arrowtypea`\arrowtypeb`\arrowtypec;\height>}}
\def\qtrianglep<#1>[#2`#3`#4;#5`#6`#7]{{
\settriparms[#1]
\width=\height                                
\diagram
\putqtrianglep<\arrowtypea`\arrowtypeb`
\arrowtypec;\height>
(0,0)[#2`#3`#4;#5`#6`{#7}]
\enddiagram
}}

\def\putdtrianglep<#1>(#2,#3)[#4`#5`#6;#7`#8`#9]{{%
\settriparms[#1]%
\xpos=#2 \ypos=#3
\puthmorphism(\xpos,\ypos)[#5`#6`{#9}]{\height}{\arrowtypec}b%
\advance\xpos by \height \advance\ypos by\height
\putmorphism(\xpos,\ypos)(-1,-1)[``{#7}]{\height}{\arrowtypea}l%
\putvmorphism(\xpos,\ypos)[#4``{#8}]{\height}{\arrowtypeb}r%
}}

\def\putdtriangle{\@ifnextchar <{\putdtrianglep}{\putdtrianglep
   <\arrowtypea`\arrowtypeb`\arrowtypec;\height>}}
\def\dtriangle{\@ifnextchar <{\dtrianglep}{\dtrianglep
   <\arrowtypea`\arrowtypeb`\arrowtypec;\height>}}
\def\dtrianglep<#1>[#2`#3`#4;#5`#6`#7]{{
\settriparms[#1]
\width=\height                                
\diagram
\putdtrianglep<\arrowtypea`\arrowtypeb`
\arrowtypec;\height>
(0,0)[#2`#3`#4;#5`#6`{#7}]
\enddiagram
}}

\def\putbtrianglep<#1>(#2,#3)[#4`#5`#6;#7`#8`#9]{{%
\settriparms[#1]%
\xpos=#2 \ypos=#3
\puthmorphism(\xpos,\ypos)[#5`#6`{#9}]{\height}{\arrowtypec}b%
\advance\ypos by\height
\putmorphism(\xpos,\ypos)(1,-1)[``{#8}]{\height}{\arrowtypeb}r%
\putvmorphism(\xpos,\ypos)[#4``{#7}]{\height}{\arrowtypea}l%
}}

\def\putbtriangle{\@ifnextchar <{\putbtrianglep}{\putbtrianglep
   <\arrowtypea`\arrowtypeb`\arrowtypec;\height>}}
\def\btriangle{\@ifnextchar <{\btrianglep}{\btrianglep
   <\arrowtypea`\arrowtypeb`\arrowtypec;\height>}}
\def\btrianglep<#1>[#2`#3`#4;#5`#6`#7]{{
\settriparms[#1]
\width=\height                               
\diagram
\putbtrianglep<\arrowtypea`\arrowtypeb`
\arrowtypec;\height>
(0,0)[#2`#3`#4;#5`#6`{#7}]
\enddiagram
}}

\def\putAtrianglep<#1>(#2,#3)[#4`#5`#6;#7`#8`#9]{{%
\settriparms[#1]%
\xpos=#2 \ypos=#3
{\multiply \height by2
\puthmorphism(\xpos,\ypos)[#5`#6`{#9}]{\height}{\arrowtypec}b}%
\advance\xpos by\height \advance\ypos by\height
\putmorphism(\xpos,\ypos)(-1,-1)[#4``{#7}]{\height}{\arrowtypea}l%
\putmorphism(\xpos,\ypos)(1,-1)[``{#8}]{\height}{\arrowtypeb}r%
}}

\def\putAtriangle{\@ifnextchar <{\putAtrianglep}{\putAtrianglep
   <\arrowtypea`\arrowtypeb`\arrowtypec;\height>}}
\def\Atriangle{\@ifnextchar <{\Atrianglep}{\Atrianglep
   <\arrowtypea`\arrowtypeb`\arrowtypec;\height>}}
\def\Atrianglep<#1>[#2`#3`#4;#5`#6`#7]{{
\settriparms[#1]
\width=\height                                     
\diagram
\putAtrianglep<\arrowtypea`\arrowtypeb`
\arrowtypec;\height>
(0,0)[#2`#3`#4;#5`#6`{#7}]
\enddiagram
}}

\def\putAtrianglepairp<#1>(#2)[#3;#4`#5`#6`#7`#8]{{%
\settripairparms[#1]%
\setpos(#2)%
\settokens`#3`%
\puthmorphism(\xpos,\ypos)[\tokenb`\tokenc`{#7}]{\height}{\arrowtyped}b%
\advance\xpos by\height
\puthmorphism(\xpos,\ypos)[\phantom{\tokenc}`\tokend`{#8}]%
{\height}{\arrowtypee}b%
\advance\ypos by\height
\putmorphism(\xpos,\ypos)(-1,-1)[\tokena``{#4}]{\height}{\arrowtypea}l%
\putvmorphism(\xpos,\ypos)[``{#5}]{\height}{\arrowtypeb}m%
\putmorphism(\xpos,\ypos)(1,-1)[``{#6}]{\height}{\arrowtypec}r%
}}

\def\putAtrianglepair{\@ifnextchar <{\putAtrianglepairp}{\putAtrianglepairp%
   <\arrowtypea`\arrowtypeb`\arrowtypec`\arrowtyped`\arrowtypee;\height>}}
\def\Atrianglepair{\@ifnextchar <{\Atrianglepairp}{\Atrianglepairp%
   <\arrowtypea`\arrowtypeb`\arrowtypec`\arrowtyped`\arrowtypee;\height>}}

\def\Atrianglepairp<#1>[#2;#3`#4`#5`#6`#7]{{
\settripairparms[#1]
\settokens`#2`
\width=\height                                
\diagram
\putAtrianglepairp                            
<\arrowtypea`\arrowtypeb`\arrowtypec`
\arrowtyped`\arrowtypee;\height>
(0,0)[{#2};#3`#4`#5`#6`{#7}]
\enddiagram
}}

\def\putVtrianglep<#1>(#2,#3)[#4`#5`#6;#7`#8`#9]{{%
\settriparms[#1]%
\xpos=#2 \ypos=#3
\advance\ypos by\height
{\multiply\height by2
\puthmorphism(\xpos,\ypos)[#4`#5`{#7}]{\height}{\arrowtypea}a}%
\putmorphism(\xpos,\ypos)(1,-1)[`#6`{#8}]{\height}{\arrowtypeb}l%
\advance\xpos by\height
\advance\xpos by\height
\putmorphism(\xpos,\ypos)(-1,-1)[``{#9}]{\height}{\arrowtypec}r%
}}

\def\putVtriangle{\@ifnextchar <{\putVtrianglep}{\putVtrianglep
   <\arrowtypea`\arrowtypeb`\arrowtypec;\height>}}
\def\Vtriangle{\@ifnextchar <{\Vtrianglep}{\Vtrianglep
   <\arrowtypea`\arrowtypeb`\arrowtypec;\height>}}
\def\Vtrianglep<#1>[#2`#3`#4;#5`#6`#7]{{
\settriparms[#1]
\width=\height                                 
\diagram
\putVtrianglep<\arrowtypea`\arrowtypeb`
\arrowtypec;\height>
(0,0)[#2`#3`#4;#5`#6`{#7}]
\enddiagram
}}

\def\putVtrianglepairp<#1>(#2)[#3;#4`#5`#6`#7`#8]{{
\settripairparms[#1]%
\setpos(#2)%
\settokens`#3`%
\advance\ypos by\height
\putmorphism(\xpos,\ypos)(1,-1)[`\tokend`{#6}]{\height}{\arrowtypec}l%
\puthmorphism(\xpos,\ypos)[\tokena`\tokenb`{#4}]{\height}{\arrowtypea}a%
\advance\xpos by\height
\puthmorphism(\xpos,\ypos)[\phantom{\tokenb}`\tokenc`{#5}]%
{\height}{\arrowtypeb}a%
\putvmorphism(\xpos,\ypos)[``{#7}]{\height}{\arrowtyped}m%
\advance\xpos by\height
\putmorphism(\xpos,\ypos)(-1,-1)[``{#8}]{\height}{\arrowtypee}r%
}}

\def\putVtrianglepair{\@ifnextchar <{\putVtrianglepairp}{\putVtrianglepairp%
    <\arrowtypea`\arrowtypeb`\arrowtypec`\arrowtyped`\arrowtypee;\height>}}
\def\Vtrianglepair{\@ifnextchar <{\Vtrianglepairp}{\Vtrianglepairp%
    <\arrowtypea`\arrowtypeb`\arrowtypec`\arrowtyped`\arrowtypee;\height>}}
\def\Vtrianglepairp<#1>[#2;#3`#4`#5`#6`#7]{{
\settripairparms[#1]
\settokens`#2`
\diagram
\putVtrianglepairp                             
<\arrowtypea`\arrowtypeb`\arrowtypec`
\arrowtyped`\arrowtypee;\height>
(0,0)[{#2};#3`#4`#5`#6`{#7}]
\enddiagram
}}

\def\putCtrianglep<#1>(#2,#3)[#4`#5`#6;#7`#8`#9]{{%
\settriparms[#1]%
\xpos=#2 \ypos=#3
\advance\ypos by\height
\putmorphism(\xpos,\ypos)(1,-1)[``{#9}]{\height}{\arrowtypec}l%
\advance\xpos by\height
\advance\ypos by\height
\putmorphism(\xpos,\ypos)(-1,-1)[#4`#5`{#7}]{\height}{\arrowtypea}l%
{\multiply\height by 2
\putvmorphism(\xpos,\ypos)[`#6`{#8}]{\height}{\arrowtypeb}r}%
}}

\def\putCtriangle{\@ifnextchar <{\putCtrianglep}{\putCtrianglep
    <\arrowtypea`\arrowtypeb`\arrowtypec;\height>}}
\def\Ctriangle{\@ifnextchar <{\Ctrianglep}{\Ctrianglep
    <\arrowtypea`\arrowtypeb`\arrowtypec;\height>}}
\def\Ctrianglep<#1>[#2`#3`#4;#5`#6`#7]{{
\settriparms[#1]
\width=\height                               
\diagram
\putCtrianglep<\arrowtypea`\arrowtypeb`
\arrowtypec;\height>
(0,0)[#2`#3`#4;#5`#6`{#7}]
\enddiagram
}}                                           
\def\putDtrianglep<#1>(#2,#3)[#4`#5`#6;#7`#8`#9]{{%
\settriparms[#1]%
\xpos=#2 \ypos=#3
\advance\xpos by\height \advance\ypos by\height
\putmorphism(\xpos,\ypos)(-1,-1)[``{#9}]{\height}{\arrowtypec}r%
\advance\xpos by-\height \advance\ypos by\height
\putmorphism(\xpos,\ypos)(1,-1)[`#5`{#8}]{\height}{\arrowtypeb}r%
{\multiply\height by 2
\putvmorphism(\xpos,\ypos)[#4`#6`{#7}]{\height}{\arrowtypea}l}%
}}

\def\putDtriangle{\@ifnextchar <{\putDtrianglep}{\putDtrianglep
    <\arrowtypea`\arrowtypeb`\arrowtypec;\height>}}
\def\Dtriangle{\@ifnextchar <{\Dtrianglep}{\Dtrianglep
   <\arrowtypea`\arrowtypeb`\arrowtypec;\height>}}
\def\Dtrianglep<#1>[#2`#3`#4;#5`#6`#7]{{
\settriparms[#1]
\width=\height                              
\diagram
\putDtrianglep<\arrowtypea`\arrowtypeb`
\arrowtypec;\height>
(0,0)[#2`#3`#4;#5`#6`{#7}]
\enddiagram
}}                                          
\def\setrecparms[#1`#2]{\width=#1 \height=#2}%

\def\recursep<#1`#2>[#3;#4`#5`#6`#7`#8]{{\m@th
\width=#1 \height=#2
\settokens`#3`
\settowidth{\tempdimen}{$\tokena$}
\ifdim\tempdimen=0pt
  \savebox{\tempboxa}{\hbox{$\tokenb$}}%
  \savebox{\tempboxb}{\hbox{$\tokend$}}%
  \savebox{\tempboxc}{\hbox{$#6$}}%
\else
  \savebox{\tempboxa}{\hbox{$\hbox{$\tokena$}\times\hbox{$\tokenb$}$}}%
  \savebox{\tempboxb}{\hbox{$\hbox{$\tokena$}\times\hbox{$\tokend$}$}}%
  \savebox{\tempboxc}{\hbox{$\hbox{$\tokena$}\times\hbox{$#6$}$}}%
\fi
\ypos=\height
\divide\ypos by 2
\xpos=\ypos
\advance\xpos by \width
\bfig
\putCtrianglep<-1`1`1;\ypos>(0,0)[`\tokenc`;#5`#6`{#7}]%
\puthmorphism(\ypos,0)[\tokend`\usebox{\tempboxb}`{#8}]{\width}{-1}b%
\puthmorphism(\ypos,\height)[\tokenb`\usebox{\tempboxa}`{#4}]{\width}{-1}a%
\advance\ypos by \width
\putvmorphism(\ypos,\height)[``\usebox{\tempboxc}]{\height}1r%
\efig
}}

\def\recurse{\@ifnextchar <{\recursep}{\recursep<\width`\height>}}

\def\puttwohmorphisms(#1,#2)[#3`#4;#5`#6]#7#8#9{{%
\puthmorphism(#1,#2)[#3`#4`]{#7}0a
\ypos=#2
\advance\ypos by 20
\puthmorphism(#1,\ypos)[\phantom{#3}`\phantom{#4}`#5]{#7}{#8}a
\advance\ypos by -40
\puthmorphism(#1,\ypos)[\phantom{#3}`\phantom{#4}`#6]{#7}{#9}b
}}

\def\puttwovmorphisms(#1,#2)[#3`#4;#5`#6]#7#8#9{{%
\putvmorphism(#1,#2)[#3`#4`]{#7}0a
\xpos=#1
\advance\xpos by -20
\putvmorphism(\xpos,#2)[\phantom{#3}`\phantom{#4}`#5]{#7}{#8}l
\advance\xpos by 40
\putvmorphism(\xpos,#2)[\phantom{#3}`\phantom{#4}`#6]{#7}{#9}r
}}

\def\puthcoequalizer(#1)[#2`#3`#4;#5`#6`#7]#8#9{{%
\setpos(#1)%
\puttwohmorphisms(\xpos,\ypos)[#2`#3;#5`#6]{#8}11%
\advance\xpos by #8
\puthmorphism(\xpos,\ypos)[\phantom{#3}`#4`#7]{#8}1{#9}
}}

\def\putvcoequalizer(#1)[#2`#3`#4;#5`#6`#7]#8#9{{%
\setpos(#1)%
\puttwovmorphisms(\xpos,\ypos)[#2`#3;#5`#6]{#8}11%
\advance\ypos by -#8
\putvmorphism(\xpos,\ypos)[\phantom{#3}`#4`#7]{#8}1{#9}
}}

\def\putthreehmorphisms(#1)[#2`#3;#4`#5`#6]#7(#8)#9{{%
\setpos(#1) \settypes(#8)
\if a#9 %
     \vertsize{\tempcounta}{#5}%
     \vertsize{\tempcountb}{#6}%
     \ifnum \tempcounta<\tempcountb \tempcounta=\tempcountb \fi
\else
     \vertsize{\tempcounta}{#4}%
     \vertsize{\tempcountb}{#5}%
     \ifnum \tempcounta<\tempcountb \tempcounta=\tempcountb \fi
\fi
\advance \tempcounta by 60
\puthmorphism(\xpos,\ypos)[#2`#3`#5]{#7}{\arrowtypeb}{#9}
\advance\ypos by \tempcounta
\puthmorphism(\xpos,\ypos)[\phantom{#2}`\phantom{#3}`#4]{#7}{\arrowtypea}{#9}
\advance\ypos by -\tempcounta \advance\ypos by -\tempcounta
\puthmorphism(\xpos,\ypos)[\phantom{#2}`\phantom{#3}`#6]{#7}{\arrowtypec}{#9}
}}

\def\setarrowtoks[#1`#2`#3`#4`#5`#6]{%
\def\toka{#1}
\def\tokb{#2}
\def\tokc{#3}
\def\tokd{#4}
\def\toke{#5}
\def\tokf{#6}
}
\def\hex{\@ifnextchar <{\hexp}{\hexp<1000`400>}}
\def\hexp<#1`#2>[#3`#4`#5`#6`#7`#8;#9]{%
\setarrowtoks[#9]
\yext=#2 \advance \yext by #2
\xext=#1 \advance\xext by \yext
\bfig
\putCtriangle<-1`0`1;#2>(0,0)[`#5`;\tokb``\tokd]
\xext=#1 \yext=#2 \advance \yext by #2
\putsquare<1`0`0`1;\xext`\yext>(#2,0)[#3`#4`#7`#8;\toka```\tokf]
\advance \xext by #2
\putDtriangle<0`1`-1;#2>(\xext,0)[`#6`;`\tokc`\toke]
\efig
}
\makeatother

\numberwithin{equation}{subsection}
\begin{document}
\maketitle
\def\gk{G_k(V)}

\begin{abstract}
We prove a $C^0$ {\itshape a priori} estimate on a solution of the quaternionic Calabi problem
on an arbitrary compact connected HKT-manifold. This generalizes earlier results \cite{alesker-verbitsky-IJM} and \cite{alesker-shelukhin-13}
where this result was proven under certain extra assumptions on the manifold.
\end{abstract}

\tableofcontents

\section{Introduction.}\label{S:introduction}

\subsection{Preliminaries and statement of the main result.}\label{preliminaries}
The main result of the paper is Theorem \ref{T:main-estimate} below.
In recent years the classical theories of (real) convex and complex plurisubharmonic functions and real and complex Monge-Amp\`ere
equations have been extended in several directions. Thus in \cite{alesker-bsm-03} the first author has introduced a class of plurisubharmonic (psh) functions
of quaternionic variables
on the flat quaternionic space $\HH^n$ and defined a quaternionic version of the Monge-Amp\`ere operator on continuous psh functions on $\HH^n$; applications of these notions
to convexity (valuations theory) were given in \cite{alesker-AIM-05}.
In \cite{alesker-jga-03} the first author studied further the properties of this operator and proved a solvability of the Dirichlet problem for the quaternionic Monge-Amp\`ere equation
in a domain in $\HH^n$ under appropriate assumptions; in a very recent preprint \cite{Zhu-Jingyong} J. Zhu proved the solvability of it under more general assumptions. In \cite{alesker-verbitsky-IJM}
M. Verbitsky and the first author generalized the class of quaternionic psh functions and the quaternionic Monge-Amp\`ere operator to the more general context of hypercomplex manifolds
(the definition is reminded below); there the psh functions were also related to the so called HyperK\"ahler with Torsion (HKT) geometry. Furthermore the first author \cite{alesker-jgp} extended the notion of
psh function and quaternionic Monge-Amp\`ere operator to still more general class of quaternionic manifolds. In \cite{alesker-octonions} the first author introduced psh functions and
the Monge-Amp\`ere operator on functions of two octonionic variables; applications to the valuations theory were also given there.

A parallel development of extending the class of psh functions and Monge-Amp\`ere equations in different, though partly overlapping, directions was initiated by Harvey and Lawson
in a recent series of articles. Thus in \cite{harvey-lawson-calibrated-AJM-09} they introduced and studied psh functions
in calibrated geometry; this approach was generalized further in \cite{harvey-lawson-PSH-general-geom-context-2011}, \cite{harvey-lawson-AIM-2012}.
In \cite{comm-pure-appl-math-09} Harvey and Lawson introduce and study the Dirichlet problem for a rather general class of Monge-Amp\`ere type equations which covers the
Dirichlet problem for the homogeneous
(i.e. with the vanishing right hand side) Monge-Amp\`ere equations over $\RR,\CC,\HH$. In a very recent preprint \cite{harvey-lawson-plis} they and Pli\'s study psh functions
on almost complex manifolds.

\hfill

In order to formulate the main result of this paper we have to remind some notions.

\begin{definition}
A smooth manifold $M^{4n}$ with a triple of complex structures $I,J,K$ such that $IJ=-JI=K$ will be called a \textit{hypercomplex manifold}.
\end{definition}

\begin{remark}
Hypercomplex manifolds were explicitly introduced by Boyer \cite{boyer}. We assume throughout the paper that $I,J,K$ act
on the tangent bundle $TM$ by the right multiplication. This induces the right action of the field $\HH$ of quaternions on $TM$.
\end{remark}

\begin{definition}
A hypercomplex manifold $(M^{4n},I,J,K)$ with a Riemannian
metric $\rho$ that is invariant under the three complex structures $(I,J,K)$ will be called a \textit{hyper-Hermitian} manifold.
\end{definition}

\begin{definition}
Consider the complex manifold $(M,I)$. Denote by $\Lam^{p,q}_I(M)$
the space of $(p,q)$-forms on $(M,I)$. Let $\del:\Lam^{p,q}_I(M) \to
\Lam^{p+1,q}_I(M)$ and $\delbar:\Lam^{p,q}_I(M) \to \Lam^{p,q+1}_I(M)$
be the usual $\del$ and $\delbar$ operators on the complex manifold
$(M,I)$. Define $$\del_J := J^{-1}\circ \delbar \circ J.$$ Then by
\cite{hkt-mflds-supersymmetry} we have
\begin{eqnarray*}
(a)\; J:\Lam^{p,q}_I(M) \to
\Lam^{q,p}_I(M),\\(b)\;
\del_J:\Lam^{p,q}_I(M) \to \Lam^{p+1,q}_I(M),\\
(c)\; \del\del_J=-\del_J\del.
\end{eqnarray*}
\end{definition}

\begin{remark}
The operator $\pt\pt_J$ should be considered as a quaternionic analogue of the operator $dd^c$ on complex manifolds.
Moreover  the operator $\pt\pt_J$
on functions is a quaternionic version of the Hessian.
\end{remark}

Given a hyper-Hermitian manifold $(M,I,J,K,\rho)$ consider the differential form $$\Om:= - \om_J + \sqrt{-1}\om_K$$ where $\om_L(A,B):=\rho(A,B \cdot L)$ for any $L \in \HH$ with $L^2= -1$ and any two vector fields $A,B$ on $M$. We use the following definition of an HKT-metric.

\begin{definition}\label{definition HKT metric} The metric $\rho$ on $M$ is called an \textit{HKT-metric} if $$\del \Om = 0.$$
\end{definition}

\begin{remark} HKT manifolds were introduced in the physical literature by
Howe and Papadopoulos \cite{howe-papa}. For the mathematical treatment see Grantcharov-Poon \cite{_Gra_Poon_} and Verbitsky \cite{hkt-mflds-supersymmetry}. The original definition of HKT-metrics in \cite{howe-papa}
was different but equivalent to Definition \ref{definition HKT metric} (the latter was given in \cite{_Gra_Poon_}).
\end{remark}

\begin{remark}
The classical hyper-K\"ahler metrics (i.e. Riemannian metrics with
the holonomy of the Levi-Civita connection contained in the group $Sp(n)$) form a subclass of
HKT-metrics. It is well known that a hyper-Hermitian metric $\rho$
is hyper-K\"ahler if and only if the form $\Om$ is closed, or
equivalently $\del \Om  = 0$ and $\delbar \Om = 0.$
\end{remark}

\begin{definition} A form $\om \in \Lam^{2k,0}_I(M)$ that satisfies $J \om = \overline{\om}$ will be called a \textit{real (2k,0)-form} \cite{hkt-mflds-supersymmetry}. The space of real $(2k,0)$-forms will be denoted by $\Lam^{2k,0}_{I,\R}(M).$
\end{definition}

\begin{definition}\label{D:t-isomorphism}({\it The $t$-isomorphism} \cite{hkt-mflds-supersymmetry}) Given a right vector space $V$ over $\HH$, there is an $\RR$-linear
isomorphism $t:\Lambda^{2,0}_{I,\R}(V) \to S_\HH(V)$ from the space
of the real $(2,0)$-forms on $V$ to the space of symmetric forms on
$V$ that are invariant with respect to the action of the group
$SU(2)$ of unit quaternions, given by $t(\eta)(A,A):=\eta(A,A\circ
J)$ (\cite{hkt-mflds-supersymmetry}). We shall denote by the same
letter the induced isomorphism $t:\Lam^{2,0}_{I,\R}(M) \to S_\HH(M)$,
where the target is the space of global sections of the bundle with
the fiber $S_\HH(M)_x = S_\HH(T_x M)$ over an arbitrary point $x \in
M$.
\end{definition}

\begin{definition} Forms $\eta \in \Lam^{2,0}_{I,\R}(M)$ satisfying $t(\eta) > 0$ (resp. $t(\eta) \geq 0$) will be called \textit{strictly positive forms}
(resp. \textit{positive forms}).
\end{definition}

\begin{remark}
On a hyper-Hermitian manifold the metric $\rho$ and the form $\Om =  - \om_J + \sqrt{-1}\om_K$ satisfy $t(\Om)=\rho$ and are therefore mutually defining.
\end{remark}

Thus we will denote an HKT-manifold as a 5-tuple $(M,I,J,K,\Omega)$, where $\Omega$ is the form corresponding to the HKT-metric $\rho$ via the $t$-isomorphism; $\Omega$ is called an HKT-form.

\hfill

The following quaternionic version of the Calabi problem was introduced by M. Verbitsky and the first author in \cite{alesker-verbitsky-IJM}; it is the main object of study of the present paper.
The problem formulated in analogy with the famous complex Calabi problem
solved by Yau \cite{yau}.
Let $(M^{4n},I,J,K,\Omega)$ be a compact HKT-manifold. Let $F\in C^\infty(M,\RR)$ be a smooth real
valued function. Then it is conjectured \cite{alesker-verbitsky-IJM} that there exists a constant $A>0$ such that
the equation
\begin{eqnarray}\label{E:CY}
(\Om + \del\del_J \phi)^n = A e^F \Om^n
\end{eqnarray}
has a solution $\phi \in C^\infty(M,\R)$.

This equation is
non-linear (for $n>1$) and elliptic of second order. The uniqueness of a
solution (up to an additive constant) is shown in \cite{alesker-verbitsky-IJM}.
M. Verbitsky \cite{verbitsky-2009} has found a geometric
interpretation of this equation. The quaternionic Calabi problem (\ref{E:CY}) was shown to have a positive solution by the first author \cite{alesker-AIM-13}
under rather restrictive assumption: it was assumed that the hypercomplex manifold $(M,I,J,K)$ admits a locally flat hyperK\"ahler metric (which does not have to be related to $\Omega$).

The main result of this paper is the following theorem which contributes to understanding the general case.
\begin{theorem}\label{T:main-estimate}
Let $(M^{4n},I,J,K,\Omega)$ be a compact connected HKT-manifold, and let $f\in C^\infty(M,\RR)$ is a smooth real
valued function.
Let $\phi\in C^\infty(M,\RR)$ be a solution of (\ref{E:CY}) with the constant $A=1$ normalized such that $\max_{M}\phi=0$.
\footnote{The assumption $A=1$ is made for convenience only. Clearly $A$ can be absorbed by redefining $F$ to be $F+\log A$.}
Then there exists a constant $C>0$ depending
on $M,I,J,K,\Omega$ and $||F||_{L^\infty(M)}$ such that
$$||\phi||_{L^\infty(M)}\leq C.$$
\end{theorem}

\begin{remark}
This result, i.e. the uniform estimate on the solution, was previously proved in \cite{alesker-verbitsky-IJM} under the extra assumption that the holonomy group of the Obata connection
(defined in Theorem \ref{T:obata} below) of the hypercomplex manifold $(M,I,J,K)$
is contained in the group $SL_n(\HH)$ (while in general it is contained in the larger group $GL_n(\HH)$ only). The same uniform estimate was proved in the paper \cite{alesker-shelukhin-13} under a different assumption that the hypercomplex structure $(M,I,J,K)$ is locally flat. Thus Theorem \ref{T:main-estimate} contains the two previous estimates as special cases.
Moreover it is strictly stronger than the union of them: there exist compact HKT-manifolds such that the hypercomplex structure is not locally flat and the holonomy of the Obata
connection is not contained in $SL_n(\HH)$. Indeed Joyce \cite{joyce-1992} has constructed on many compact Lie groups, including the 8-dimensional group $SU(3)$,
a left-invariant hypercomplex structure. Each of these hypercomplex manifolds admits an HKT-metric by the work \cite{_Gra_Poon_} by
Grantcharov and Poon.  Recently Soldantenkov has shown that the holonomy of the Obata connection of the Joyce hypercomplex structure on $SU(3)$ is
equal to $GL_2(\HH)$ thus providing a required example.
In fact it is expected that many of the other compact Lie groups with the Joyce's hypercomplex structure have also maximal holonomy of
the Obata connection thus providing more examples of manifolds covered by Theorem \ref{T:main-estimate}, but covered neither by \cite{alesker-verbitsky-IJM} nor by \cite{alesker-shelukhin-13}.
\end{remark}




\subsection{Remarks on the methods of the paper.}\label{Ss:methods}
The method of the paper \cite{alesker-verbitsky-IJM} to obtain the uniform estimate in the case of HKT-manifolds with $SL_n(\HH)$-holonomy of the Obata connection
was a quaternionic modification of the original Yau's method in the complex case
as presented in \cite{joyce-book}. The method of \cite{alesker-shelukhin-13} to obtain the uniform estimate in the case of locally flat hypercomplex structure was heavily based on
the Blocki's method \cite{blocki}  of the proof
of the uniform estimate in the {\itshape complex} Calabi-Yau theorem.
The proof of Theorem \ref{T:main-estimate} in the current paper is heavily based again on the Blocki's method \cite{blocki}. However there is a new subtlety: at some point
we have to show, in particular, that given a 2-jet $\bar f$ of a function at a point $z$ such that $\pt\pt_J \bar f(z)=0$ then there exists a smooth function $f$ at a neighborhood of $z$ such that
$\pt\pt_Jf\equiv 0$ and its 2-jet at $z$ is equal to $\bar f$. A slightly more precise statement is formulated in Theorem \ref{Penrose key proposition} below. The complex analogue of this result
is trivial; also this result is trivial in the case of locally flat hypercomplex manifold. The proof of Theorem \ref{Penrose key proposition} occupies the main physical body of the paper and uses
a cohomological interpretation of the kernel of $\pt\pt_J$ and an application of the Andreotti-Grauert theory.
The former is obtained by identifying  $\pt\pt_J$ with one of the so called Baston operators which, in turn, are constructed using the notion of the Penrose transform.

\hfill

In principal yet another approach might be possible to the study of quaternionic Calabi problem. In \cite{alesker-verbitsky-IJM} it was interpreted
as a special case of the complex Hessian equation. The latter equation has been studied rather actively in the literature, see e.g. \cite{blocki-hessian},
\cite{dinew-kolodziej}).

\hfill

{\bf Acknowledgements.} We are very grateful to V. Palamodov and M. Verbitsky  for answering our numerous questions. The first author is much obliged to G. Henkin who has explained him
the relation between the quaternionic Hessian and the Penrose transform in the flat case many years ago.  We thank Z. Blocki for very useful discussions.

\section{Some background on quaternionic linear algebra.}\label{S:background}\setcounter{subsection}{1}
For a quadratic form $Q$ on a right vector space $V$ over the quaternions we shall denote by $\av{Q}_u$ its average over the action of the group $SU(2)$ of unit quaternions, that is $$\av{Q}_u(x):=\int_{SU(2)} Q(x \circ L) d\nu(L),$$ for the Haar probability measure $\nu$ on $SU(2)$. We shall denote by $Q_+$ the projection $Q_+:=Q - \av{Q}_u$ of $Q$.

\begin{lemma}\label{average over SU(2) of a quadratic form}
For a quadratic form $Q$ on a right vector space $V$ over the fields $\HH$ of quaternions $$\av{Q}_u(x) = \frac{1}{4}(Q(x) + Q(x \circ I) + Q(x \circ J) + Q(x \circ K)).$$
\end{lemma}

The lemma follows by noting that both sides of the equality are $SU(2)$-invariant on one hand, and both their averages over $SU(2)$ equal $\av{Q}_u(x)$ on the other.

\begin{definition} For a quaternion $q\in \HH$ written in the
standard form
\[q= t + x \cdot i + y \cdot j + z \cdot k\] define the Cauchy-Riemann-Fueter-Dirac operator $\frac{\del}{\del \overline{q}}$
on an $\HH$-valued function $F$ by
\[\frac{\del}{\del \overline{q}}F = \frac{\del}{\del t}F + \frac{\del}{\del x}F + \frac{\del}{\del y}F + \frac{\del}{\del z}F,\] and its quaternionic conjugate $\frac{\del}{\del q}$ by \[\frac{\del}{\del q}F = \overline{\frac{\del}{\del \overline{q}}\overline{F}} = \frac{\del}{\del t}F + \frac{\del}{\del x}F + \frac{\del}{\del y}F + \frac{\del}{\del z}F.\]
\end{definition}

Note that the quaternionic Hessian $Hess_\HH (u) := (\frac{\del^2
u}{\del\bar q_i \del q_j})$ is the average of the quadratic
form $D^2 u$ over $SU(2)$, as a short computation shows. As for
different $i$ and $j$ the operators $\frac{\del}{\del q_i}$ and
$\frac{\del}{\del \overline{q_j}}$ commute, this matrix is
hyper-Hermitian, namely satisfies the following definition.

\begin{definition} We call a quaternionic matrix $A = (a_{ij})_{i,j = 1}^n$ \textit{hyper-Hermitian} if $a_{ji} = \overline{a_{ij}}$ for the quaternionic conjugation $\HH \ni q \mapsto \overline{q} \in \HH.$
\end{definition}

We shall also use a version of a determinant defined for
hyper-Hermitian quaternionic matrices referring for further details,
properties, and references to \cite{alesker-bsm-03}. Considering $\HH$
as an $\R$-linear vector space of dimension $4$ we have the
embedding $Mat(n,\HH) \inj Mat(4n,\R)$ of matrix $\R$-algebras. For $A\in Mat(n,\HH)$ we denote by
$^\R A$ its image in $ Mat(4n,\RR)$.
Denote by $\mathcal{H}_n$ the image of the subspace of
hyper-Hermitian matrices. In the following proposition we denote by $\det\colon Mat(4n,\RR)\to \RR$ the usual determinant of real matrices.

\begin{proposition}
There exists a unique polynomial $P\colon \mathcal{H}_n\to \RR$ such that $\det(^\R A)=P^4( A)$ for any $A\in \mathcal{H}_n$ and $P(Id)=1$.
Moreover $P$ is homogenous of degree $n$.
\end{proposition}

\begin{definition}
For a hyper-Hermitian matrix $A$ its \textit{Moore determinant} is $$\det(A):= P ( A) \in \R.$$
\end{definition}

\begin{theorem}
The Moore determinant restricts to the usual determinant on the $\R$-subspace of complex Hermitian matrices.
\end{theorem}

\begin{remark}\label{moore-non-linear}
On the flat HKT manifold $(\HH^n,I,J,K)$, the Moore determinant of $Hess_\HH(u)$ can be
naturally identified (up to a positive multiplicative
constant) with $(\del\del_J u)^n$. More precisely \cite[Corollary
4.6]{alesker-verbitsky-JGA-06} says
$$(\pt\pt_J f)^n=\frac{n!}{4^n}\det\left(\frac{\pt^2f}{\pt\bar q_i\pt q_j}\right)dz_1\wedge dz_2\wedge\dots\wedge dz_{2n-1}\wedge dz_{2n}.$$
\end{remark}

\section{Proofs of the  main result and of some of the auxiliary results.}\label{S:main}

\subsection{A technical proposition on quaternionic Hessian.}\label{Ss:tech-prop-Hessian}
Let $(M^{4n},I,J,K)$ be a hypercomplex manifold. Let $z\in M$. Let us choose an open neighborhood $D$ of $z$ and a biholomoprhic diffeomorphism of
$(D,I)$ with an open subset of $\CC^{2n}$.
\begin{proposition}\label{flattening nonlinear Hessian}
On an $I$-complex coordinate chart $D \subset M$ of a hypercomplex manifold $(M^{4n},I,J,K)$ as above with complex coordinates $z_1,...,z_{2n},$
the $(2,0)$-form $\del\del_J u(z)$ for a function $u \in \sm{M,\R}$ at a point $z \in D$ depends only on the second derivatives $\del_{ij} u$ of $u$ and on the complex structure $J(z)$ at the point $z$ (and not on the derivatives of $J$ at $z$).
\end{proposition}

{\bf Proof.}
Indeed computing in local coordinates we have \[\displaystyle\del\del_J u = \del \circ J^{-1} \circ \delbar u = \del \circ J^{-1} (\Sum_j \del_{\overline{z}_j}(u)\; d{\overline{z}_j}) = \del (\Sum_{j,k} \del_{\overline{z}_j}(u)\; (J^{-1})^j_k \; d{z_k}) =\]
\[= \del^2_{\overline{z}_j{z_l}}(u) \; (J^{-1})^j_k \; dz_l \wedge dz_k + \del_{\overline{z}_j}(u) \; \frac{\del(J^{-1})^j_k}{\del{z_l}} \; dz_l \wedge dz_k.\]

Therefore it is enough to show that the second summand vanishes. To this end we compute $\del\del_J u$ for $u$ to be a linear combination of $z_k$ and $\overline{z}_j.$ We obtain \[\del\del_J z_k = \del \circ J^{-1} (\delbar z_k) = 0\] and \[\del\del_J \overline{z}_j = -\del_J (\del \overline{z}_j) = 0.\] The proof is concluded by plugging these results into the previous computation of $\del\del_J u.$
\qed

\hfill

For $z$ and $D$ as above let us denote by $\tilde J_z$ the complex structure on $D$ which is uniquely characterized by the two properties:

$\bullet$
under the above chosen diffeomorphism of $D$ with an open subset of $\CC^{2n}$ the complex structure $\tilde J_z$ is translation invariant;

$\bullet$ at the point $z$ the complex structure $\tilde J_z$ is equal to the complex structure $J$.

It is easy to see that $I$ and $\tilde J_z$ anti-commute, and thus $(D,I,\tilde J_z, I\cdot \tilde J_z)$ is a hypercomplex manifold. We have the following immediate corollary.

\begin{corollary}\label{Cor:q-hessians-comparison}
Under the above notation for any smooth function $u\colon D\to \RR$ one has
$$\pt\pt_Ju(z)=\pt\pt_{\tilde J_z}u(z).$$
\end{corollary}

\subsection{Proof of Theorem \ref{T:main-estimate}.}
Here we present the proof of Theorem \ref{T:main-estimate} assuming several important propositions
that we prove later. The method of the proof is heavily based on the paper \cite{blocki} of Z. Blocki.
We use the generic notation $const$ to denote any positive constant that depends only on $n,M,I,J,K,\Ome$.


{\it Step 1: the $L^1$ bound.}

It was shown in \cite[Proposition 2.3]{alesker-shelukhin-13} that
if $\max_M\phi=0$ then the norm $||\phi||_{L^1}$ is bounded by a constant depending on $M,I,J,K,\Ome$ only.


{\it Step 2: The Taylor expansion of the potential of the metric.}

We first formulate and prove the exact statements that allow us to make our calculations uniform in an appropriate sense, and then make the calculations.
Below we denote by $B(z,r)$ the open ball of radius $r$ centered at point $z\in M$ with respect to the HKT-metric.

\begin{lemma}\label{good r, a and g}
Given an HKT-manifold $(M,I,J,K,\Ome)$ there exist $r>0$ and $a>0$ such that $\forall z_0 \in M$ $\exists g \in \sm{B(z_0,2r)}$ such that
\begin{enumerate}
\item $g < 0$ on $B(z_0,2r)$,
\item $\Ome = \pt\pt_J g$ on $B(z_0,2r)$,
\item $\inf_{B(z_0,2r) \setminus B(z_0,r)}g \geq \inf_{B(z_0,r)}g + a$,
\item $\inf_{B(z_0,r)}g = g(z_0),$
\item $\normC{g}{10}{B(z_0,2r)} \leq const.$
\end{enumerate}
\end{lemma}

{\bf Proof of Lemma \ref{good r, a and g}.}
First we state the following result which will be needed soon.
\begin{theorem}\label{Penrose key proposition}
Let $(M^{4n},I,J,K)$ be a hypercomplex manifold. Let $\rho$ be an arbitrary infinitely smooth Riemannian metric on $M$. Let $\mathcal{C}\subset M$ be a compact subset.
Then there exists $\delta_0>0$ such that for any $z_0\in \mathcal{C}$ and for any 2-jet $\bar f$ of a function at $z_0$ satisfying $\pt\pt_J\bar f(z_0)=0$,
there exists a function $f\in C^\infty(B(z_0,\delta_0),\RR)$ such that
$$\pt\pt_J f=0 \mbox{ on } B(z_0,\delta_0),$$
and the 2-jet of $f$ at $z_0$ is equal to $\bar f$, where $B(z_0,\delta_0)$ is the open ball with respect to the metric $\rho$.
\end{theorem}

We postpone the proof of this theorem till Section \ref{S:proof-of-proposition}. We have the following claim.
\begin{claim}\label{Claim-inside} There exists a finite open cover $\{U_i\}_{i \in I}$ of $M$ such that
\begin{enumerate}
\item $U_i \Subset V_i$ for open sets $V_i$ isomorphic as $I$-complex manifolds to open subsets in $\mathbb{C}^{2n}$,
\item there exist functions $\tilde g_i \in \sm{V_i,\R}$ such that $\Ome = \pt\pt_J \tilde g_i$ on $V_i$,
\item $\normC{\tilde g_i}{20}{U_i} < const,$
\item there exists $r_0>0$ such that for any point $z_0\in M$ there
exists $i_0\in I$ satisfying $B(z_0,3r_0)\Subset U_{i_0}$,
and for any $z_0\in M$, for any 2-jet $\bar f$ of a function at $z_0$ satisfying $\pt\pt_J\bar f(z_0)=0$, there exists
$f\in C^\infty(B(z_0,3r_0),\RR)$ such that
$$\pt\pt_Jf=0 \mbox{ on } B(z_0,3r_0)$$
and the 2-jet of $f$ at $z_0$ is equal to $\bar f$.
\end{enumerate}
\end{claim}
Claim \ref{Claim-inside} is an easy consequence of
\cite[Proposition 1.15]{alesker-verbitsky-JGA-06} and Theorem \ref{Penrose key proposition}.
By Corollary \ref{Cor:q-hessians-comparison} one has $\pt\pt_J \tilde g_i (z)=\pt \pt_{\widetilde{J}_z} \tilde g_i (z)$, where the complex structure
$\tilde J_z$ was defined right before the statement of Corollary \ref{Cor:q-hessians-comparison}. Under the $t$-isomorphism (see Definition \ref{D:t-isomorphism}) at the point $z$ this expression goes to $Hess_\HH(\tilde g_i)_z$ (as a quadratic form) with respect to the flat quaternionic structure given by $I$ and $\widetilde{J}_z.$ Hence we shall write $Hess_\HH(\tilde g_i)_z$ for $t^{-1}(\pt \pt_J \tilde g_i (z)).$

There exists a positive constant $\epsilon$ such that $Hess_\HH(\tilde g_i) > \ep \rho$ on $U_i$ for all $i \in I$, where $\rho$ is the HKT-metric on $M$.

Indeed, since $\tilde g_i$ is strictly plurisubharmonic on $V_i$ we have
$Hess_\HH(\tilde g_i) > 0$ and as $U_i \Subset V_i$, there exist constants
$\ep_i$ such that $Hess_\HH(\tilde g_i) > \ep_i \rho$ on $U_i$. Take $\ep :=
min_{i \in I} \ep_i$.

Given a point $z_0 \in M$, we shall now construct the function $g$ and
choose $r > 0$, $a > 0$. Consider $(U,\widetilde{g})=(U_{i_0},\tilde g_{i_0})$ with $B(z_0,3r_0)\Subset U_{i_0}$.

The Taylor expansion of $\widetilde{g}$ about $z_0$ (in local coordinates on $V_i \Supset U_i$) gives $$\widetilde{g}(z_0 + h) = \widetilde{g}(z_0) + d_{z_0}\widetilde{g}(h) + D^2_{z_0}\widetilde{g}(h) + \Theta(h),$$ where $\Theta(h)=o(|h|^2).$ Now we split the quadratic form $D^2_{z_0}\widetilde{g}(h)$ into its invariant and complementary parts with respect to the action of $SU(2)$: $$D^2_{z_0}\widetilde{g}(h) = \av{D^2_{z_0}\widetilde{g}}_u(h) + (D^2_{z_0}\widetilde{g})_+(h) = Hess_\HH(\widetilde{g})_{z_0}(h) + (D^2_{z_0}\widetilde{g})_+(h).$$

We will now define a function $g$ such that $\pt\pt_J g \equiv \pt\pt_J \widetilde{g}$ on $U$
with the additional property that $D^2_{z_0} g = Hess_\HH(g)_{z_0}.$
To this end let us apply Theorem \ref{Penrose key proposition} to the $2$-jet \[\overline{f}(h) = \widetilde{g}(z_0) + d_{z_0}\widetilde{g}(h) + (D^2_{z_0}\widetilde{g})_+(h),\] that satisfies $\pt \pt_J \overline{f}(z_0) = 0$ by definition.

Let $PH(U):=\{f\in C^\infty(U)|\, \pt\pt_Jf=0\}$ denote the space of real valued pluriharmonic functions in $U$.
Let $j^2_{z_0}$ denote the (finite dimensional) space of 2-jets of functions at $z_0$.
The map
$$PH(U)\to \{\bar f\in j^2_{z_0}|\, \pt\pt_J\bar f(z_0)=0\}$$
which sends $f$ to its 2-jet at $z_0$ is continuous in $C^{10}$-topology on $PH(U)$ and is onto by Theorem \ref{Penrose key proposition}. The Banach inverse theorem implies that
for any $\bar f\in j^2_{z_0}$ such that $\pt\pt_J\bar f(z_0)$ one  can choose $f\in PH(U)$ such that
\begin{eqnarray}\label{E:estimateC10}
||f||_{C^{10}(B(z_0,2r))}\leq const\cdot ||\bar f||,
\end{eqnarray} where $||\cdot ||$
is any continuous norm on $j^2_{z_0}$.

Then we are going to show that with the choice of $f$ satisfying (\ref{E:estimateC10}) the function \[g:=\widetilde{g} - f\] satisfies the required properties. The function $g$ satisfies Condition 5 by (\ref{E:estimateC10}) and property 3 of Claim \ref{Claim-inside}.

Note that by construction $D^2_{z_0}g = Hess_\HH(g)_{z_0} =
Hess_\HH(\widetilde{g})_{z_0} > \ep \rho$ where $\ep$ is as above. We shall choose $r < r_0$ small enough so that $D^2_{(z_0 +
h)}g \geq \frac{\ep}{2} \rho$ for all $h \leq 2r$. Consider the Taylor
expansion of $D^2_{(z_0 + h)}g$ in $h$ about $h=0$: $$D^2_{(z_0 +
h)}g = D^2_{(z_0)}g + \Theta_1(h).$$ The function $\Theta_1$
satisfies $|\Theta_1| < \kappa |h| \rho$ for a constant $\kappa$
depending only on $\normC{\widetilde{g}}{10}{U}$. Therefore, by
Claim \ref{Claim-inside}, $\kappa$ depends on $M,I,J,K,\Ome$ only. Consequently
$$D^2_{(z_0 + h)}g \geq \ep \rho - \kappa |h| \rho = (\ep - \kappa |h|)
\rho.$$ Hence for $|h| < \frac{\ep}{2\kappa}$ we have $D^2_{(z_0 + h)}g
\geq \frac{\ep}{2} \rho$. Therefore we can choose $r =
\min\{r_0,\frac{\ep}{2\kappa}\}/4$. By the
choice of $r$, the function $g$ is convex (in $I$-complex coordinates as above) with $D^2 g \geq
\frac{\ep}{2} \rho$ in $B(z_0,2r)$ with minimum in $z_0$, thus
satisfying condition 4.

By a straightforward computation we have then the estimate $g(z_0 +
h) > const \cdot \frac{\ep h^2}{2}$. Hence $a = const \cdot
\frac{\ep r^2}{2} = const$ can be chosen to satisfy condition 3.

The function $g$ can be
modified by adding a constant to satisfy condition 1.
Lemma \ref{good r, a and g} is proved.

{\it Step 3: End of proof of the main Theorem \ref{T:main-estimate}.}

Choose $z_0 \in M$ at which the function $\phi$ attains its minimum:
$\phi(z_0)= \min_M \phi$. Lemma \ref{good r, a and g} provides us
then with an appropriate $r,a >0$ and $g \in \sm{B(z_0,2r)}$. Then
the function $u = \phi + g$ in the domain $D = B(z_0,2r)$ satisfies
the conditions of Lemma \ref{key lemma} below. Indeed $\pt\pt_J u
= \Ome + \pt\pt_J \phi > 0$ whence $u$ is plurisubharmonic and $u$
is negative on $D$ as both $\phi$ and $g$ are. The set $\{u < \inf_D
u + a\}$ is contained in $B(z_0,r) \Subset B(z_0,2r)$ and hence is
relatively compact.

Lemma \ref{key lemma} gives us $$\normL{\phi}{\infty}{M} -
\normL{g}{\infty}{D} \leq \normL{u}{\infty}{D} \leq a + const \cdot
(r/a)^{4n} \normL{u}{1}{D} \normL{f}{\infty}{D}^4.$$  Moreover
$\normL{u}{1}{D} \leq \normL{\phi}{1}{M} + \normL{g}{\infty}{D} \leq
const$ by Step 1 and property 5 of $g$. Hence $$\normL{\phi}{\infty}{M} \leq const +
const \cdot \normL{f}{\infty}{D}^4,$$ which finishes the proof.
\qed

\subsection{The Alexandrov-Bakelman-Pucci inequality.}\label{Ss:ABP}

We are going to prove the following lemma which was used above.
\begin{lemma}\label{key lemma}
Let $D \subset (M^{4n},I,J,K)$ be an $I$-complex coordinate chart in an HKT manifold, with complex
coordinates $z_1,...,z_{2n}.$ Let $u$ be a negative $C^2$ strictly plurisubharmonic function in $D$ and $a > 0$ a
constant such that the set $\{u < \inf_D u + a\}$ is relatively
compact in $D$. Denote by $f$ the nonlinear Moore
determinant $f = (\pt \pt_J u)^n/ \varkappa_n\bigwedge_{j=1}^{2n} dz_j$ for the constant $\varkappa_n = n!/4^n$.
Then $$\normL{u}{\infty}{D} \leq a +
const \cdot \frac{\mathrm{diam}(D)^{4n}}{a^{4n}} \normL{u}{1}{D}
\normL{f}{\infty}{D}^4.$$
\end{lemma}

As we will show below, this lemma is a formal consequence of the following proposition:

\begin{proposition}\label{key proposition}
Let $D \subset (M^{4n},I,J,K)$ be abounded open subset of an $I$-complex coordinate chart in an HKT manifold, with complex
coordinates $z_1,...,z_{2n}.$ Let $u\in C^2(D)\cap C(\overline{D})$ be a non-positive strictly plurisubharmonic
function in $D$, that vanishes on $\pt D$. Denote by $f$ the nonlinear Moore
determinant $f = (\pt \pt_J u)^n/ \varkappa_n\bigwedge_{j=1}^{2n} dz_j$ for the constant $\varkappa_n = n!/4^n$. Then
$$\normL{u}{\infty}{D} \leq const \cdot \mathrm{diam}(D) \normL{f}{4}{D}^{1/n}.$$
\end{proposition}

\begin{remark}
In both Lemma \ref{key lemma} and Proposition \ref{key proposition} the diameter $\mathrm{diam}(-)$ is measured in the Euclidean metric induced by the real coordinates underlying the $z_j.$
\end{remark}

We will first show how the lemma follows from the proposition and then prove the latter.

\begin{proof}(Lemma \ref{key lemma})
Define $v:= u -\inf_D u -a$. Set $D'=\{v < 0\}$. Then $D'\Subset D$
by assumption. By Proposition \ref{key proposition},
$$a = \normL{v}{\infty}{D'} \leq const \cdot \mathrm{diam}(D') \normL{f}{4}{D'}^{1/n} \leq const \cdot \mathrm{diam}(D')(Vol(D'))^{1/4n} \normL{f}{\infty}{D'}^{1/n}.$$
On the other hand $$Vol(D') \leq \frac{\normL{u}{1}{D'}}{|\inf_D u + a|} = \frac{\normL{u}{1}{D'}}{\normL{u}{\infty}{D} - a}\leq \frac{\normL{u}{1}{D}}{\normL{u}{\infty}{D} - a}.$$ The lemma now follows by
direct substitution.
\end{proof}

{\bf Proof of Proposition \ref{key proposition}.} Consider the $4n$ real coordinates underlying the complex coordinates $z_1,...,z_{2n}.$ By the Alexandrov-Bakelman-Pucci inequality (Lemma 9.2 in \cite{GilbargTrudinger}) we have
\begin{equation}\label{Alexandrov-Bakelman-Pucci}\normL{u}{\infty}{D} \leq const \cdot \mathrm{diam}(D) (\int_\Gamma det D^2 u)^{1/4n},\end{equation}
where $$\Gamma := \{x \in D |\; u(x) + \langle Du(x),y-x \rangle
\leq u(y) \;\forall y\in D\} \subset \{D^2 u \geq 0\},$$
$D^2u$ denotes the usual real Hessian of $u$ on the flat space $\CC^{2n}\simeq \RR^{4n}$ which is a real symmetric $4n\times 4n$ matrix, and $\det D^2u$ denotes its usual determinant.



We claim that at any point $z \in D$ where $D^2 u (z) \geq 0$ we have the following inequality.
\begin{equation}\label{det D^2 and nonlinear Moore} det D^2 u(z)
\leq const \cdot f(z)^4.\end{equation}

Plugging Equation (\ref{det D^2 and nonlinear Moore}) into
Equation (\ref{Alexandrov-Bakelman-Pucci}) yields the proposition.

We now prove inequality (\ref{det D^2 and nonlinear Moore}). By Corollary \ref{Cor:q-hessians-comparison} we can replace $J|_D$ by the flat complex structure $\widetilde{J}_z$ induced by $J(z)$ and by the trivialization of $TD$ given by the complex coordinates (see Section \ref{Ss:tech-prop-Hessian}), with the effect that

\begin{equation}\label{flattened nonlinear Hessian}\pt\pt_J u (z) = \pt\pt_{\widetilde{J}_z} u (z).\end{equation}

Since $I$ and $\widetilde{J}_z$ induce a flat hypercomplex structure on $D,$ we have the following inequality established
in \cite{alesker-shelukhin-13} by use of standard quaternionic linear algebra (we use Remark \ref{moore-non-linear} to cast it in the terms of the present paper).

\begin{equation}\label{det D^2 and linear Moore} det D^2 u(z)
\leq const \cdot (\pt\pt_{\widetilde{J}_z} u (z)^n/\varkappa_n \bigwedge_{j=1}^{2n} dz_j)^4.\end{equation}

By Equation (\ref{flattened nonlinear Hessian}) the right hand side equals $f(z)^4,$ and inequality (\ref{det D^2 and nonlinear Moore}) is proven.
\qed

\section{Proof of Theorem \ref{Penrose key proposition}.}\label{S:proof-of-proposition}
The goal of this section is to prove Theorem \ref{Penrose key proposition}. We will need however some preparations and review background
about twistor space of a hypercomplex manifold, Douady spaces, Penrose transform and Baston operators.
Roughly put, the goal is to show that the operator $\pt\pt_J$ coincides (up to a normalizing constant) with one of the Baston differential operators
arising from the Penrose transform. The kernel of this operator has a cohomological interpretation which will allow us to use tools from
the sheaf cohomology theory of complex analytic manifolds.

\subsection{Twistor space of a hypercomplex manifold.}\label{Ss:twistor}
In \cite{salamon} S. Salamon has attached to any so called quaternionic manifold a complex analytic
manifold called its twistor space. In the special case of locally flat quaternionic manifold the twistor space was previously constructed by Gindikin and Henkin \cite{gindikin-henkin}.
We will review this construction in the special case of hypercomplex manifolds. We refer to Ch. 14G of the book \cite{besse} for details in the general case.

\hfill

Let $(M^{4n},I,J,K)$ be a hypercomplex manifold. For a point $a=(a_1,a_2,a_3)$ from the standard 2-dimensional unit Euclidean sphere $S^2$,
i.e. $a_1^2+a_2^2+a_3^2=1$, we denote by
$$L_a=a_1I+a_2J+a_3K$$
the new almost complex structure on $M$; it turns out to be integrable.

The twistor space $Z$ of $M$ as a smooth manifold is defined by
$$Z:=S^2\times M.$$
Let us define the almost complex structure $\cj$ on $Z$ as follows. For any point
$(a,x)\in Z=S^2\times M$ its tangent space is
$$T_{(a,x)}Z=T_aS^2\oplus T_xM.$$
Define $\cj$ at the point $(a,x)$ to preserve this direct sum, and let it act on $T_aS^2$ as the standard complex structure on the
sphere (which is uniquely defined by the standard metric and orientation), and let it act on $T_xM$ as $L_a$. It is clear that $\cj$ is an almost complex structure. The main property of $\cj$ is the following.
\begin{claim}
The almost complex structure $\cj$ is integrable.
\end{claim}
The following claim is well known and can be easily checked.
\begin{claim}
The natural projection $$Z=S^2\times M\to S^2$$
is a holomorphic map,
where the unit sphere $S^2$ is equipped with the standard complex structure.
\end{claim}
\begin{remark}
Observe that the fibers of the natural projection $Z\to M$ are complex curves $S^2$ which are complex submanifolds of $Z$
and which we will often identify with the complex projective line $\CC\PP^1$ in the standard way. Note that the projection $Z\to M$ is
not a holomorphic map in any reasonable sense.
\end{remark}

\begin{theorem}[Obata \cite{obata}]\label{T:obata}
Let $(M,I,J,K)$ be a hypercomplex manifold. Then there exists a unique torsion free connection $\nabla^{Ob}$ on the tangent bundle
of $M$ which preserves all the complex structures:
$$\nabla^{Ob}I=\nabla^{Ob}J=\nabla^{Ob}K.$$
\end{theorem}

\begin{corollary}\label{Cor:Obata-cor}
Let $(M^{4v},I,J,K)$ be a hypercomplex manifold. Let $\rho$ be an arbitrary infinitely smooth Riemannian metric on $M$. Let $C\subset M$ be a compact subset. Then there
exists $\eps>0$ such that for any point $x_0\in C$ there exists a real analytic diffeomorphism $\phi$ of the open ball $B(x_0,\eps)$ with respect to $\rho$
onto the $\eps$-neighborhood of 0 in $\HH^n$ (with respect to the standard Euclidean metric) which maps $x_0$ to 0 and such that
$$I(x)=I_0+O(x^2), J(x)=J_0+O(x^2),K(x)=K_0+O(x^2),$$
where $I_0,J_0,K_0$ are the standard complex structures on $\HH^n$, and $I,J,K$ denote also the complex structures on $\HH^n$ obtained using the
identification of the ball in $\HH^n$ with $B(x_0,\eps)\subset M$ via the diffeomorphism $\phi$, and $O(x^2)$ are uniform with respect to $x\in C$.
\end{corollary}
{\bf Proof.} Use the geodesic coordinates on $M$ with respect to the Obata connection with the center at $x_0$. \qed

\subsection{The Douady space.}\label{Ss:Douady}
In this section we review basic properties of the Douady spaces which are a complex analytic analogue of the Hilbert schemes.
This material will be useful to define and to work with the Penrose transform.

Let $Z$ be a complex analytic space (not necessarily compact, and local rings may contain nilpotent elements).
\begin{theorem}[Douady \cite{douady}]
Consider the contravariant functor from the category of complex analytic spaces to the category of sets given by
\begin{eqnarray}\label{Def:Douady}
S\mapsto \{A\subset Z\times S\big| A\mbox{ is closed subspace, flat and proper over } S \},
\end{eqnarray}
where flatness and properness are with respect to the obvious projection $$p\colon A\to S.$$
Then this functor is representable by a complex analytic space we denote by $Dou(Z)$.
\end{theorem}

Informally speaking, $Dou(Z)$ parameterizes all compact subspaces of $Z$.
The following proposition is proven by Daniel Barlet in Appendix below.

\begin{proposition}\label{P:Inclusion-Douady}
Let $Z$ be a complex analytic space, and $Y\subset Z$ is a closed subspace.
Then morphism of spaces $Dou(Y)\to Dou(Z)$ given by the obvious morphism of functors they represent, is a closed imbedding.
\end{proposition}

Now we will need a space which parameterizes compact subspaces of $Z$ together with a point on it. More precisely we have the following well known result.
\begin{theorem}\label{T:scheme+point}
Consider the contravariant functor from the category of complex analytic spaces to the category of sets given by
\begin{eqnarray}\label{Def:Douady+point}
S\mapsto\{A\subset Z\times S,\, f\colon S\to A|\, A \mbox{ is as in (\ref{Def:Douady}) and } p\circ f=Id_S\},
\end{eqnarray}
where, as previously, $p\colon A\to S$ is the obvious projection.

(1) This functor is representable by a complex analytic space we denote by $F(Z)$.

(2) We have two natural morphisms of spaces
$$\tau\colon F(Z)\to Dou(Z) \mbox{ and } \eta\colon F(Z)\to Z$$
induced by morphisms of corresponding functors:
the first is given by forgetting the map $f$ in the pair $(A,f)$, the second is given by forgetting the $A$. \footnote{More precisely, to define $\eta$ consider the composition
$S\overset{f}{\to}A\subset Z\times S\to Z$. We got a morphism $S\to Z$. Namely we obtained a functorial in $S$ map of sets $Mor(S,F(Z))\to Mor(S,Z)$.}
Then the induced morphism $$\eta\times\tau\colon F(Z)\to Z\times Dou(Z)$$ is a closed imbedding.

(3) Let $Y\subset Z$ be a closed subspace. Then
$$F(Z)\times_{(Z\times Dou(Z))}(Y\times Dou(Y))=F(Y),$$
where we have used the natural closed imbedding $Y\times Dou(Y)\subset Z\times Dou(Z)$.
\end{theorem}
{\bf Proof.} We present a proof for the sake of completeness.

(1) Let $F\subset Z\times Dou(Z)$ be the so the called universal family, namely $F$ is a closed subspace, flat and proper over $Dou(Z)$ corresponding to the
identity morphism $Id\colon Dou(Z)\to Dou(Z)$. By the abuse of notation we will also denote by $F$ the contravariant functor from complex spaces to sets given by
the Yoneda construction, i.e. by $S\mapsto Hom(S,F)$. $F$ has the obvious universal property that if $h\in Hom(S,Dou(Z))$ corresponds to a closed subspaces $A\subset Z\times Dou(Z)$
then $A=F\times_{Dou(Z),h}S\subset Z\times S$.

Let us denote by $F'$ the functor (\ref{Def:Douady+point}). We have the morphism of functors
$$\cg\colon F'\to F$$ as follows. For any complex space $S$ we have
\begin{eqnarray*}
F'(S)=\{A\subset Z\times S,\, f\colon S\to A|\, A \mbox{ is as in (\ref{Def:Douady}) and } p\circ f=Id_S\}=\\
\{f\colon S\to F\times_{Dou(Z),h}S\subset Z\times S|\,p\circ f=Id_S\},
\end{eqnarray*}
where $h$ in the last row is the morphism $h\colon S\to Dou(Z)$ corresponding to $A$.
Then define $\cg(S)\colon F'(S)\to F(S)$ by $\cg(f):=p_F\circ f$, where $$p_F\colon F\times_{Dou(Z),h}S\to F$$ is the obvious morphism.
The definition of $\cg$ on morphisms is obvious.
Clearly $\cg$ is a morphism of functors. To see that $\cg$ is an isomorphism of functors, we have to check that $\cg(S)\colon F'(S)\to F(S)$
is an isomorphism of sets. Let us describe the inverse map $F(S)\to F'(S)$. Let $g\in F(S)=Hom(S,F)$. Let us consider the commutative diagram
$$\square[S`S`F`Dou(Z);id`g`\tilde g`],$$
where $\tilde g$ is the composition $S\overset{g}{\to}F\inj Z\times Dou(Z)\to Dou(Z)$, and
the lower horizontal morphism is the obvious projection. This diagram induces the morphism
$$\hat g\colon S\to F\times_{Dou(Z)}S.$$
Hence we have $\cg(S)^{-1}(g)=(F\times_{Dou(Z)}S, \hat g)\in F'(S)$. Part (1) is proved.

\hfill

(2) This immediately follows from the fact that the $F$ from part (1) is a closed subspace of $Z\times Dou(Z)$.

\hfill

(3) We have
\begin{eqnarray*}
\left(F(Z)\times_{(Z\times Dou(Z))}(Y\times Dou(Y))\right)(S)=\\
F(Z)(S)\times_{(Hom(S,Z)\times Dou(Z)(S))}\left(Hom(S,Y)\times Dou(Y)(S)\right)=\\
\{A\subset Z\times S,\,f\colon S\to A\}\times_{(Hom(S,Z)\times\{B\subset Z\times S\})}\left(Hom(S,Y)\times\{C\subset Y\times S\}\right)=\\
\{A\subset Y\times S,\, f\colon S\to A\}=F(Y)(S),
\end{eqnarray*}
where the subspaces $A,B,C$ are assumed to be closed in the corresponding spaces and flat and proper over $S$.
\qed


\hfill

Let us apply these general constructions for the twistor space $Z$ of a connected hypercomplex manifold $(M^{4n},I,J,K)$.
But first we remind the following well known fact.
\begin{proposition}\label{P:real-analyticity}
For all complex structures $L_a,\, a\in S^2$ the underlying real analytic structures on complex manifolds $(M,L_a)$ coincide
with each other. Thus $M$ has a canonical real analytic structure.
\end{proposition}

\begin{proposition}\label{P:complexification}
Consider the map $\iota\colon M\to Dou(Z)$ defined such that for $x\in M$ the image $\iota(x)$ is the fiber over $x$ of the natural
map $Z\to M$, which is isomorphic to $\CC\PP^1$.

Then there exists a neighborhood of the image $\iota(M)$ inside $Dou(Z)$ which is a smooth complex manifold of complex dimension $4n$. The map $\iota$,
as a map into this neighborhood, is a map of real analytic manifolds which is closed imbedding. Moreover $\iota(M)$ is totally real submanifold of $Dou(Z)$,
hence a small neighborhood of  $\iota(M)$ can be considered as a complexification of $M$.
\end{proposition}

We failed to find a reference to Propositions \ref{P:real-analyticity} and \ref{P:complexification}, so we will outline the argument for
the sake of completeness. First the normal bundle $\cn$ of any fiber $\CC\PP^1$ over any point $x\in M$ is isomorphic to $\co(1)^{\oplus 2n}$ (see \cite{HKLR}, p. 556).
Hence $H^1(\CC\PP^1,\cn)=0,\, \dim H^0(\CC\PP^1,\cn)=4n$. These conditions and the Kodaira theorem \cite{kodaira} imply that $\iota(M)$ has an open neighborhood in $Dou(Z)$
which is a smooth complex manifold of complex dimension $4n$. The map $\iota$ from $M$ to this neighborhood is smooth injective immersion (in the sense of real manifolds) which
is moreover proper,
hence $\iota$ is a diffeomorphism onto its image $\iota(M)$ which is a smooth submanifold.
By the abuse of notation, we will denote by $Dou(Z)$ throughout the rest of the paper an open neighborhood of $\iota(M)$ which is smooth; it can be made arbitrarily smaller
in further arguments.

In order to define a real analytic structure on $M$, consider the involution $\sigma\colon Z\to Z$ given by
$$\sigma(a,x)=(-a,x),\, \mbox{ for any } a\in S^2,\, x\in M.$$
$\sigma$ is an anti-holomorphic involution. It induces an anti-holomorphic involution $\tilde \sigma$ on $Dou(Z)$ (roughly put, it sends a subspace to its image under $\sigma$
which is well defined). Consider the subset of fixed points $Dou(Z)^{\tilde\sigma}$. It is not hard to see that in a small enough neighborhood of $\iota(M)$ this set coincides with
$\iota (M)$. Any anti-holomorphic map of a complex manifold is real analytic. The set of fixed points of an anti-holomorphic involution of any complex manifold is a real analytic smooth submanifold
of half dimension, whose small neighborhood is its complexification. Hence $\iota(M)$ is a real analytic submanifold of $Dou(Z)$. Hence, via the diffeomorphism $\iota\colon M\tilde\to\iota(M)$,
we get a real analytic structure on $M$.
It remains to show that this real analytic structure coincides with the real analytic structure underlying the complex manifold structure on $(M,L_a)$ for any $a\in S^2$.

Consider the diagram of maps
$$S^2=\CC\PP^1\overset{\xi}{\leftarrow}Z\overset{\eta}{\leftarrow}F(Z)\overset{\tau}{\to}Dou(Z)\supset \iota(M)=Dou(Z)^{\tilde\sigma},$$
where $F(Z)$ is the scheme from Theorem \ref{T:scheme+point}(1) parameterizing  subschemes of $Z$ with a point on them, $Dou(Z)$ denotes a small neighborhood of $M$
as we have agreed above, $\tau$ and $\eta$ are the natural holomorphic maps from Theorem \ref{T:scheme+point}(2), and $\xi$ is the obvious projection which is also holomorphic.
For such a small $Dou(Z)$, the scheme $F(Z)$ is also smooth complex manifold.
Since $\iota(M)$ is a real analytic submanifold, and $\tau$ is complex analytic, in particular real analytic, submersion, the preimage $\tau^{-1}(\iota(M))$ is a real analytic submanifold.
Let $\eta'$ denote the restriction of $\eta$ to $\tau^{-1}(\iota(M))$.
Then $$\eta'\colon \tau^{-1}(\iota(M))\to Z$$ is an isomorphism of real analytic manifolds. Fix any $a\in S^2$. Then $\xi^{-1}(a)$ is the complex submanifold $(M,L_a)\subset Z$.
It follows that $\eta'^{-1}((M,L_a))\subset \tau^{-1}(\iota(M))$ is a real analytic submanifold which is isomorphic to $(M,L_a)$ as a real analytic manifold.

But the real analytic map $$\tau\colon \eta'^{-1}((M,L_a))\to \iota(M)$$
is a diffeomorphism. Hence it is an isomorphism of real analytic manifolds. Hence $$\tau\circ \eta'^{-1}\colon (M,L_a)\to \iota(M)$$
is an isomorphism of real analytic manifolds. But $\tau\circ \eta'^{-1}=\iota$. Hence we see that all the real analytic structures on $M$ defined
by the complex structures $L_a$ are equal to each other and to the previously defined real analytic structure. \qed

\subsection{Andreotti-Grauert theorem and twistor spaces.}\label{Ss:AG}
In this section we remind the Andreotti-Grauert theorem and prove its consequence for twistor
spaces of hypercomplex manifolds. We refer for the book \cite{demailly-book}, Ch. IX, \S\S 2A and 4F, and for more thorough treatment of
the Andreotti-Grauert theory to \cite{henkin-leiterer}.

\begin{definition}\label{Def:q-complete}
A complex manifold $X$ of complex dimension $n$ is called {\itshape strongly $q$-complete} if there exists a $C^\infty$-smooth
function $\psi\colon X\to\RR$ such that

$\bullet$ $\psi$ is an exhaustion function, i.e. all sublevel sets $\{z\in X|\, \psi(z)<c\},\, c\in \RR$ are relatively compact in $X$;

$\bullet$ $\psi$ is strongly $q$-convex, namely for any point $z\in X$ the complex Hessian $\frac{\pt^2\psi(z)}{\pt z_i\pt \bar z_j}$
has at least $n-q+1$ positive eigenvalues.
\end{definition}
A consequence of the Andreotti-Grauert theorem says:
\begin{theorem}\label{T:AG}
Let $X$ be a stongly $q$-complete complex manifold. Then for any coherent analytic sheaf $\cs$ on $X$ one has
$$H^k(X,\cs)=0 \mbox{ for } k\geq q.$$
\end{theorem}

We will apply this result in the following situation.

\begin{proposition}\label{P:twistors-2-complete}
Let $(M^{4n},I,J,K)$ be a hypercomplex manifold. Let $g$ be an arbitrary infinitely smooth Riemannian metric on $M$.
Let $C\subset M$ be a compact subset. Then there exists $\eps_0>0$ such that for any $\eps\in (0,\eps_0)$ and any $z\in C$
the twistor space of the open ball $B(z,\eps)$ is strongly 2-complete.
\end{proposition}
{\bf Proof.} Let us assume the contrary. Then there exist a sequence of points $\{z_i\}\subset K$ and a sequence of positive numbers $\{\eps_i\}$
such that $z_i\to z$, $\eps_i\to 0$, and $B(z_i,\eps_i)$ is not strongly 2-complete for any $i\geq 1$.

Let $\tilde \psi_i$ (resp. $\tilde\psi$) be the square of the distance function from the point $z_i$ (resp. $z$) on $M$. Then there exists an open neighborhood $U$ of $z$ such that
$\tilde\psi_i$ and $\tilde\psi$ are $C^\infty$-smooth functions on $U$ for $i\gg 1$. It is clear that
\begin{eqnarray}\label{E:sm-convergence1}
\tilde\psi_i\to \tilde\psi \mbox{ in } C^\infty(U).
\end{eqnarray}
Let $\psi_i\colon S^2\times U\to \RR$ (resp. $\psi\colon S^2\times U\to \RR$) be the composition of $\tilde\psi_i$ (resp. $\tilde\psi$) with the natural projection
$S^2\times U\to U$. Obviously $\psi_i|_{S^2\times B(z_i,\eps_i)}$ is an exhaustion function on $S^2\times B(z_i,\eps_i)$. Let
us show that the function $\psi_i|_{S^2\times B(z_i,\eps_i)}$ is strongly 2-complete for $i\gg 1$ with respect to the standard complex structure on the twistor space, i.e.
the complex Hessian $\left(\frac{\pt^2\psi_i}{\pt w_p\pt\bar w_q}\right)$ has at least $2n$ positive eigenvalues on $S^2\times B(z_i,\eps_i)$, where $w_1,\dots,w_{2n+1}$
are local complex coordinates on the twistor space.

By (\ref{E:sm-convergence1}) we have
\begin{eqnarray}\label{E:sm-convergence2}
\psi_i\to \psi \mbox{ in } C^\infty(S^2\times U).
\end{eqnarray}
Hence, by the continuity of eigenvalues of Hermitian matrices, it suffices to check that
$\left(\frac{\pt^2\psi}{\pt w_p\pt\bar w_q}(x)\right)$ has at least $2n$ positive eigenvalues for any $x\in S^2\times\{z\}$.
Thus let $x=(a,z)$, $ a\in S^2$ is a complex structure on $M$. By the definition of the complex structure on the twistor space, it suffices to show that the $a$-complex Hessian of $\tilde\psi$
at the point $z$ is positive definite. Since $d\tilde\psi(z)=0$, the real Hessian $D^2\tilde\psi(z)\in Sym^2T^*_zM$ is well defined and independent of any coordinates.
Moreover, since $\tilde\psi$ is the square of the distance function from $z$, $D^2\tilde\psi(z)$ is a positive definite quadratic form on $T_zM$.

Again, since $d\tilde \psi(z)=0$, the $a$-complex Hessian of $\tilde\psi$ at $z$, which we denote by $Hess_{a-\CC}\tilde\psi(z)$, depends only on $D^2\tilde\psi(z)$, more precisely:
\begin{eqnarray*}
Hess_{a-\CC}\tilde\psi(z)(V)=\frac{1}{2}\int_0^{2\pi}D^2\tilde\psi(z)(V\cdot e^{at})dt
\end{eqnarray*}
for any $V\in T_zM$. The last formula also implies that
$$Hess_{a-\CC}\tilde\psi(z)>0.$$
\qed

\subsection{The Penrose transform.}\label{Ss:Penrose-Baston}
We recall that, by abuse of notation, throughout the rest of the paper we denote by $Dou(Z)$ a small neighborhood of $\iota(M)$ in the Douady space.
Recall also that we have diagram
$$\CC\PP^1\overset{\xi}{\leftarrow}Z\overset{\eta}{\leftarrow}F(Z)\overset{\tau}{\to}Dou(Z),$$
where now $F(Z)$ denotes the preimage under $\tau$ of this neighborhood.

Let $\cl$ be a coherent sheaf of $Z$. The {\itshape Penrose transform} of $\cl$ is $R\tau_*\eta^{-1}\cl\in D^+(Sh_{Dou(Z)})$, where $Sh_Y$ denote the category of sheaves
of vector spaces on a topological space $Y$,
$\eta^{-1}$ denotes the pull-back in the category of sheaves of vector spaces under the map $\eta$, $D^+(Sh_Y)$ denotes the bounded from below
derived category of $Sh_Y$, and $R\tau_*\colon D^+(Sh_{F(Z)})\to D^+(Sh_{Dou(Z)})$ in the derived push-forward map.

To describe the Penrose transform of some locally free coherent $\co_Z$-modules $\cl$,
Baston \cite{baston} has suggested to use the hypercohomology spectral sequence. Differentials of the sequence are differential operators. \footnote{Unfortunately the word "differential"
has here two completely different meanings.} For us it will be important that for some specific choice of $\cl$ one of these operators is proportional to our quaternionic Hessian $\pt\pt_J$.
This cohomological interpretation of this operator will be useful in the proof of our main proposition. This specific choice of $\cl$ and the cohomological interpretation
will be discussed in subsequent sections.

\hfill

Let us assume now that $\cl$ is a locally free $\co_Z$-module of finite rank. Let us consider the resolution $\Ome^\bullet(\cl)$
of the sheaf $\eta^{-1}(\cl)$ by the relative de Rham complex:
\begin{eqnarray}\label{E:deRham}
0\to \eta^{-1}(\cl)\to \eta^*\cl\overset{d}{\to}\eta^*\cl\otimes_{\co}\Ome^1_{F/Z}\overset{d}{\to}\eta^*\cl\otimes_{\co}\Ome^2_{F/Z}\overset{d}{\to}\dots,
\end{eqnarray}
where $\eta^*\cl$ denotes the pull-back of $\cl$ in the category of coherent analytic sheaves, the morphism $\eta^{-1}(\cl)\to \eta^*\cl$ is the
obvious one, $\Ome^i_{F/Z}$ is the sheaf of $\eta$-relative $i$-forms on $F(Z)$, and $d$ is the relative de Rham differential.
Thus $R\tau_*(\eta^{-1}(\cl))$ is isomorphic to $R\tau_*(\Ome^\bullet(\cl))$. The sequence of sheaves (\ref{E:deRham}) is exact.

To compute the push-forward $R\tau_*$ of a complex of sheaves $$\cf=[\cdots\to F^i\to F^{i+1}\to \cdots]$$ there exists a standard hypercohomology spectral sequence
$E^{p,q}_r(\cf)$ converging to $R\tau_*^{p+q}(\cf)$ whose first term is equal to
$$E^{p,q}_1(\cf)=R^q\tau_*(F^p).$$

In our situation $\tau$ is a proper map with 1-dimensional fibers, and all the sheaves $F^i$ are coherent. It follows that
$$E^{p,q}_1(\Ome^\bullet(\cl))=0 \mbox{ unless } q=0,1 \mbox{ and } p\geq 0.$$

In this paper we will need only the sheaves $\cl:=\xi^*(\co(-l)),\, l\geq 2$. The case $l=2$ will be the most important.
Below we will discuss it in detail.

\hfill

In all our discussions involving the Penrose transform and the Baston operators it will usually be convenient to shrink the manifold $M$ to be a small enough neighborhood of some point
and also shrink the neighborhood of $\iota (M)$ in $Dou(Z)$ such that the fibers of the map $\eta\colon F(Z)\to Z$ will be contractible. This is always possible due to the following
proposition which is probably a folklore.


\begin{proposition}\label{P:contractible}
Let $(X_1,\rho_1),\, (X_2,\rho_2)$ be infinitely smooth Riemannian manifolds. Let $F\subset X_1\times X_2$ be a closed subset which is
a smooth submanifold. Assume that the natural projections
$$p_{1,2}\colon F\to X_{1,2}$$
are submersions, and $p_2$ is proper. Let $C\subset X_2$ be a compact subset.

Then there exists $\eps_0>0$ such that any $\eps\in (0,\eps_0)$ has the following properties: for any $a\in C$ and any $b$ belonging to the $\eps$-neighborhood
(with respect to the metric $\rho_1$) of the compact set $p_1(p_2^{-1}(a))$ the manifold $F\cap \left(\{b\}\times B(a,\eps)\right)$ is diffeomorphic to $\RR^k$, where
$B(a,\eps)$ is the open ball with respect to the metric $\rho_2$.
\end{proposition}
{\bf Proof.} Assume that the proposition is wrong. Then there exist:

$\bullet$ a sequence $\{\eps_i\}$ of positive real numbers such that $\lim_{i\to\infty}\eps_i=0$;

$\bullet$ a sequence $\{a_i\}\subset C$ such that $\lim_{i\to\infty}a_i=a\in C$;

$\bullet$ a sequence $\{b_i\}\subset X_1$ such that $b_i$ belongs to $\eps_i$-neighborhood of $p_1(p_2^{-1}(a_i))$ and $\lim_{i\to\infty}b_i=b$ belonging
necessarily to $p_1(p_2^{-1}(a))$, and such that $F\cap \left(\{b_i\}\times B(a_i,\eps_i)\right)$ is {\itshape not} diffeomorphic to $\RR^k$ for all $i\geq 1$.

\hfill

Let us define functions $R_i,R\colon X_2\to \RR$ by
\begin{eqnarray*}
R_i(x)=dist^2(x,a_i),\\
R(x)=dist^2(x,a),
\end{eqnarray*}
where $dist$ denotes the distance with respect to the metric $\rho_2$.
Then there exists $\delta>0$ such that

$\bullet$ $R_i,R\in C^\infty(B(a,\delta))$ for all $i\gg 1$;

$\bullet$ $R_i\to R$ in $C^\infty(B(a,\delta))$;

$\bullet$ the only critical point of $R$ in $B(a,\delta)$ is $a$; it is non-degenerate and it is a point of global minimum of $R$.

\hfill

Next there exists an open neighborhood $U\subset X_1$ of $b$, an open neighborhood $V\subset X_1\times X_2$ of $(b,a)$, and a diffeomorphism
$$\phi\colon V\cap F\tilde\to U\times \RR^k$$
such that the following diagram is commutative:
$$\qtriangle[V\cap F`U\times\RR^k`U;\phi`p_1`p_U],$$
where $p_U$ is the natural projection. Let us denote
$$q:=p_2\circ \phi^{-1}\colon U\times \RR^k\to X_2.$$
Notice that for any $u\in U$ the restriction of $q$ to $\{u\}\times \RR^k$ is
an injective immersion; we denote this restriction by $q_u$ and consider it as a map
$$q_u\colon \RR^k\to X_2.$$
Evidently we may and will assume that
$$\phi(b,a)=(b,0)\in U\times \RR^k.$$
Let us define functions $f_{u,i},f_u\colon \RR^k\to \RR$ by
$$f_{u,i}:=R_i\circ q_u,\, f_u:=R\circ q_u.$$
It is clear that there exists $\kappa>0$ such that for all $u\in U$ and $i\gg 1$ the functions
$f_{u,i}, f_{u}$ are $C^\infty$-smooth on the Euclidean ball $B(0,\kappa)\subset \RR^k$ and
\begin{eqnarray}\label{E:convergence000}
f_{b_i,i}\to f_b \mbox{ in } C^\infty(B(0,\kappa)) \mbox{ as } i\to \infty.
\end{eqnarray}
Since $q_b(B(0,\kappa))\subset X_2$ is a smooth submanifold containing the point $a$, the restriction of the function $R$ to
$q_b(B(0,\kappa))\cap B(a,\eps)$ for $0<\eps\ll 1$ has the unique critical point which is $a$, and it is non-degenerate.

Let us denote for $0<\eps\ll 1$
$$\Ome_\eps:=q_b^{-1}(B(a,\eps))\subset B(0,\kappa)\subset \RR^k.$$
Clearly $\pt \Ome_\eps=q_b^{-1}(\pt B(a,\eps))$ is a smooth compact manifold. Then the point $0\in \Ome_\eps$ and it is the only critical point of the function
$$f_b\colon \bar\Ome_\eps\to \RR;$$
moreover this critical point is non-degenerate and $f_b|_{\pt \Ome_\eps}\equiv \eps$. Then making $\eps$ small enough we see that
(\ref{E:convergence000}) implies that for $i\gg 1$ each $f_{b_i,i}$ has a unique critical point $c_i$ in $\bar \Ome_\eps$ which belongs to $\Ome_\eps$
and is non-degenerate; moreover $\lim_{i\to \infty}c_i=0$.

Notice that the set $\Ome_i:=\{x\in \RR^k|\, f_{b_i,i}(x)<\eps_i\}$ is diffeomorphic (via $\phi$) to $$F\cap \left(\{b_i\}\times B(a_i,\eps_i)\right).$$
Hence we will get a contradiction if we apply the following general lemma to the function
$$f_{b_i,i}\colon \bar\Ome_i\to \RR$$
for $i\gg 1$.
\begin{lemma}\label{L:Morse}
Let $N$ be a compact manifold with boundary $\pt N$. Let $$f\colon N\to \RR$$ be a smooth function such that

(i) $f$ is constant on the boundary $\pt N$.

(ii) $f$ has a single critical point. Assume also that it does not belong to $\pt N$ and is non-degenerate.

Then $N$ is diffeomoprhic to the closed Euclidean ball.
\end{lemma}
A proof of the lemma uses the basic Morse theory and is left to the reader.
Proposition \ref{P:contractible} is proved.



\subsection{Baston operators.}\label{Ss:Baston-operators} Recall that for a hypercomplex manifold $M$ we have a diagram of holomorphic maps
$$\CC\PP^1\overset{\xi}{\leftarrow}Z\overset{\eta}{\leftarrow}F(Z)\overset{\tau}{\to}Dou(Z),$$
where, as previously, $Dou(Z)$ denotes a small neighborhood of $\iota(M)$ is the Douady space of the twistor space $Z$.

\hfill








In this section we will study in greater detail the Penrose transform of the sheaves $\cl=\xi^*(\co(-l))$ with $l\geq 2$; the case
$l=2$ will be especially important.

First let us remind couple of general facts from complex analysis. For a complex analytic space $X$
and its closed analytic subspace $Y\subset X$ we denote by $Y^{(n)}$ the $n$th infinitesimal neighborhood of $Y$. $Y^{(n)}$ is a closed analytic subspace of $X$.
The following result is essentially due to Grauert. It is a special case of a result mentioned on p. 119 of the book
\cite{banica-stanasila} (in the proof of Theorem 3.4 there).
\begin{theorem}\label{T:grauert}
Let $f\colon X\to Y$ be a proper holomorphic submersive map of smooth complex analytic manifolds.
Let $\cf$ be a locally free sheaf of $\co_X$-modules of finite rank. Assume that for some $i$ the function
$$y\mapsto \dim H^i(X_y,\cf_y)$$
is constant on $Y$, where $X_y=f^{-1}(y)$ and $\cf_y$ is the pull-back of $\cf$ to $X_y$.

Then for any $y\in Y$ and $n\in \ZZ_{\geq 0}$ the canonical map
$$R^if_*\cf\otimes_{\co_Y}\co_{y^{(n)}}\to H^i(X_y^{(n)},\cf_n)$$
is an isomorphism, where $\cf_n$ denotes the pull-back of $\cf$ to $X_y^{(n)}=X\times_Yy^{(n)}$.\footnote{Notice that $X\times_Yy^{(n)}$ coincides
with the $n$th infinitesimal neighborhood of the fiber $X_y$.}
\end{theorem}

We will also use the following result which is an immediate consequence of Demailly's theorem \cite{demailly}.
\begin{theorem}\label{T:demailly}
Let $f\colon X\to Y$ be a proper holomorphic submersive map of smooth complex analytic manifolds.
Let $\cf$ be a sheaf of locally free $\co_X$-modules of finite rank.  Then, in the notion of Theorem \ref{T:grauert}, the function
$$y\mapsto \dim H^i(X_y,\cf_y)$$
is upper semi-continuous for the analytic Zariski topology, more explicitly for any integer $N$ the set
$$\{y\in Y|\, \dim H^i(X_y,\cf_y)\geq N\}$$
is a closed analytic subset.
\end{theorem}

\begin{proposition}\label{P:constant-dim}
Let $M$ be a hypercomplex manifold. Let $p,i\geq 0$ and $l$ be integers. For $x\in Dou(Z)$ let us denote
$F_x:=\tau^{-1}(x)$. Then the function
\begin{eqnarray}\label{E:diff-forms-relat}
x\mapsto \dim H^i\left(F_x,\left((\xi\circ\eta)^*\co(l)\otimes_{\co_F}\Ome^p_{F/Z}\right)_x\right)
\end{eqnarray}
is constant in a neighborhood of $M$ in $Dou(Z)$.
\end{proposition}
{\bf Proof.} Since $Dou(Z)$ is a complexification of $M$, Theorem \ref{T:demailly} easily implies that it suffices
to show that the function (\ref{E:diff-forms-relat}) is constant on $M$; then it will be automatically constant
in a neighborhood of $M$.

Thus let $x\in M$. Then $\xi\circ\eta\colon F_x\to \CC\PP^1$ is an isomorphism. Using this isomorphism we
will identify $F_x=\CC\PP^1$. Also $[F_x]=\iota(x)$. By \cite{HKLR}, p. 556, the normal bundle $N_{F_x}Z$ of $F_x$ in $Z$ is
canonically isomorphic to
\begin{eqnarray}\label{E:normal-bundle}
N_{F_x}Z\simeq \co(1)\otimes_\CC T_xM_I,
\end{eqnarray}
where $T_xM_I$ denotes the tangent space at $x$ to the complex manifold $(M,I)$. For any $z\in \CC\PP^1$ we have
\begin{eqnarray*}
\Ome^1_{F/Z}|_{(z,x)}=\{\psi\in H^0(F_x,N_{F_x}Z)|\, \psi(z)=0\}^*\simeq\\
\{\zeta\in H^0(\CC\PP^1,\co(1))|\,\zeta(z)=0\}^{*}\otimes T^*_xM_I.
\end{eqnarray*}
Since the first factor in the last expression is one dimensional we have
\begin{eqnarray}\label{E:forms-fiber}
\Ome^p_{F/Z}|_{(z,x)}=\{\zeta \in H^0(\CC\PP^1,\co(1))|\,\zeta(z)=0\}^{*\otimes p}\otimes \wedge^p T^*_xM_I.
\end{eqnarray}
It follows that the restriction of $\Ome_{F/Z}^p$ on $F_x$ is canonically isomorphic to
\begin{eqnarray*}
\Ome^p|_{F_x}\simeq \co(1)^{\otimes p}\otimes_\CC \wedge^p T^*_XM_I=\co(p)\otimes_\CC \wedge^pT^*_xM_I.
\end{eqnarray*}
Now
$$\left((\xi\circ \eta)^*\co(l)\otimes_{\co_F}\Ome^p_{F/Z}\right)|_{F_x}\simeq \co(p+l)\otimes_\CC \wedge^pT^*_xM_I.$$
Hence we get the canonical isomorphism
\begin{eqnarray}\label{E:cohomology-isomorphism}
H^i(F_x,(\xi\circ \eta)^*\co(l)\otimes_{\co_{F}} \Ome^p_{F/Z})\simeq H^i(\CC\PP^1,\co(p+l))\otimes\wedge^pT^*_xM_I.
\end{eqnarray}
The proposition follows. \qed

\hfill

Proposition \ref{P:constant-dim} and Theorem \ref{T:grauert} immediately imply the following corollary.
\begin{corollary}\label{Cor:grauert-twistor}
Let $M$ be a hypercomplex manifold. Let $p,i\geq 0$ and $l$ be integers.
There exists a neighborhood of $\iota(M)$ in $Dou(Z)$ such that for any $y$ from this neighborhood and any $n,p\in \ZZ_{\geq 0}, l\in \ZZ$ the canonical morphism
$$R^i\tau_*((\xi\circ\eta)^*\co(l)\otimes_{\co_F}\Ome_{F/Z}^p)\otimes_{\co_X}\co_{y^{(n)}}\to H^i\left(F_y^{(n)},\left((\xi\circ\eta)^*\co(l)_{\co_F}\Ome_{F/z}^p\right)_n\right)$$
is an isomorphism.
\end{corollary}
From the formula (\ref{E:cohomology-isomorphism}) and Corollary \ref{Cor:grauert-twistor} we immediately deduce
\begin{proposition}\label{P:push-identification}
Let $M$ be a hypercomplex manifold. Let $p,i\geq 0$ and $l$ be integers. Then the sheaf of $\co_{Dou(Z)}$-modules $R^i\tau_*((\xi\circ\eta)^*\co(l)\otimes_{\co_F}\Ome_{F/Z}^p)$
is the sheaf of holomorphic sections of a vector bundle over $Dou(Z)$ whose restriction to $M$ is canonically isomorphic to sheaf of real analytic sections of the vector bundle
over $M$
$$H^i(\CC\PP^1,\co(p+l))\otimes_\CC \wedge^{p,0}_IM,$$
where $\wedge^{p,0}_IM$ is the bundle of $(p,0)$-forms on the complex manifold $(M,I)$.
\end{proposition}

Let us recall that
\begin{eqnarray*}
H^0(\CC\PP^1,\co(k))=0 \mbox{ iff } k\leq -1;\\
H^1(\CC\PP^1,\co(k))=0 \mbox{ iff } k\geq -1.
\end{eqnarray*}
Moreover if we fix an isomorphism of sheaves $\ome\simeq \co(-2)$ where $\ome$ is the sheaf of holomorphic 1-forms on $\CC\PP^1$, and denote $\cv:=H^0(\CC\PP^1,\co(1))$,
then in the cases when the cohomology groups do not vanish one has
\begin{eqnarray*}
H^0(\CC\PP^1,\co(k))\simeq Sym^k\cv\mbox{ for }k\geq 0,\\
H^1(\CC\PP^1,\co(k))\simeq Sym^{-k-2}\cv^*\mbox{ for } k\leq -2,
\end{eqnarray*}
where in the second isomorphism we have used the Serre's duality and the above chosen isomorphism $\ome\simeq \co(-2)$.
Recall also that $\dim \cv=2$.

Since the vector bundles $\Lam^{p,0}_IM$ are real analytic on $M$, we can consider their complexifications on $Dou(Z)$
which are holomorphic vector bundles on $Dou(Z)$. We will denote them by the same notation $\Lam^{p,0}_IM$.

Let now $-2n\leq l \leq -2$, $\cl=\co(l)$. Then Proposition \ref{P:push-identification} implies that the first
term of the spectral sequence computing $R\tau_*(\Ome^\bullet(\cl))$ is equal to
\begin{eqnarray*}
E^{p,0}_1(\Ome^\bullet(\cl))\simeq \Lambda^{p,0}_IM\otimes Sym^{p+l}\cv \mbox{ for } l+p\geq 0,\\
E^{p,1}_1(\Ome^\bullet(\cl))\simeq \Lambda^{p,0}_IM\otimes Sym^{-p-l-2}\cv^*\mbox{ for } l+p\leq -2,\\
E^{p,q}_1(\Ome^\bullet(\cl))=0 \mbox{ otherwise.}
\end{eqnarray*}
Thus the first term of the spectral sequence looks as on Fig. 1:
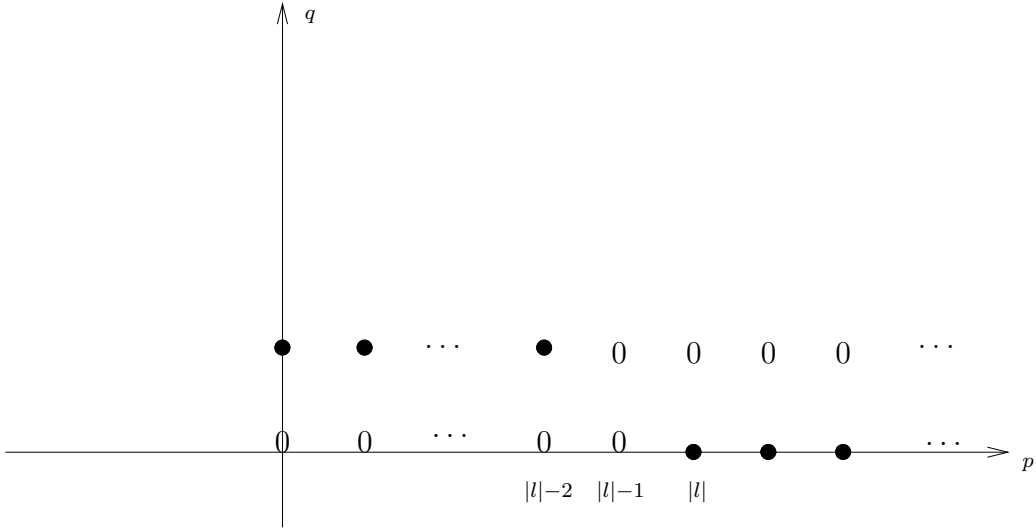
\begin{figure}[h]
\setlength{\unitlength}{0.00087489in}
\begingroup\makeatletter\ifx\SetFigFontNFSS\undefined%
\gdef\SetFigFontNFSS#1#2#3#4#5{%
  \reset@font\fontsize{#1}{#2pt}%
  \fontfamily{#3}\fontseries{#4}\fontshape{#5}%
  \selectfont}%
\fi\endgroup%
{\renewcommand{\dashlinestretch}{30}
\begin{picture}(6254,3210)(0,-10)
\put(1677,1092){\blacken\ellipse{90}{90}}
\put(1677,1092){\ellipse{90}{90}}
\put(2172,1092){\blacken\ellipse{90}{90}}
\put(2172,1092){\ellipse{90}{90}}
\put(3252,1092){\blacken\ellipse{90}{90}}
\put(3252,1092){\ellipse{90}{90}}
\put(4152,462){\blacken\ellipse{90}{90}}
\put(4152,462){\ellipse{90}{90}}
\put(4602,462){\blacken\ellipse{90}{90}}
\put(4602,462){\ellipse{90}{90}}
\put(5052,462){\blacken\ellipse{90}{90}}
\put(5052,462){\ellipse{90}{90}}
\path(1677,12)(1677,3162)
\path(1707.000,3042.000)(1677.000,3162.000)(1647.000,3042.000)
\path(12,462)(6042,462)
\path(5922.000,432.000)(6042.000,462.000)(5922.000,492.000)
\put(1632,462){\makebox(0,0)[lb]{\smash{{\SetFigFontNFSS{12}{14.4}{\rmdefault}{\mddefault}{\updefault}0}}}}
\put(2127,462){\makebox(0,0)[lb]{\smash{{\SetFigFontNFSS{12}{14.4}{\rmdefault}{\mddefault}{\updefault}0}}}}
\put(3207,462){\makebox(0,0)[lb]{\smash{{\SetFigFontNFSS{12}{14.4}{\rmdefault}{\mddefault}{\updefault}0}}}}
\put(2577,552){\makebox(0,0)[lb]{\smash{{\SetFigFontNFSS{12}{14.4}{\rmdefault}{\mddefault}{\updefault}$\ldots$}}}}
\put(3657,462){\makebox(0,0)[lb]{\smash{{\SetFigFontNFSS{12}{14.4}{\rmdefault}{\mddefault}{\updefault}0}}}}
\put(3657,1002){\makebox(0,0)[lb]{\smash{{\SetFigFontNFSS{12}{14.4}{\rmdefault}{\mddefault}{\updefault}0}}}}
\put(4557,1002){\makebox(0,0)[lb]{\smash{{\SetFigFontNFSS{12}{14.4}{\rmdefault}{\mddefault}{\updefault}0}}}}
\put(5007,1002){\makebox(0,0)[lb]{\smash{{\SetFigFontNFSS{12}{14.4}{\rmdefault}{\mddefault}{\updefault}0}}}}
\put(4107,1002){\makebox(0,0)[lb]{\smash{{\SetFigFontNFSS{12}{14.4}{\rmdefault}{\mddefault}{\updefault}0}}}}
\put(2532,1092){\makebox(0,0)[lb]{\smash{{\SetFigFontNFSS{12}{14.4}{\rmdefault}{\mddefault}{\updefault}$\ldots$}}}}
\put(5547,507){\makebox(0,0)[lb]{\smash{{\SetFigFontNFSS{12}{14.4}{\rmdefault}{\mddefault}{\updefault}$\ldots$}}}}
\put(5502,1092){\makebox(0,0)[lb]{\smash{{\SetFigFontNFSS{12}{14.4}{\rmdefault}{\mddefault}{\updefault}$\ldots$}}}}
\put(1812,3072){\makebox(0,0)[lb]{\smash{{\SetFigFontNFSS{12}{14.4}{\rmdefault}{\mddefault}{\updefault}$\scriptstyle q$}}}}
\put(6132,372){\makebox(0,0)[lb]{\smash{{\SetFigFontNFSS{12}{14.4}{\rmdefault}{\mddefault}{\updefault}$\scriptstyle p$}}}}
\put(3132,192){\makebox(0,0)[lb]{\smash{{\SetFigFontNFSS{12}{14.4}{\rmdefault}{\mddefault}{\updefault}$\scriptstyle |l|-2$}}}}
\put(3567,192){\makebox(0,0)[lb]{\smash{{\SetFigFontNFSS{12}{14.4}{\rmdefault}{\mddefault}{\updefault}$\scriptstyle |l|-1$}}}}
\put(4052,192){\makebox(0,0)[lb]{\smash{{\SetFigFontNFSS{12}{14.4}{\rmdefault}{\mddefault}{\updefault}
$\scriptstyle |l|$}}}}
\end{picture}
}
\caption{The first term $E^{p,q}_1(\Ome^\bullet({\mathcal O}(l)))$ of
the spectral sequence with $-2n\leq l\leq -2$.}
\end{figure}

Hence the only non-zero differential $d_r$ for $r=2$ in the spectral sequence
is a morphism of sheaves
$$\Delta\colon \Lambda^{-l-2,0}_IM\to \Lambda^{-l,0}_IM.$$
$\Delta$ is called the Baston operator. It is the second order differential operator with holomorphic coefficients \cite{gindikin-henkin}, \cite{baston}.
All other non-zero differetials $d_r$ in the spectral sequence have $r=1$ and are first order differential operators with holomorphic coefficients \cite{gindikin-henkin}, \cite{baston}.
By the construction of the spectral sequence we obtain the following complex of locally free $\co_{Dou(Z)}$-modules for any $-2n\leq l\leq -2$:
\begin{eqnarray*}
0\to Sym^{-l-2}\cv^*\overset{d_1}{\to}\Lam^{1,0}_IM\otimes Sym^{-l-3}\cv^*\overset{d_1}{\to}\dots\overset{d_1}{\to}\Lam^{-l-3,0}_IM\otimes \cv^*\overset{d_1}{\to}\\
\overset{d_1}{\to} \Lam^{-l-2,0}_IM\overset{\Delta}{\to}\Lam^{-l,0}M\overset{d_1}{\to}\\
\overset{d_1}{\to}\Lam^{-l+1,0}_IM\otimes\cv\overset{d_1}{\to}\Lam^{-l+2,0}_IM\otimes Sym^2\cv\overset{d_1}{\to}\dots\overset{d_1}{\to}\Lam^{2n,0}_IM\otimes Sym^{2n+l}\cv\to 0.
\end{eqnarray*}
We will denote this complex of sheaves on $Dou(Z)$ by $\BC_l$. The complex of global sections of $\BC_l$ is often called in the literature
the Baston complex.

We claim that $R\tau_*\Ome^\bullet(\cl)[1]$ is functorially isomorphic to the Baston complex; this immediately follows from the following general lemma whose proof
is left to the reader.
\begin{lemma}\label{L:homol-alg}
Let $F\colon \ca\to \cb$ be a left exact functor between two abelian categories, such that $\ca$ has sufficiently many injective objects.
Let $RF\colon D^+\ca\to D^+\cb$ denote the derived functor between the bounded below derived categories. Let us fix an integer $a\in \ZZ$. Let us denote by $\ct$ the full subcategory of
the category $Kom^+(\ca)$ of bounded below complexes, consisting of complexes $A^\bullet$ satisfying the following vanishing properties:
\begin{eqnarray*}
R^qA^a=0 \mbox{ for every } q;\\
\mbox{ for every } p< a\,\, R^qF(A^p)=0 \mbox{ for every } q\ne 1;\\
\mbox{ for every } p> a\,\, R^qF(A^p)=0 \mbox{ for every } q\ne 0.
\end{eqnarray*}
Define the functor $G\colon\ct\to D^+(\cb)$ by
$$G(A^\bullet)=[\dots\to R^1F(A^{a-2})\to  R^1F(A^{a-1})\overset{\Delta}{\to} F(A^{a+1})\to F(A^{a+2})\to \dots],$$
where $\Delta$ is the only non-trivial differential in the standard spectral sequence whose first term is $E^{p,q}_1=R^qF(A^p)$ and which converges to
cohomology of $RF(A^\bullet)$; the other differentials in $G(A^\bullet)$ are induced by the corresponding differentials of $A^\bullet$, and the whole complex is shifted such that the term
$R^1F(A^{a-1})$ is located in the degree $a$.

Then the functors $RF,G\colon \ct\to D^+\cb$ are isomorphic.
\end{lemma}

Returning back to the Baston complex, notice that if we restrict this complex to $M$, we will get a complex of sheaves of real analytic sections
of vector bundles, and the differentials of the complex are differential operators with real analytic coefficients of order 1 or 2. We will normalize the Baston complex so that that the first non-zero term
$Sym^{-l-2}\cv^*$ is at cohomological degree 0.

The following theorem was proven in \cite{gindikin-henkin} in the flat case and in \cite{baston} for general quaternionic manifolds.
\begin{theorem}\label{T:Baston-cohomology}
Let us assume that the twistor space $Z$ of a hypercomplex manifold $M^{4n}$ is strongly 2-complete.
Let $Dou(Z)$ be a neighborhood of $\iota(M)$ in the Douady space; assume that $Dou(Z)$ is a Stein manifold.
Assume also that the fibers of the map
$\eta\colon F(Z)\to Z$ are contractible. Let $-2n\leq l\leq -2$ be an integer. Let $\cl:=\xi^*(\co(l))$. Then we have

(1) The Baston complex on $Dou(Z)$ is exact everywhere except cohomological degree 0.

(2) The 0th cohomology of the Baston complex on $Dou(Z)$ is isomorphic to $H^1(Z,\cl)$ via the following maps
\begin{eqnarray*}
H^1(Z,\cl)\tilde\to H^1(Z,\eta_*\eta^{-1}\cl)\simeq H^1(F(Z),\eta^{-1}\cl)\simeq\\
\HH^1(F(Z),\Ome^\bullet(\cl))\simeq \mbox{ 0th cohomology of the Baston complex}.
\end{eqnarray*}

(3) $H^0(Z,\cl)=0$.
\end{theorem}
{\bf Proof.} We present a proof for the sake of completeness. Since $\eta$ has contractible fibers we have for any $i$
$$H^i(Z,\cl)\tilde\to \HH^i(Z,R\eta_*\eta^{-1}\cl)\simeq H^i(F(Z),\eta^{-1}\cl)\simeq \HH^i(F(Z),\Ome^\bullet(\cl)),$$
where $\HH^i$ denotes the $i$th hypercohomology of space with coefficients in a complex.

We also have
\begin{eqnarray*}
\HH^i(F(Z),\Ome^\bullet(\cl))\tilde\to\HH^i(Dou(Z), R\tau_*\Ome^\bullet(\cl))=\HH^{i-1}(Dou(Z), \BC_l).
\end{eqnarray*}
Hence we obtain an isomorphism for every integer $i$
\begin{eqnarray}\label{E:isomorphism-cohomol}
H^i(Z,\cl)\simeq \HH^{i-1}(Dou(Z), \BC_l).
\end{eqnarray}
Since $\BC_l$ lies in the non-negative degrees, the right hand side of equality (\ref{E:isomorphism-cohomol})
vanishes for $i\leq 0$. Hence $H^0(Z,\cl)=0$. This proves part (3).

Since $Z$ is strongly 2-complete, $H^i(Z,\cl)=0$ for $i\geq 2$ by Theorem \ref{T:AG}. Hence
\begin{eqnarray}\label{E:hypercohom-Baston1}
\HH^j(Dou(Z), \BC_l)=0 \mbox{ for } j\geq 1.
\end{eqnarray}
Taking in (\ref{E:isomorphism-cohomol}) $i=1$ we get
\begin{eqnarray}\label{E:hypercohom-Baston2}
H^1(Z,\cl)\simeq \HH^{0}(Dou(Z), \BC_l).
\end{eqnarray}

Finally let us use the assumption that $Dou(Z)$ is the Stein manifold. This implies that positive cohomology of $Dou(Z)$ with coefficients in any coherent sheaf vanish,
i.e. every coherent
sheaf is acyclic with respect to the functor of taking global sections. Since
in $\BC_l$ all sheaves are coherent, and hence acyclic, $\HH^j(Dou(Z), \BC_l)$ is equal to the $j$th cohomology of the complex of global sections of
$\BC_l$, i.e. of the Baston complex.
Hence parts (1) and (2) of the proposition follow from (\ref{E:hypercohom-Baston1}) and (\ref{E:hypercohom-Baston2})
respectively. \qed

\subsection{A technical result on relative de Rham complex.}\label{Ss:relative-de-Rham}
The goal of this section is to prove the following result.
\begin{theorem}\label{T:relative-deRham}
Let $X_1,X_2$ be a smooth complex analytic manifold. Let $F\subset X_1\times X_2$ be a closed complex analytic submanifold such that the projections
$p_{1,2}\colon F\to X_{1,2}$ are submersions, and $p_2$ be proper. Let $a\in X_2$ be a point. Let $C:=p_2^{-1}(a)$. For any non-negative integer $n$ let us denote by
$C^{(n)}$ the $n$th infinitesimal neighborhood of $C$ inside $F$; thus $C^{(n)}\subset F$ is a closed analytic subspace. Denote by
$i_n\colon C^{(n)}\inj F$ the natural imbedding. Let $p_1^{(n)}$ denote the morphism of analytic spaces $p_1^{(n)}\colon C^{(n)}\to (p_1(C))^{(n)}$ induced by $p_1$.
Let $\cl$ be a holomorphic vector bundle over $X_1$.

Then for any non-negative integer $n$ the complex of sheaves on $F$
\begin{eqnarray}\label{E:relative-deRham-explicit}
0\to i_{n*}\left((p_1^{(n)})^{-1}(\cl)\right)\to i_{n*}i_n^*p_1^*(\cl)\to i_{n-1*}i_{n-1}^*(\Omega^1_{F/X_1}\otimes_{\co_F}p_1^*\cl)
\to \\ i_{n-2*}i_{n-2}^*(\Omega^2_{F/X_1}\otimes_{\co_F}p_1^*\cl)\to \dots\to
\end{eqnarray}
is exact. Here $\cl$ denotes also for brevity the pull-back of $\cl$ to $(p_1(C))^{(n)}$ in the category of coherent sheaves;
the first non-zero arrow is the obvious inclusion;
the differentials are induced by the relative de Rham differentials.
\end{theorem}
\begin{remark}
By convention, in complex (\ref{E:relative-deRham-explicit}) the terms $i_{n-q*}i_{n-q}^*(\dots)$ vanish for $n-q<0$.
\end{remark}

Before the proof we will need some preparations.

\begin{lemma}\label{L:submflds}
Let $X_1,X_2$ be a smooth complex analytic manifold. Let $F\subset X_1\times X_2$ be a closed complex analytic submanifold such that the projections
$p_{1,2}\colon F\to X_{1,2}$ are submersions. Let $(a_1,a_2)\in F$.

(1) Then there exist neighborhoods $U_1,U_2$ of $a_1,a_2$ respectively and a smooth closed complex analytic submanifold $R\subset U_2$ such that

$\bullet$ $a_2\in R$;

$\bullet$ $p_1\colon \left((p_2|_F)^{-1}R\right)\cap (U_1\times U_2)\to U_1$ is a biholomorphic diffeomorphism onto its image which is an open subset of $U_1$.

\hfill

(2) Furthermore the point $(a_1,a_2)\in F$ has an open neighborhood $\cu\subset F$ and a biholomorphic diffeomorphism
$$\phi\colon ((p_2|_F)^{-1}(R)\cap \cu)\times B^m\tilde\to \cu,$$
where $B^m$ is an open ball in $\CC^m$, and such that the following diagram commutes:
\begin{eqnarray*}
\square<1`1`1`1;1000`500>[\left((p_2|_F)^{-1}(R)\cap \cu\right)\times B^m`\cu`(p_2|_F)^{-1}(R)\cap \cu`X_1;\phi``p_1`p_1],
\end{eqnarray*}
where the first vertical arrow in the obvious projection.
\end{lemma}
{\bf Proof.} Denote $C:=(p_2|_F)^{-1}(a_2)$. Clearly
$$T_{(a_1,a_2)}C\cap T_{(a_1,a_2)}(p_1|_F)^{-1}(a_1)=\{0\}.$$
Let us choose a linear subspace $L\subset T_{(a_1,a_2)}F$ such that
$$T_{(a_1,a_2)}C\oplus \left(T_{(a_1,a_2)}(p_1|_F)^{-1}(a_1)\right)\oplus L=T_{(a_1,a_2)}F.$$
Since $L\cap T_{(a_1,a_2)}C=\{0\}$, the linear map
$$dp_2\colon L\to dp_2(L)$$
is an isomorphism of vector spaces.

Let us choose a complex analytic submanifold $R\subset X_2$ such that $a_2\in R$ and $T_{a_2}R=dp_2(L)$.
It is easy to see that
$$T_{(a_1,a_2)} \left((p_2|_F)^{-1}(R)\right)=T_{(a_1,a_2)}C\oplus L.$$
Hence
$$p_1\colon (p_2|_F)^{-1}(R)\to X_1$$
is a biholomorphic diffeomorphism of a small neighborhood of $(a_1,a_2)$ onto its image.
This implies part(1) of the lemma.

\hfill

(2) Let us denote $m:=\dim_{\CC}(p_1|_F)^{-1}(a_1)$. Let us choose a holomorphic map
$$f\colon \co\to \CC^m$$
defined on a neighborhood
$\co\subset F$ of $(a_1,a_2)$ such that
$$Ker\, df_{(a_1,a_2)}\cap T_{(a_1,a_2)}\left((p_1|_F)^{-1}(a_1)\right)=\{0\}.$$
Consider the holomorphic map
$$p_1|_F\times f\colon \co\to X_1\times \CC^m.$$
Its differential at $(a_1,a_2)$ is an isomorphism. Then there exists a neighborhood $\cu\subset \co\subset F$
of $(a_1,a_2)$ and an open ball $B^m\subset \CC^m$ such that
$$p_1|_F\times f\colon \cu\to X_1\times B^m$$
is a biholomorphic diffeomorphism. We take
$$\phi:=(p_1|_{\cu\cap p_2^{-1}(R)})^{-1}\times f.$$
It is easy to see that $\cu$ and $\phi$ satisfy part(2) of the lemma.
\qed

\hfill

Let us denote by $\Ome^i_{pol}(V)$ the space of polynomial $i$-forms on a finite dimensional complex vector space $V$.
The de Rham differential $d$ acts on $\Ome_{pol}^\bullet(V)$. The following lemma is a well known folklore. Nevertheless we will include the proof for convenience
of the reader.
\begin{lemma}\label{L:pol-de-rham}
(1) The polynomial de Rham complex $(\Ome^\bullet_{pol}(V),d)$ is exact everywhere but the zero degree.
Its zero degree cohomology equals $\CC$.

(2) Let us write $\Ome_{pol}^i(V)=Sym^\bullet(V^*)\otimes \wedge^iV^*$. If for $\ome\in Sym^p(V^*)\otimes \wedge^iV^*,\, i>0,$ one has
$$d\ome=0$$
then there exists $\eta\in Sym^{p+1}(V^*)\otimes \wedge^{i-1}V^*$ such that
$$\ome=d\eta.$$
\end{lemma}
{\bf Proof}. The $Eu$ denote the Euler vector field on $V$; recall that $Eu$ is the generator of the group of homotheties of $V$:
$(\lam,z)\mapsto \lam z$ where $\lam\in \CC,z\in V$. In coordinates
$$Eu=\sum_{i=1}^{\dim V}z_i\frac{\pt}{\pt z_i}.$$
Let us denote by $d^*$ the operator of contraction of a differential form with $Eu$.
One can easily check the following properties:

\begin{eqnarray}
(d^*)^2=0;\\
d\left(Sym^p(V^*)\otimes \wedge^iV^*\right)\subset Sym^{p-1}(V^*)\otimes \wedge^{i+1}V^*;\label{E:deRham-deg}\\
d^*\left(Sym^p(V^*)\otimes \wedge^iV^*\right)\subset Sym^{p+1}(V^*)\otimes \wedge^{i-1}V^*;\label{E:deRham-star-dg}\\
\mbox{ for }\ome\in Sym^p(V^*)\otimes \wedge^i(V^*) \mbox{ one has }
(dd^*+d^*d)\ome=(p+i)\ome.
\end{eqnarray}
Thus if $d\ome=0$ then $dd^*\ome=(p+i)\ome$. Hence if $p+i>0$ then $\ome=\frac{1}{p+i}dd^*\ome$ is exact; this implies (1).
Part (2) also follows from the last formula and (\ref{E:deRham-star-dg}).
\qed

\begin{lemma}\label{L:de-rham-deg-corol}
Let $(\Ome_{pol}^\bullet(V),d)$ be the polynomial de Rham complex on a complex finite dimensional vector space $V$. Let $W$ be another complex finite dimensional vector space.
Consider the complex $\left(\Ome^\bullet_{pol}(V)\otimes Sym^\bullet(W^*),d\right)$ where $d$, by the abuse of notation, denotes  the only $Sym^\bullet(W^*)$-linear
extension of the de Rham differential. Let $A\in \ZZ_{\geq 0}$. Then the natural imbedding of complexes
$$\left(\oplus_{p,i,q:\,p+i+q\leq A}\left(Sym^pV^*\otimes\wedge^iV^*\otimes Sym^q(W^*)\right),d\right)\hookrightarrow (\Ome^\bullet_{pol}(V)\otimes Sym^\bullet(W^*),d)$$
is a quasi-isomorphism, namely induces isomorphism on the cohomology.
\end{lemma}
{\bf Proof.} Since $d$ preserves the degree in $Sym^\bullet(W^*)$ the lemma follows from the following immediate
consequence of Lemma \ref{L:pol-de-rham}: for fixed $A,q\in \ZZ_{\geq 0}$ the cohomology of the complex
$$\left(\oplus_{p,i:\,p+i\leq A-q}Sym^p(V^*)\otimes\wedge^iV^*\otimes Sym^q(W^*),d \right)$$
vanished in the positive cohomological degrees (with respect to the index $i$) and equals $Sym^q(W^*)$ in the zero degree. \qed

\hfill

{\bf Proof of Theorem \ref{T:relative-deRham}.} Since the statement is local, we may assume that the bundle $\cl$ is trivial of rank one.
By Lemma \ref{L:submflds} it suffices to prove the following statement. Let $\car$ be a smooth complex manifold. Let $C\subset \car$ be a closed
complex analytic submanifold. Let $V$ be an $m$-dimensional complex vector space, and $B^m\subset V$ a centered ball. We will identify $\car$ with
$\car\times\{0_V\}\subset \car\times V$. Let $i_n\colon C^{(n)}\hookrightarrow \car\times B^m $ be the $n$th infinitesimal neighborhood of $C$ inside $\car\times B^m$.
Let $\tilde C^{(n)}\subset \car$ be the $n$th infinitesimal neighborhood of $C$ inside $\car$. Let $p_1\colon \car\times B^m\to \car$ be the natural projection.
Let $p_1^{(n)}\colon C^{(n)}\to \tilde C^{(n)}$ be the restriction of $p_1$ to $C^{(n)}$. Under these notation and assumptions, our theorem will immediately follow
from the exactness of the complex
\begin{eqnarray}\label{E:simple-de-rham}
0\to (p_1^{(n)})^{-1}\co_{\tilde C^{(n)}}\to i_{n*}i_n^*\co_{\car\times B^m}\to i_{n-1*}i_{n-1}^*\Ome^1_{(\car\times B^m)/\car}\to\\
i_{n-2*}i_{n-2}^*\Ome^2_{(\car\times B^m)/\car}\to \dots.
\end{eqnarray}

It is clear that after making $\car$ smaller we may assume that $\car$ is a open subset of $\CC^N=W\times \CC^l$,
where $W:=\CC^{N-l}$ and $C=(\{0_W\}\times \CC^l)\cap \car$.

Then it is easy to see that
\begin{eqnarray*}
(p_1^{(n)})^{-1}\co_{\tilde C^{(n)}}=\co_C\otimes_\CC \left(\oplus_{q=0}^n Sym^qW^*\right);\\
i_{n-i*}i_{n-i}^*\Ome^i_{(\car\times B^m)/\car}=\co_C\otimes_\CC \left(\oplus_{p+q\leq n-i}Sym^pV^*\otimes Sym^qW^*\right)\otimes\wedge^iV^*.
\end{eqnarray*}
Hence the exactness of the complex (\ref{E:simple-de-rham}) follows from immediately from Lemmas \ref{L:de-rham-deg-corol} and \ref{L:pol-de-rham} after tensoring
the complexes there with $\co_C$ over $\CC$. \qed

\subsection{Few technical lemmas on infinitesimal neighborhoods of a subvariety.}\label{Ss:lemmas-on-infi}
We start with the following
\begin{lemma}\label{L:2-nbhd-twistors}
Let $(M^{4n},I,J,K)$ be a hypercomplex manifold. Let $Z=\CC\PP^1\times M$ denote the twistor space of $M$. Let $a\in M$ be a point, and
$C:=\CC\PP^1\times \{a\}$ be the complex curve in $Z$.  Let $M_0:=(\HH^n,I_0,J_0,K_0)$ be the standard flat hypercomplex space, and let $Z_0:=\CC\PP^1\times M_0$
denote its twistor space. Let $C_0:=\CC\PP^1\times \{0\}$ be the complex curve in $Z_0$. Then there exists an isomorphism of complex analytic spaces
of the second infinitesimal neighborhoods of these curves in the ambient twistor spaces
$$C^{(2)}\simeq C_0^{(2)}.$$
\end{lemma}
{\bf Proof.} Making $M$ smaller we may assume that $M=(U,I,J,K)$ where $U\subset \HH^n$ is an open neighborhood of 0, and the
complex structures $I,J,K$ satisfy
\begin{eqnarray}\label{E:complex-str}
I(x)=I_0+O(x^2),
\end{eqnarray}
and similarly for $J,K$. Hence considering $\CC\PP^1\times U$ as the twistor space of either $M$ or $M_0$ we get two complex structures on it
denoted respectively by $\cj$ and $\cj_0$. Their construction and (\ref{E:complex-str}) imply that for any $z_0\in \CC\PP^1\times \{0\}$
and any $z\in \CC\PP^1\times U$ one has
\begin{eqnarray}\label{E:c-structures}
\ci(z)=\ci_0(z)+O((z-z_0)^2).
\end{eqnarray}

But as previously we denote the curve $\CC\PP^1\times \{0\}$ by $C$ if it is considered as a curve in the complex manifold $Z$, and by $C_0$
if it is considered as a curve in $Z_0$. Similarly $C^{(2)}$ and $C_0^{(2)}$ are the second infinitesimal neighborhoods of
$C$ and $C_0$ in $\CC\PP^1\times U$ equipped with the appropriate complex structure.

Let $\underline{C}^\infty(\CC\PP^1\times U)$ denote the sheaf of smooth $\CC$-valued functions on the manifold $\CC\PP^1\times U$.
Let $\ci$ be the subsheaf of ideals of the sheaf $\underline{C}^\infty(\CC\PP^1\times U)$ consisting of functions vanishing on
$\CC\PP^1\times \{0\}$. The sheaves of regular (holomorphic) functions on $C^{(2)}$ and $C_0^{(2)}$ are subsheaves of the sheaf
$\underline{C}^\infty(\CC\PP^1\times U)/\ci^3$ (the latter can be informally considered as the sheaf of smooth functions on either
$C^{(2)}$ or $C_0^{(2)}$). Let us show that these two subsheaves coincide; clearly that will imply the lemma.

Let $\underline{\Ome}^1$ be the sheaf of 1-forms on $\CC\PP^1\times U$. We have two morphisms of sheaves induced by the
Dolbeault differentials
corresponding to the complex structures $\cj,\cj_0$:
\begin{eqnarray}\label{E:partial}
\bar\pt,\bar\pt_0\colon \underline{C}^\infty(\CC\PP^1\times U)/\ci^3\to \underline{\Ome}^1/\ci^2\cdot\underline{\Ome}^1.
\end{eqnarray}

The subsheaf of regular functions on $C^{(2)}$ (resp. $C_0^{(2)}$) is the kernel of the operator $\bar \pt$ (resp. $\bar\pt_0$) from (\ref{E:partial}).
Hence it suffices to show that the operator
$$\bar\pt-\bar\pt_0\colon \underline{C}^\infty(\CC\PP^1\times U)/\ci^3\to \underline{\Ome}^1/\ci^2\cdot\underline{\Ome}^1$$
vanishes. This easily follows from (\ref{E:c-structures}). \qed

\hfill

Let $Z$ be a complex analytic space. Let $C\subset Z$ be a compact analytic subspace.
We denote by $[C]$ the corresponding closed point in the Douady space $Dou(Z)$. Let $n$ be a non-negative integer.
Let $C^{(n)}$ be the $n$th infinitesimal neighborhood of $C$ in $Z$; thus $C^{(n)}\subset Z$ is a closed analytic subspace.
Hence by Proposition \ref{P:Inclusion-Douady} $i\colon Dou(C^{(n)})\inj Dou(Z)$ is a closed analytic subspace.
\begin{lemma}\label{L:closed-of-Dou}
The $n$th infinitesimal neighborhood $[C]^{(n)}$ of the closed point $[C]\in Dou(Z)$ is a closed analytic subspace
of $Dou(C^{(n)})$ (which, in turn, a closed subspace of $Dou(Z)$). Moreover $[C]^{(n)}$ coincides with the $n$th infinitesimal neighborhood
of $[C]$ in $Dou(C^{(n)})$.
\end{lemma}
{\bf Proof.} First, tautologically,
$i\in Mor([C]^{(n)},Dou(Z))$. Let us denote by $A$ the corresponding closed analytic subspace
$$A\subset Z\times [C]^{(n)}$$
which is flat and proper over $[C]^{(n)}$, and such that the fiber in $A$ of the point $[C]\in [C]^{(n)}$ is equal
to $C$. We claim that in fact $A$ is a closed analytic subspace of $C^{(n)}\times [C]^{(n)}\subset Z\times [C]^{(n)}$ (the latter
inclusion is also a closed imbedding of analytic spaces obviously); this is a special case of the next lemma.
\begin{lemma}\label{L:infinitesimal-inclus}
Let $Z_1,T$ be affine complex analytic schemes. Assume moreover that the ring of regular functions $\co[T]$ is a finite dimensional (over $\CC$)
local ring with maximal ideal $\frak{m}$ such that $\frak{m}^{n+1}=0$. (It follows that $\co[T]/\frak{m}=\CC$.)

Let $A_1\subset Z_1\times T$ be a closed analytic subspace. Let $C_1=A_1\times_T T_0$ where $T_0$ is the subspace of $T$ whose algebra of functions is the field $\CC$,
i.e. the only point of $T$.

Then $A_1$ is a closed subscheme of $C_1^{(n)}\times T$, where $C_1^{(n)}$ denotes the $n$th infinitesimal neighborhood of $C_1$ in $Z_1$.
\end{lemma}
Let us postpone the proof of this lemma and finish the proof of Lemma \ref{L:closed-of-Dou}.
Thus $A\subset C^{(n)}\times [C]^{(n)}$ is a closed analytic subspace. This corresponds to a morphism
$\tilde i\in Mor([C]^{(n)},Dou(C^{(n)}))$. It has the property that $i=j\circ\tilde i$ where $j\colon Dou(C^{(n)})\to Dou(Z)$
is the natural closed imbedding. Since $i$ and $j$ are closed imbeddings, it follows that $\tilde i$ is also a closed imbedding. This implies the first part of
Lemma \ref{L:closed-of-Dou}.

Let us prove the second part. We prove a slightly more general statement. Let $A$ be a closed analytic subspace of an analytic space $B$, and let $p\in A$ be a point.
Let $p^{(n)}$ be the $n$th infinitesimal neighborhood of $p$ in $B$. Assume that $p^{(n)}$ is a closed analytic subspace of $A$. Then the claim is that the
$n$th infinitesimal neighborhood of $p$ in $A$ is equal to $p^{(n)}$. We may and will assume that $A,B$ are affine analytic spaces. Let $\ca$ and $\cb$ denote the
algebras of functions on $A$ and $B$ respectively. Let $\cj\subset \cb$ be the ideal such that $\ca=\cb/\cj$. Let $\frak{m}\subset\cb$ be the ideal of $p$.
Since $p^{(n)}$ is a closed analytic subspace of $A$ then
\begin{eqnarray}\label{E:inclusion-ideals}
\cj\subset \frak{m}^{n+1}.
\end{eqnarray}
We have the obvious epimorphism of algebras $\cb\to \ca/(\frak{m}/\cj)^{n+1}$. To finish the proof of the lemma one has to show that its kernel
is equal to $\frak{m}^{n+1}$. But this follows from the inclusion (\ref{E:inclusion-ideals}).
Thus Lemma \ref{L:closed-of-Dou} is proved. \qed

\hfill


{\bf Proof of Lemma \ref{L:infinitesimal-inclus}.}
Let $I_{A_1}\subset \co[Z_1\times T]=\co[Z_1]\otimes_\CC \co[T]$
be the ideal defining $A_1$. Let $J_{C_1}\subset \co[Z_1]$ be the ideal defining $C_1$.
We have to show that
\begin{eqnarray}\label{E:show1}
(J_{C_1}\otimes_\CC \co[T])^{n+1}\subset I_{A_1}.
\end{eqnarray}

First we have
\begin{eqnarray}\label{E:show2}
\co[C_1]=\co[Z_1]/J_{C_1}=(\co[Z_1]\otimes \co[T])/(J_{C_1}+\frak{m}).
\end{eqnarray}
On the other hand since $C_1$ is the fiber over the closed point of the projection $A_1\to T$, we also have
\begin{eqnarray}\label{E:show3}
\co[C_1]=(\co[Z_1]\otimes \co[T])/(I_{A_1}+\frak{m}).
\end{eqnarray}
Then (\ref{E:show2}) and (\ref{E:show3}) imply
\begin{eqnarray*}
J_{C_1}+\frak{m}=I_{A_1}+\frak{m}.
\end{eqnarray*}
In particular any $x\in J_{C_1}$ has a presentation
$$x=a+y \mbox{ with } a\in I_{A_1},\, y\in \frak{m}.$$
We have
$$x^{n+1}=(a+y)^{n+1}=\left(\sum_{l=1}^{n+1}\binom{n+1}{l}a^ly^{n+1-l}\right)+y^{n+1}.$$
The first sum in brackets belongs to $I_{A_1}$, while $y^{n+1}=0$. Thus $x^{n+1}\in I_{A_1}$ and (\ref{E:show1}) is proved.
Hence Lemma \ref{L:infinitesimal-inclus} is proved, and Lemma \ref{L:closed-of-Dou} follows too. \qed


\subsection{Identification of a Baston operator with $\pt\pt_J$.}\label{Basto=d-dj}
The main result of this section is
\begin{theorem}\label{T:Baston}
Consider the Penrose transform of $\cl=\xi^*(\co(-2))$.
Then the Baston operator
$$\Delta\colon C^\infty(M)\to \wedge^{(2,0)}_IM$$
is proportional to $\pt\pt_J$ with a non-zero constant of proportionality.
\end{theorem}
To prove this theorem we will need some preparations. Let us first prove this theorem in the flat case.
\begin{lemma}\label{L:Baston-flat case}
Theorem \ref{T:Baston} is true for the flat space $(\HH^n,I_0,J_0,K_0)$.
\end{lemma}
{\bf Proof.} By Proposition 7.1(i) in \cite{alesker-jgp} the Baston operator $\Delta$ can be identified
(up to a normalizing constant) with the quaternionic Hessian $\Delta'$ introduced in \cite{alesker-bsm-03}.
The latter can be identified (again after a normalization) with $\pt\pt_J$ (see \cite{alesker-verbitsky-JGA-06}, Proposition 4.1). \qed

\hfill

Now we are going to reduce the general case of Theorem \ref{T:Baston} to the flat one. Fix $a\in M$. By Corollary \ref{Cor:Obata-cor} there exists a real analytic diffeomoprhism of
a neighborhood of $a\in M$ with a neighborhood of $0\in\HH^n$ such that $a$ goes to $0$ and under such an identification
$$I(x)=I_0+O(x^2),$$
and similarly for $J,K$. As we have seen in Lemma \ref{L:2-nbhd-twistors}, this induces an isomorphism of complex analytic spaces
$$\phi\colon C^{(2)}\tilde\to C_0^{(2)},$$
where $C^{(2)}\subset Z$ (resp. $C_0^{(2)}\subset Z_0$) is the second infinitesimal neighborhood of
$C=\CC\PP^1\times\{a\}$ (resp. $C_0=\CC\PP^1\times\{0\}$) in the twistor space $Z$ of $M$ (resp. $Z_0$ of $\HH^n$).

This isomorphism $\phi$ induces obviously isomorphism
$$\tilde\phi\colon Dou(C^{(2)})\tilde\to Dou(C_0^{(2)})$$
such that $\tilde\phi([C])=[C_0]$. Then $\tilde\phi$ induces isomorphism of the second infinitesimal
neighborhood of $[C]$ in $Dou(C^{(2)})$ with the second infinitesimal neighborhood of $[C_0]$ in $Dou(C_0^{(2)})$.
Then by Lemma \ref{L:closed-of-Dou} they coincide with the second infinitesimal neighborhoods of $[C]$ in $Dou(Z)$ and of $[C_0]$ in $Dou(Z_0)$ respectively,
namely we get
\begin{claim}\label{Cl:tilde-phi-iso}
(i) $\tilde\phi$ induces an isomorphism $[C]^{(2)}\tilde\to [C_0]^{(2)}$.

(ii) For any $0\leq l\leq 2$ the restriction of $\tilde\phi$ to $[C]^{(l)}$ is an isomorphism
$$\tilde\phi_l\colon [C]^{(l)}\tilde\to [C_0]^{(l)}.$$
\end{claim}
\begin{remark}
Observe that for $0\leq l\leq 2$ the scheme $[C]^{(l)}$ (resp. $[C_0]^{(l)}$) coincides with the $l$th infinitesimal neighborhood of
$[C]$ (resp. $[C_0]$) in $Dou(C^{(2)})$ (resp. $Dou(C_0^{(2)})$) This again follows from Lemma \ref{L:closed-of-Dou}.
\end{remark}


As in Theorem \ref{T:scheme+point} for any analytic space $Z$ we will denote by $F(Z)$ the analytic space representing the contravariant functor (denoted again by $F(Z)$)
from analytic spaces to sets given by
\begin{eqnarray*}\label{Def:scheme+pt}
F(Z)(S)=\{A\subset Z\times S,\, f\colon S\to A|\\ A\mbox{ is closed analytic subspace, flat and proper over } S,\\
f \mbox{ is a section of } A\to S\}.
\end{eqnarray*}
\begin{lemma}\label{L:fiber-product}
Let $Y$ be a closed analytic subspace of an analytic space $Z$. Then
$$F(Y)=F(Z)\times_{Dou(Z)}Dou(Y).$$
\end{lemma}
{\bf Proof.} For an analytic space $W$ we will also denote below by the same letter $W$ also
the contravariant functor from analytic spaces to sets given by
$$S\mapsto Mor(S,W).$$
Then we have
\begin{eqnarray*}
\left(F(Z)\times_{Dou(Z)}Dou(Y)\right)(S)=F(Z)(S)\times_{Dou(Z)(S)}Dou(Y)(S)=\\
\left\{A\subset Z\times S, \,f\colon S\to A \mbox{ as above}\right\}\times_{\{A\subset Z\times S\mbox{ as above}\}}
\{A\subset Y\times S\mbox{ as above}\}=\\
\{A\subset Y\times S,f\colon S\to A\mbox{ as above}\}=F(Y)(S).
\end{eqnarray*}
Lemma is proved. \qed

\hfill

\def\ttphi{\tilde{\tilde\phi}}

Lemma \ref{L:fiber-product} implies in particular
\begin{eqnarray}\label{E:fiber-prod-Dou}
F(C^{(2)})=F(Z)\times_{Dou(Z)}Dou(C^{(2)}).
\end{eqnarray}
Since $\phi\colon C^{(2)}\tilde \to C_0^{(2)}$ is an isomorphism, it induces an isomorphism
$$\tilde{\tilde \phi}\colon F(C^{(2)})\tilde\to F(C_0^{(2)}).$$
Using (\ref{E:fiber-prod-Dou}) it is equivalent to
$$\tilde{\tilde\phi}\colon F(Z)\times_{Dou(Z)}Dou(C^{(2)})\tilde\to  F(Z_0)\times_{Dou(Z_0)}Dou(C_0^{(2)}).$$
Clearly the following diagram commutes

\begin{eqnarray}
\square<2`1`1`2;1300`500>[F(Z)\times_{Dou(Z)}Dou(C^{(2)})`F(Z_0)\times_{Dou(Z_0)}Dou(C_0^{(2)})`Dou(C^{(2)})`Dou(C_0^{(2)});\ttphi`\tau`\tau_0`\tilde\phi],
\end{eqnarray}
where the vertical arrows are the natural projections.

Recall that $[C]^{(2)}\subset Dou(C^{(2)})$, $[C_0]^{(2)}\subset Dou(C_0^{(2)})$ are closed imbedding, and
for $[C]^{(l)}\subset [C]^{(2)}$, $[C_0]^{(l)}\subset [C_0]^{(2)}$ are also closed imbedding for $0\leq l\leq 2$. Hence making
a base change, Claim \ref{Cl:tilde-phi-iso} implies that $\phi$ induces isomorphisms
\begin{eqnarray*}
\hat\phi_l\colon F(Z)\times_{Dou(Z)}[C]^{(2)}\tilde\to F(Z_0)\times_{Dou(Z_0)}[C_0]^{(l)},\ 0\leq l\leq 2.
\end{eqnarray*}
Clearly the following diagrams commutes for $0\leq l\leq 2$:
\begin{eqnarray}\label{Diam:diagram-A}
\square<2`1`1`2;1300`500>[F(Z)\times_{Dou(Z)}[C]^{(l)}`F(Z_0)\times_{Dou(Z_0)}[C_0]^{(l)}`[C]^{(l)}`[C_0]^{(l)};\hat\phi_l`\tau`\tau_0`\tilde\phi_l],
\end{eqnarray}
\begin{eqnarray}\label{Diam:diagram-B}
\square<2`1`1`2;1300`500>[F(Z)\times_{Dou(Z)}[C]^{(l)}`F(Z_0)\times_{Dou(Z_0)}[C_0]^{(l)}`C^{(l)}`C_0^{(l)};\hat\phi_l`\eta`\eta_0`\phi_l],
\end{eqnarray}
where $\phi_l\colon C^{(l)}\tilde\to C_0^{(l)}$ is the restriction of $\phi$ to $C^{(l)}$; clearly it is an isomorphism.
Note also that $[C]^{(1)}$ is a closed analytic subspace of $[C]^{(2)}$ and $\phi\colon [C]^{(1)}\tilde\to [C_0]^{(1)}$ is an isomorphism.

\hfill

Before we formulate our next lemma let us note that for any open subset $U\subset C^{(n)}$ and for $f\in H^0(U,\co_{C^{(n)}})$ the differential
$df$ is well defined as an element of $H^0(U,i_{n-1}^*\Ome^1)$, i.e. the order of the infinitesimal neighborhood is smaller by 1 (in the last expression $U$ is also considered as an
open subset of $C^{(n-1)}$.
Analogous statement is true for the relative de Rham differential:
$$d_{F/Z}f\in H^0(U,i_{n-1}^*\Ome^1_{F/Z}).$$
\begin{lemma}\label{L:deRham-commute}
Let $F\overset{\eta}{\rightarrow}Z$, $F_0\overset{\eta_0}{\rightarrow}Z_0$
be holomorphic submersive maps of complex analytic smooth manifolds (not necessarily coming from quaternionic manifolds).
Let $\tilde C\subset F$, $\tilde C_0\subset F_0$ be smooth complex submanifolds such that $\eta(\tilde C)$ and $\eta_0(\tilde C_0)$
are smooth complex submanifolds of $Z$ and $Z_0$ respectively. Assume in addition that $\eta\colon\tilde C\to \eta(\tilde C)$ and $\eta_0\colon\tilde C_0\to \eta(\tilde C_0)$ are isomorphisms.

Let for a fixed $n\in \NN$ $$\psi\colon (\eta(C))^{(n)}\tilde\to (\eta(C_0))^{(n)},\tilde\psi\colon \tilde C^{(n)}\tilde\to \tilde C_0^{(n)}$$
be isomorphisms of complex analytic spaces such that the diagram commutes
$$
\square<1`1`1`1;700`500>[\tilde C^{(n)}`\tilde C_0^{(n)}`(\eta(\tilde C))^{(n)}`(\eta(\tilde C_0))^{(n)};\tilde\psi`\eta`\eta_0`\psi].
$$

Then the following holds.

(1) For $0\leq k\leq n$ the restriction of $\tilde\psi$ to $\tilde C^{(k)}$ is an isomorphism
$\tilde C^{(k)}\tilde\to \tilde C_0^{(k)}$; this restriction will be denoted by $\tilde\psi$ again.

(2) There exists an isomorphism $\Psi_{n-1}$ of sheaves of graded algebras of $\co_{\tilde C^{(n-1)}}$-modules
$$\Psi_{n-1}\colon \tilde\psi^*i_{n-1}^*\Ome_{F_0/Z_0}^\bullet\tilde\to i_{n-1}^*\Ome_{F/Z}^\bullet$$
which is uniquely characterised by the following additional property:
\footnote{Here $i_{n-1}$ denotes two closed imbeddings $\tilde C^{(n-1)}\inj Z$ and $\tilde C^{(n-1)}_0\inj Z_0$.}

$\bullet$ for any open set $U\subset C_0^{(n-1)}$, any $f\in H^0(U,\co_{\tilde C_0^{(n-1)}})$ and any $g\in H^0(U,\co_{\tilde C_0^{(n)}})$ one has
$$\Psi_{n-1}(\tilde\psi^*(f\cdot d_{F_0/Z_0}g))=\tilde\psi^* f\cdot d_{F/Z}(\tilde\psi^*g).$$

(3) Moreover $\Psi_{n-1}$ is compatible with the relative de Rham differential: for $\ome\in H^0(U,i_{n-1}^*\Ome_{F_0/Z_0}^\bullet)$ one has
$$\Psi_{n-2}(\tilde\psi^*(d_{F_0/Z_0})\ome)=d_{F/Z}\left(\Psi_{n-1}(\tilde\psi^*\ome)\right).$$
(Note that in the left hand side $d_{F_0/Z_0}\ome\in H^0(U,\co_{\tilde C^{(n-2)}})$, and hence $\tilde\psi$ has the meaning of the
isomorphism $\tilde C^{(n-1)}\tilde\to \tilde C_0^{(n-1)}$, and $\Psi_{n-2}$ is the isomorphism whose existence was claimed in part (2),
but for $n-2$ instead of $n-1$.
\end{lemma}
We leave the proof of the lemma to the reader.
Let us apply Lemma \ref{L:deRham-commute} in the situation of Theorem \ref{T:Baston}. We denote $\tilde C=\tau^{-1}([C])$, then $\tilde C^{(n)}=F(Z)\times_{Dou(Z)}[C]^{(n)}$; similarly for $C_0$ instead of $C$.
Then we obtain the following commutative diagram of sheaves:
\begin{eqnarray}\label{E:isomor-comlexes}
\bfig \putsquare<1`1`1`1;900`500>(0,500)[\phi_2^*\co_{\tilde C_0^{(2)}}`\phi_1^*i_{1*}i_1^*\Ome^1_{F_0/Z_0}`\co_{\tilde C^{(2)}}`i_{1*}i_1^*\Ome^1_{F/Z};
d_{F_0/Z_0}`id`\Psi_1`d_{F/Z}]\putsquare<1`0`1`1;900`500>(900,500)[\phantom{\phi^*_1i_{1*}i_1^*\Ome^1_{F_0/Z_0}}`\phi^*_0i_{0*}i^*_0\Ome^2_{F_0/Z_0}`
\phantom{i_1*i_1^*\Ome^1_{F/Z}}`i_{0*}i_0^*\Ome^2_{F/Z};d_{F_0/Z_0}``\Psi_0`d_{F/Z}]\efig
\end{eqnarray}
where the vertical arrows are isomorphisms, $i_n$ denotes both the closed imbeddings $\tilde C^{(n)}=F(Z)\times_{Dou(Z)}[C]^{(n)}\inj Z$ and
$\tilde C^{(n)}_0=F(Z_0)\times_{Dou(Z_0)}[C_0]^{(n)}\inj Z_0$.

Applying $R\tau_*$ and $R\tau_{0*}$ to the first and second rows of the diagram (\ref{E:isomor-comlexes}) respectively and using Corollary \ref{Cor:grauert-twistor}
we get a spectral sequence for each complex in either row. Since the diagram (\ref{E:isomor-comlexes}) commutes, we get isomorphisms
of the first terms of these spectral sequences commuting with the second differentials which we will denote respectively
by $\Delta^{(2)}$ and $\Delta_0^{(2)}$:
\begin{eqnarray*}
\square<1`1`1`1;700`500>[\co_{[C_0]^{(2)}}`\Lambda^{(2,0)}_I\HH^n|_0`\co_{[C]^{(2)}}`\Lambda_I^{(2,0)}M|_a;\Delta^{(2)}_0`\phi^*`d\phi(a)^*`\Delta^{(2)}],
\end{eqnarray*}
where in the second column we have isomorphism of vector spaces of $(2,0)$-forms at the point $a\in M$ and $0\in \HH^n$ respectively induced by the pull-back under the differential
$d\phi(a)$.

To finish the proof of Theorem \ref{T:Baston} it remains to show that $\Delta^{(2)}$ and $\Delta_0^{(2)}$ are restrictions of $\Delta$ and $\Delta_0$ to 2-jets of functions at $a$ and $0$ respectively.
Since the second case is a special case of the first one, let us do it for $\Delta^{(2)}$. We have the following commutative diagram:
\begin{eqnarray}\label{Diagr:dR-vs-infinitesim}
\bfig \putsquare<1`1`1`1;700`500>(0,500)[\co_{F(Z)}`\Ome^1_{F/Z}`\co_{F(Z)\times_{Dou(Z)}}C^{(2)}`i_{1*}i_1^*\Ome^1_{F/Z};d_{F/Z}```d_{F/Z}]
\putsquare<1`0`1`1;700`500>(700,500)[\phantom{\Ome^1_{F/Z}}`\Ome^2_{F/Z}`\phantom{i_{1*}i_1^*\Ome^1_{F/Z}}`i_0^*i_{0*}\Ome^2_{F/Z};d_{F/Z}```d_{F/Z}]
\putsquare<1`0`1`1;700`500>(1400,500)[\phantom{\Ome^2_{F/Z}}`\Ome^3_{F/Z}`\phantom{i_0^*i_{0*}\Ome^2_{F/Z}}`0;d_{F/Z}```]
\putsquare<1`0`0`1;500`500>(2100,500)[\phantom{\Ome^3_{F/Z}}`\dots`\phantom{0}`\dots;```]\efig
\end{eqnarray}
where the vertical arrow are the adjunction morphisms; they form a morphism of complexes staying in the rows of the last diagram.
Thus applying $R\tau_*$ we get a morphism of spectral sequences, in particular we get the commutative diagram for the second differential:
\begin{eqnarray}\label{Diag-p28}
\square<1`1`1`1; 800`500>[\co_{Dou(Z)}`\tau_*\Ome^2_{F/Z}`\co_{[C]^{(2)}}`\Lambda^{(2,0)}_IM_a;\Delta```\Delta^{(2)}],
\end{eqnarray}
where the first vertical arrow is the obvious map of taking the second jet of a function, and the second
vertical map is the evaluation map at the point $a$. Thus we see that $\Delta^{(2)}$ is the restriction
of $\Delta$ to 2-jets of functions at the point $a$. Theorem \ref{T:Baston} is proved. \qed

\hfill

Before we finish this section let us prove another lemma about hypercomplex manifold we will use later.
\begin{lemma}\label{L:star-lemma}
Let us assume that the twistor space $Z$ of a hypercomplex manifold $M^{4n}$ is strongly 2-complete.
Assume that $Dou(Z)$ is a Stein manifold.
Assume also that the fibers of the map
$\eta\colon F(Z)\to Z$ are contractible.
Let $a\in M$, $C=\CC\PP^1\times\{a\}\subset Z$. Let $\kappa\colon C^{(2)}\inj Z$ denote the natural imbedding.
Let $\cl=\xi^*(\co(-2))$.

Then there exists an isomorphism of vector spaces $h\colon H^1(C^{(2)},\kappa^*\cl)\tilde\to Ker \Delta^{(2)}$ such that the following diagram commutes:
\begin{eqnarray}\label{magend}
\square<1`1`1`1; 800`500>[H^1(Z,\cl)`Ker(\Delta)`H^1(C^{(2)},\kappa^*\cl)`Ker \Delta^{(2)};```h],
\end{eqnarray}
where the first horizontal arrow is the isomorphism from Theorem \ref{T:Baston-cohomology}(2), the first vertical arrow is induced by the adjunction
morphism $Id\to \kappa_*\kappa^*$, and the second vertical arrow sends function to its 2-jet at the point $a$.
\end{lemma}
{\bf Proof.} Let $\tilde C:=\tau^{-1}(a)$. Let us denote
$$\eta^{(2)}\colon \tilde C^{(2)}\to C^{(2)}$$
the restriction of $\eta$ to $\tilde C^{(2)}$. Let $i_n\colon \tilde C^{(n)}\to F(Z)$ be the natural imbedding. We have the commutative diagram of sheaves on $F(Z)$:
\begin{eqnarray}\label{Diag:long}
\bfig \putsquare<1`0`1`1;450`500>(0,500)[0`\eta^{-1}\cl`0`i_{2*}(\eta^{(2)})^{-1}(\kappa^*\cl);```]
\putsquare<1`0`1`1;650`500>(450,500)[\phantom{\eta^{-1}\cl}`\eta^*\cl`\phantom{i_{2*}(\eta^{(2)})^{-1}(\kappa^*\cl)}`i_{2*}i_2^*(\eta^*\cl);```]
\putsquare<1`0`1`1;800`500>(1100,500)[\phantom{\eta^*\cl}`\eta^*\cl\otimes\Ome^1_{F/Z}`\phantom{i_2^*i_{2*}(\eta^*\cl)}`i_{1*}i_1^*(\eta^*\cl\otimes \Ome^1_{F/Z});```]
\putsquare<1`0`1`1;880`500>(1900,500)[\phantom{\eta^*\cl\otimes\Ome^1_{F/Z}}`\eta^*\cl\otimes\Ome^2_{F/Z}`\phantom{i_{1*}i_1^*(\eta^*\cl\otimes \Ome^1_{F/Z})}`i_{0*}i_0^*(\eta^*\cl\otimes \Ome^2_{F/Z});```]
\putsquare<1`0`1`1;650`500>(2780,500)[\phantom{\eta^*\cl\otimes\Ome^2_{F/Z}}`\eta^*\cl\otimes \Ome^3_{F/Z}`\phantom{i_{0*}i_0^*(\eta^*\cl\otimes \Ome^2_{F/Z})}`0;```]
\putsquare<1`0`0`1;440`500>(3430,500)[\phantom{\eta^*\cl\otimes \Ome^3_{F/Z}}`\dots`\phantom{0}`\dots;```]
\efig
\end{eqnarray}
where the vertical arrows are the adjunction morphisms, and the horizontal arrows are the relative de Rham differentials. The first row in this diagram is exact as the
relative de Rham complex. The second row is exact by Theorem \ref{T:relative-deRham}. Let us denote
\begin{eqnarray*}
A^\bullet=[0\to \eta^*\cl\to \eta^*\cl\otimes \Ome^1_{F/Z}\to \eta^*\cl\otimes\Ome^2_{F/Z}\to\dots];\\
B^\bullet=[0\to i_{2*}i^*_2(\eta^*\cl)\to i_{1*}i_1^*(\eta^*\cl\otimes \Ome^1_{F/Z})\to i_{0*}i_0^*(\eta^*\cl\otimes\Ome^2_{F/Z})\to 0\to \dots].
\end{eqnarray*}

Thus $A^\bullet$ and $B^\bullet$ are resolutions of $\eta^{-1}\cl$ and $i_{2*}(\eta^{(2)})^{-1}(\kappa^*\cl)$ respectively.
We have a sequence of canonical maps which are isomorphisms in our case
\begin{eqnarray*}
H^1(Z,\cl)\tilde\to H^1(F(Z),\eta^{-1}\cl)\tilde\to \HH^1(F(Z),A^\bullet)=H^1(R\Gamma\circ R\tau_*(A^\bullet)),
\end{eqnarray*}
where $R\Gamma$ is the derived functor of taking global sections over $Dou(Z)$; the first map is an isomorphism since the fibers of $\eta$ are contractible.
The composition of these isomorphisms we will denote by
$$\alp\colon H^1(Z,\cl)\tilde\to H^1(R\Gamma\circ R\tau_*(A^\bullet)).$$

Next, since $\eta|_{\tilde C}\colon \tilde C\to C$ is an isomorphism, the natural map
\begin{eqnarray}\label{E:for-B1}
H^1(Z,\kappa_*\kappa^*\cl)\simeq H^1(C^{(2)},\kappa^*\cl)\to H^1(\tilde C^{(2)},(\eta^{(2)})^{-1}(\kappa^*\cl))
\end{eqnarray}
is an isomorphism. Furthermore the canonical map
\begin{eqnarray}\label{E:for-B2}
H^1(\tilde C^{(2)}, (\eta^{(2)})^{-1}(\kappa^*\cl))\to H^1(F(Z),i_{2*}(\eta^{(2)})^{-1}(\kappa^*\cl)).
\end{eqnarray}
is also an isomorphism for trivial reason. But
\begin{eqnarray}\label{E:for-B3}
H^1(F(Z),i_{2*}(\eta^{(2)})^{-1}(\kappa^*\cl))\simeq H^1(R\Gamma\circ R\tau_*(B^\bullet)).
\end{eqnarray}

Composing the isomorphisms (\ref{E:for-B1})-(\ref{E:for-B3}) we obtain an isomorphism
\begin{eqnarray*}
\bar\alp\colon H^1(Z,\kappa_*\kappa^*\cl)\tilde\to H^1(R\Gamma\circ R\tau_*(B^\bullet)).
\end{eqnarray*}

It is easy to see that the following diagram commutes:
\begin{eqnarray}\label{Diag:cohomology1}
\bfig \putsquare<1`1`1`1;900`400>(0,400)[H^1(Z,\cl)`H^1(R\Gamma\circ R\tau_*(A^\bullet))`H^1(Z,\kappa_*\kappa^*\cl)`H^1(R\Gamma\circ R\tau_*(B^\bullet));\alp```\bar\alp]\efig
\end{eqnarray}
where the first vertical arrow is the adjunction morphism, and the second vertical arrow is induced by the morphism of complexes
$\Phi\colon A^\bullet\to B^\bullet$ from diagram (\ref{Diag:long}).

By Corollary \ref{Cor:grauert-twistor}
\begin{eqnarray}\label{E:semes-1}
R^q\tau_*(i_{n*}i_n^*(\eta^*\cl\otimes \Ome^p_{F/Z}))=(R^q\tau_*(\eta^*\cl\otimes \Ome^p_{F/Z}))\otimes_{\co_{Dou(Z)}}\co_{a^{(n)}}.
\end{eqnarray}
But by Proposition \ref{P:push-identification} and the computation after it we have
\begin{eqnarray}\label{E:semes-2}
R^q\tau_*(\eta^*\cl\otimes \Ome^p_{F/Z})=0 \mbox{ unless } p=0 \mbox{ and } q=1, \mbox{ or }
p\geq 2 \mbox{ and } q=0.
\end{eqnarray}
By (\ref{E:semes-1}) and (\ref{E:semes-2}) we see that the complexes $A^\bullet$, $B^\bullet$ and the functor $R\tau_*$ satisfy the assumptions of
Lemma \ref{L:homol-alg}. Hence it follows that there exist isomorphisms in $D^+(Sh_{Dou(Z)})$ (the notation below is obvious)
\begin{eqnarray*}
\beta\colon R\tau_*(A^\bullet)\tilde\to \BC_2(A^\bullet)[-1];\\
\bar\beta\colon R\tau_*(B^\bullet)\tilde\to \BC_2(B^\bullet)[-1]
\end{eqnarray*}
and a morphism $\Psi\colon \BC_2(A^\bullet)\to \BC_2(B^\bullet)$ again in $D^+(Sh_{Dou(Z)})$ such that
the following diagram commutes:
\begin{eqnarray}\label{E:semes-3}
\bfig \putsquare<1`1`1`1;900`500>(0,500)[R\tau_*(A^\bullet)`\BC_2(A^\bullet)[-1]`R\tau_*(B^\bullet)`\BC_2(B^\bullet)[-1];\beta`R\tau_*(\Phi)`\Psi[-1]`\bar\beta]\efig
\end{eqnarray}
Let us apply to all elements of the last diagram the functor $R\Gamma$ and then apply the functor $H^1$ of the first cohomology.
\footnote{Since $\BC_2(A^\bullet),\BC_2(B^\bullet)$ consist of complexes of coherent sheaves on $Dou(Z)$, and $Dou(Z)$ is a Stein manifold by
assumption, it follows that $R\Gamma(\BC_2(A^\bullet))$ and $R\Gamma(\BC_2(B^\bullet))$ are obtained by taking global sections of $A^\bullet$ and $B^\bullet$
respectively.} Then combining this with the diagram (\ref{Diag:cohomology1}) we deduce the lemma. \qed

\subsection{End of Proof of Theorem \ref{Penrose key proposition}.}\label{Ss:End-of-proof}
By Theorem \ref{T:Baston} the operator
$$\pt\pt_J\colon C^\infty(M)\to \Lam_I^{(2,0)}(M)$$
is proportional to the Baston operator $\Delta$ corresponding to the Penrose transform
of the line bundle $\cl=\xi^*(\co(-2))$.
Hence it suffices to show that there exists $\eps>0$ such that for any $z_0\in \mathcal{C}$
and any 2-jet $\bar f$ of a function at $z_0$ such that $\Delta\bar f(z_0)=0$ there exists a real analytic function $f$ in the open ball
$B(z_0,\eps)$ with respect to the metric $\rho$ such that
$$\Delta f\equiv 0$$
and the 2-jet of $f$ at $z_0$ is equal to $\bar f$.

\hfill

By Proposition \ref{P:twistors-2-complete} there exists $\eps_1>0$ such that the twistor spaces of
the balls $B(z_0,\eps)\subset M$ are strongly 2-complete for any $\eps\in (0,\eps_1]$ and any $z_0\in \mathcal{C}$.

\hfill

It is possible to extend the Riemannian metric $\rho$ to a Riemannian metric on an open neighborhood of $M$ inside the Douady space of the
twistor space $Z(M)$ of $M$ (this neighborhood is denoted again by $Dou(Z(M))$. We denote this extension by $\rho$ again. For $x\in Dou(Z(M))$ and $\delta>0$ we denote by $B_{Dou}(x,\delta)$
the open ball in
$Dou(Z(M))$ centered at $x$ and of radius $\delta$. It is easy to see that one can choose $\eps_2>0$ such that for any $\eps\in (0,\eps_2]$ and any
$z_0\in \mathcal{C}\subset M\subset Dou(Z(M))$ the ball $B_{Dou}(z_0,\eps)$ is a Stein manifold. Indeed the function $dist^2(\cdot, z_0)$ is an exhausting plurisubharmonic function on
this ball.

\hfill

Recall that by Section \ref{Ss:Douady} we have two submersions
\begin{eqnarray*}
\eta\colon F(Z(M))\to Z(M),\\
\tau\colon F(Z(M))\to Dou(Z(M)),
\end{eqnarray*}
where $F(Z(M))\subset Z(M)\times Dou(Z(M))$ is a closed subset which is a smooth submanifold. The map $\tau$ is proper.

By Proposition \ref{P:contractible} there exists $\eps_3>0$ such that any $\eps\in (0,\eps_3]$ has the following property:
for any $z_0\in \mathcal{C}$ and any $b\in \CC\PP^1\times B(z_0,\eps)\subset Z(M)$ the fibers of $\eta$
$$F(Z(M))\cap \left(\{b\}\times B_{Dou}(z_0,\eps)\right)$$
are contractible.

\hfill

Let us fix now $\eps\in (0,\min\{\eps_1,\eps_2,\eps_3\}]$.
Fix any $z_0\in \mathcal{C}$. Let us replace $M$ with the ball $B(z_0,\eps)$; we denote it by $M$ again.
Let us choose $Dou(Z(B(z_0,\eps)))$ to be the ball $B_{Dou}(z_0,\eps)$, we denote it by $Dou (Z)$ where $Z$ is the twistor space
of $M=B(z_0,\eps)$. Hence we can apply Theorem \ref{T:Baston-cohomology} for $l=-2$ and get the Baston complex and the cohomological interpretation of $Ker \Delta\simeq H^1(Z,\cl)$.

\hfill


With this notation, we are going to show that for any 2-jet $\bar f$ of a function at $z_0$
such that
$$\Delta \bar f(z_0)=0$$
there exists a holomorphic $f$ in $Dou(Z)$ such that
$$\Delta f \equiv 0$$
and the 2-jet of $f$ at $z_0$ equals $\bar f$. Let $C=\CC\PP^1\times\{z_0\}$, and let $C^{(2)}$ denote as previously the second infinitesimal neighborhood of $C$ in $Z$, and let $\kappa\colon C^{(2)}\inj Z$
denote the natural imbedding.
Then by Lemma \ref{L:star-lemma} it suffices to show that the adjunction map
$$H^1(Z,\cl)\to H^1(C^{(2)}, \kappa^*\cl)=H^1(Z,\kappa_*\kappa^*\cl)$$
is onto. This follows from the strong 2-completeness of $Z$. Indeed let $\ck$ be the kernel of the adjunction epimorphism $\cl\to \kappa_*\kappa^*\cl$.
$\ck$ is a coherent sheaf, hence for $i\geq 2$ one has $H^i(Z,\ck)=0$ by the strong 2-completeness of $Z$.
From the long exact sequence in cohomology we have
$$H^1(Z,\cl)\to H^1(Z,\kappa_*\kappa^*\cl)\to H^2(Z,\ck)=0$$
Theorem \ref{Penrose key proposition} is proven. \qed

\appendix

\section{Appendix. Proof of Proposition \ref{P:Inclusion-Douady} (by Daniel Barlet).}

\subsection{The Douady space of a subspace.}
In this appendix the Douady space of a complex analytic space $X$ is denoted by $D(X)$ (and not by $Dou(X)$ as in the main text).
The proof of the proposition 4.2.2 will use the following lemma:

\begin{lemma}\label{trivial}
Let $f : M \to N$ be a morphism of complex spaces (not necessarily reduced, but locally finite dimensional). Then if the Zariski tangent map
 $$Tf_{x} : T_{M,x} \to T_{N,f(x)}$$
  is injective, then $f$ is  proper embedding of an open neighbourhood of $x \in M$ in an open neighbourhood of $f(y)$ in $N$.
\end{lemma}

The standard proof  is left to the reader.

\begin{lemma}\label{proper}
Let $f : M \to N$ be a morphism of complex spaces (not necessarily reduced, but locally finite dimensional). Assume that
\begin{enumerate}
\item  the map  $f$ is injective ;
\item the map  $f$ is a  local embedding near each point $x \in M$ ;
\item the set theoretic  image $f(M)$ is a closed analytic subset of $N$.
\item for any holomorphic germ $\gamma : (\C,0) \to (f(M), f(x))$ there exists a holomorphic germ $ \tilde{\gamma} : (\C,0) \to (M,x)$ , so that $\gamma = f\circ\tilde{\gamma}$.
\end{enumerate}
Then $f$ is a proper embedding of $M$ in $N$, that is to say that $f$ is an isomorphism of $M$ on the subspace $\mathcal{M}$ of $N$ defined on $f(M)$ by the sheaf $f_{*}(\mathcal{O}_{M})$ which is isomorphic to the quotient of $\mathcal{O}_{N}$ by the coherent ideal $I \subset \mathcal{O}_{N} $ of germs $h$ such that $h\circ f = 0$ in $\mathcal{O}_{M}$.
\end{lemma}

{\bf Proof.} Let $y = f(x)$ be a point in $f(M)$. Then there exists open neighbourhoods $U_{x}$ of $x$ in $M$ and $V_{y}$ of $y$ in $N$  such that $f$ induces a proper embedding of $U_{x}$ in $V_{y}$. Assuming that $V_{y}$ is small enough, then the analytic subset $f(M) \cap V_{y}$ has all its irreducible components passing through $y$. Assume that there is such an irreducible component which is not contained in $f(U_{x})$. Then we can find a germ of holomorphic curve $\gamma :  (\C,0) \to (f(M), f(x))$ such that only the point $\gamma(0) = y$ of this curve is in $f(U_{x})$. But this contradicts the condition 4. So we conclude that, for any $y = f(x)$ we can find $U_{x}$ and $V_{y}$ such that $f(Y)\cap V_{y} = f(U_{x}) \cap V_{y}$. And of course we can assume that $U_{x}$ is a relatively compact set in $M$.\\
Take any compact set $K$ in $N$ and for any $y = f(x) \in K \cap f(M)$ choose such $U_{x}$ and $V_{y}$. Extract a finite sub-covering of the open covering $(V_{y})_{y \in K \cap f(M)}$. Then let $y_{i} := f(x_{i}), i \in [1,N]$ the corresponding index of the sub-covering and define $L := \cup_{i=1}^{N} \ \bar U_{x_{i}}$. Then, thanks to the injectivity of $f$  we have $f^{-1}(K) \subset L$ which shows that $f$ is proper. As it is an  injective map and a local embedding this is a proper embedding.$\hfill \blacksquare$

\bigskip

{\bf Proof of the Proposition 4.2.2.} First recall that if $I := I_{X}$ is an ideal of $\mathcal{O}_{Z}$ defining a point $X$  in $D(Z)$ the Douady's space of $Z$, the Zariski tangent space at $I$ of $D(Z)$ is the  (finite dimensional) vector space
$$H^{0}(Z, (I/I^{2})^{*}) $$
where the star denotes the dual sheaf $Hom_{\mathcal{O}_{Z}}(-, \mathcal{O}_{Z})$. Note that, as $X$ is a compact subspace, the support of the sheaf $(I/I^{2})^{*}$ is compact, so this vector space is finite dimensional.

For the convenience of the reader  we shall sketch a proof  (see \cite{douady} p.77) :\\

Denote by $S := \{z \in \C \ / \  z^{2} = 0\}$ which is the non reduced complex space associated to the $\C-$algebra $Spec\big(\C[\varepsilon]\big/(\varepsilon^{2})\big)$. Then, by definition, the Zariski tangent space at the point  $X$ of the space $D(Z)$ corresponding to the coherent ideal sheaf $I$, where $\mathcal{O}_{X} := \mathcal{O}_{Z}\big/ I$ has a compact support, is given by the subset $Mor(S,D(Z))_{X}$  of the set $Mor(S,D(Z))$ of morphisms with ``value'' the point $X$. This corresponds to the coherent ideals $J \subset \mathcal{O}_{Z}\otimes_{\C}A $, where $A := \C[\varepsilon]\big/(\varepsilon^{2})$, inducing $\mathcal{O}_{Z}\big/ I \simeq \mathcal{O}_{X}$ by the quotient map $ id\otimes q :  \mathcal{O}_{Z}\otimes_{\C}A \to  \mathcal{O}_{Z}$, where $q := A \to \C$ is defined by $q(\varepsilon) = 0$, such that $ \mathcal{O}_{Z}\otimes_{\C}A\big/J$ is $A-$flat.\\
Recall that a $A-$module $M$ is $A-$flat if and only if we have the equivalence for $m \in M$
$$ \varepsilon.m = 0 \quad  {\rm  if \ and \ only\  if} \quad  m \in \varepsilon.M .$$
To explicite the natural bijection between $Mor(S,D(Z))_{X}$ and $H^{0}(Z, (I/I^{2})^{*})$ consider an ideal $J \subset \mathcal{O}_{Z}\otimes_{\C}A $ with image $I$ by the quotient map $id\otimes q$, such that $ \mathcal{O}_{Z}\otimes_{\C}A\big/J$ is $A-$flat. Then  $\varepsilon.I$  is contained in $J$ and, by flatness of $J$ as a $A-$module, we have $\varepsilon.I = \varepsilon.J$. And, for each $i \in I$, there exists $k \in \mathcal{O}_{Z}$ such that $i + \varepsilon.k$ lies in $J$. If we have also $i + \varepsilon.k'$ in $J$, we find that $\varepsilon.(k-k') \in J$. The flatness of the quotient by $J$ gives $k-k' \in I$. So for each $i \in I$ there exists an unique class $[k] \in \mathcal{O}_{Z}\big/I$ such that $i + \varepsilon.k' \in J $ for any $ k' \in [k]$. This defines a $\mathcal{O}_{Z}-$linear map of sheaves $\varphi : I \to \mathcal{O}_{Z}\big/I$ such that any element in $J$ is of the form $i + \varepsilon.k$ where $k \in [\varphi(i)]$. This defined our map
$$Mor(S,D(Z))_{X} \to H^{0}(Z, Hom_{\mathcal{O}_{Z}}(I, \mathcal{O}_{Z}\big/I)) \simeq H^{0}(Z, Hom_{\mathcal{O}_{Z}}(I\big/I^{2}, \mathcal{O}_{Z})).$$
Conversely, for such a $\mathcal{O}_{Z}-$linear morphism  $\varphi : I \to \mathcal{O}_{Z}\big/I$ define the  ideal $J$ as the set of $i + \varepsilon.k$ where $k \in [\varphi(i)]$. Then it is easy to verify that $\mathcal{O}_{Z}\otimes_{\C}A\big/J $ is $A-$flat and has $I$ for image by $id\otimes q$. $\hfill \blacksquare$

\bigskip

 If $Y$ is a closed complex subspace in $Z$, then $\mathcal{O}_{Y}$ is a quotient of $\mathcal{O}_{Z}$, and if $J := J_{X}$ is the ideal in $\mathcal{O}_{Y}$ image of the ideal $I := I_{X}$ of $X \subset Z$ by the quotient map, we have a surjective map of coherent $\mathcal{O}_{Z}-$sheaves (with compact supports)
$$ I/I^{2} \to J/J^{2} \to 0 $$
which induces an injection of the dual sheaves and then an injection
$$ H^{0}(Z, (J/J^{2})^{*}) \to H^{0}(Z, (I/I^{2})^{*})$$
 which shows that the natural map\footnote{Which is deduced from the universal property of $D(Z)$.} $D(Y) \to D(Z)$ induces an injection between  the Zariski tangent spaces at $I$ and $J$ and so is a local embedding. This gives the condition 2. in Lemma \ref{proper}.\\
 To see that it is a proper embedding, as injectivity is obvious, it is enough to prove that the image of $D(Y) $ is closed subset in $D(Z)$, that this image is a locally closed analytic subset of $D(Z)$ and that the condition 4. of Lemma \ref{proper} is satisfied.  \\

 To prove that the image is closed is not so easy : Let $X \in D(Z)$ such that $X \not\subset Y$. Then there exists a point $x \in X$ and a holomorphic function on an open neighbourhood $U$ of $x \in Z$, such that $f \in I_{X}$ is not a section of $I_{Y}$ on $U$. By flatness over $D(Z)$ we can extend $f$ to a neighbourhood $S \times U'$ of $(X,x)$ in $D(Z)\times Z$ with $x \in U'  \subset U$ and $S$ an open neighbourhood of $X$ in $D(Z)$,  to a section $\tilde{f}$ of the ideal $\mathcal{I}_{\mathcal{X}}$, where $\mathcal{X} \subset D(Z)\times Z$ is the graph of the universal family, and $\tilde{f}$ is not in $I_{Y}$ along $\{s\}\times U'$ for any $s \in S$\footnote{Remark that to make this argument completely explicit, following Douady, we can consider a privileged compact set $K$ for  $Y$  in a neighbourhood of $x$ in $Z$ and consider the multiplication  by $\tilde{f}$ acting on the trivial holomorphic banach bundle $S \times B(K, \mathcal{O}_{Y})$ ; see the argument below.}  . Then any $s \in S$ is not in the image of $D(Y)$ in $D(Z)$ and this image is closed.\\

 Unfortunately, in order to show that the image of $D(Y )$ in $D(Z)$ is a locally closed analytic subset in D(Z), I will have to assume that the reader has some ideas about the construction of Douady. That may be the reason why it is not easy to find a reference for this result in the literature.\\
 Points in $D(Z)$  are locally described by holomorphic banach (locally trivial) fiber bundles $F$ over an open set $S$ in $D(Z)$ which are (direct) quotients of the trivial banach bundle $S \times B(K, \mathcal{O}_{Z})$ and which are $B(K, \mathcal{O}_{Z})-$modules, where $K$ is a privileged compact set in $Z$ (for the notion of privileged compact and the notation $B(K, \mathcal{O}_{Z})$ see \cite{douady}). The main point is to see that if $f$ is a holomorphic function in a neighbourhood of $K$ the subset  of points $s$ in $S$  for which the morphism $\tilde{f}$ induced by $f$ vanishes of the fiber $F(s)$ of the bundle $F$ is an analytic subset in $S$. This is rather obvious as it is the same to ask that the fiber at $s$ of the kernel of the quotient map $S \times B(K, \mathcal{O}_{Z}) \to F$ contains the closed ideal $\{s\} \times B(K, f.\mathcal{O}_{Z})$ because this kernel is also a locally trivial holomorphic banach bundle over $S$ of $B(K, \mathcal{O}_{Z})-$modules. This implies that adding theses banach analytic conditions in Douady's construction, we defined locally an analytic subset of $D(Z)$ which is exactly the set of points corresponding to $X \in D(Z)$ which are subspaces of $Y$.\\
 So we have shown that conditions 1.-3. of the Lemma \ref{proper} are satisfied for the natural map $i : D(Y) \to D(Z)$. Let us show that the condition 4. is also satisfied :\\
 If $\gamma : (\C,0) \to i(D(Y),i(X))$ is a germ of holomorphic curve, the pull-back by $\gamma$ of the universal family over $D(Z)$ is a flat and proper family. By definition of $i(D(Y))$, the points of this (closed analytic) subset corresponds exactly to the subspaces $X \in D(Z)$ which are subspaces of $Y$. So the pull-back family is in fact a flat and proper family of subspaces in $Y$, and the universal property of $D(Y)$ gives the desired lifting $\tilde{\gamma} : (\C,0) \to (D(Y),X)$ of $\gamma$; so the condition 4. of Lemma \ref{proper} is satisfied. This complete the proof.$\hfill \blacksquare$\\

\begin{remark} Note that if $Y$ has global equations in $Z$\footnote{In fact global equations around each $X$ for each point $X $ in the open set of $D(Z)$ we want to consider would be enough.} we can avoid to go back in Douady's construction because we can consider the global  sheaf map
 $$\times f :  p_{*}(\mathcal{O}_{\mathcal{X}}) \longrightarrow  p_{*}(\mathcal{O}_{\mathcal{X}}) $$
 where $f : Z \to \C$ is a global holomorphic function vanishing on $Y$, $\mathcal{X} \subset D(Z) \times Z$ is the graph of universal family on $D(Z)$ and $p : D(Z)\times Z \to D(Z)$ the natural projection. Then, as the sheaf $p_{*}(\mathcal{O}_{\mathcal{X}})$ is coherent and flat on $D(Z)$, it is a locally free sheaf so it is the sheaf of sections of a vector bundle $F$ on $D(Z)$. Then we can consider the analytic subset $D(Z)_{f}$ of $D(Z)$ of points  $s \in D(Z)$ such that the map induced by $\times f$ on the fiber $F(s)$ vanishes. If $f_{1}, ..., f_{N}$ are global holomorphic equations for $Y$ in $Z$ the intersection of the analytic subsets $D(Z)_{f_{i}}$ for $i \in [1,N]$ will define the  set theoretical image of $D(Y)$ in $D(Z)$.\\
 Remark that the argument above is the same than the "general'' one sketched before, but using a vector bundle to replace the locally trivial holomorphic banach bundle.
 \end{remark}


\subsection{A general remark on the Douady space.}
    For any complex space $Z$ there exists a natural holomorphic map
   $$ \pi : red(D_{n}(Z)) \to \mathcal{C}_{n}(Z) $$
   of the reduced Douady's space of pure $n-$dimensional subspaces to the space of compact $n-$cycles in $Z$ which is also a reduced complex space (see \cite{barlet-75} chapter VI for this map). An interesting fact for the study of twistor's space is that when $Z$ is a complex manifold and where we restrict the map $\pi$ to the open set  $\Omega \subset D(Z)$ of reduced smooth complex compact $n-$dimensional submanifolds of $Z$ the map $\pi$ is  an isomorphism of reduced complex space $red(\Omega)$ onto an open set in $\mathcal{C}_{n}(Z)$. So it is possible to use in this case some results known of the space $\mathcal{C}_{n}(Z)$. For instance, Bishop's theorem (see
   \cite{barlet-magnusson}, Ch.IV) allows to give the following description of compact subset in $\mathcal{C}_{n}(Z)$ :\\
   A subset $\mathcal{K} \subset \mathcal{C}_{n}(Z) $ is compact if and only the following properties are satisfied :
   \begin{enumerate}
   \item There exists a compact subset $K \subset Z$ such that any $X \in \mathcal{K}$ has its support $\vert X\vert$ contained in $K$.
   \item There exists a continuous hermitian metric $\omega$ on $Z$ and a constant $C(\mathcal{K},\omega)$ such that for any $X \in \mathcal{K}$ we have
   $$ \int_{X} \ \omega^{n} \leq C(\mathcal{K},\omega) .$$
   \end{enumerate}
   Note that the condition 2. does not depend of the choice of the continuous hermitian metric $\omega$ when the condition 1. is satified.\\
   This gives, for instance, the properness missing in page 17 :\\
   Let $p : Z \to M$ the twistor map ; it is proper. Now the injection $ i : M \to D_{1}(Z)$ given by the fibers of $p$ is also proper because if $red(\Omega) \subset red(D_{1}(Z)) $ is an open set in the reduced Douady's space of $Z$ consisting of reduced smooth compact curves in $Z$, as $\pi$ induces an isomorphism of $red(\Omega)$ on the corresponding open set  in $\mathcal{C}_{1}(Z)$, it is enough to prove the properness of the map $\pi\circ i$. If $\mathcal{K}$ is a compact set in $red(\Omega)$, there exists, thanks to the characterisation above, a compact set $K$ in $Z$ such that any $X \in \mathcal{K}$ has its support in $K$. This implies that
   $$( i\circ\pi)^{-1}(\mathcal{K}) \subset p(K) $$
   which is a compact set in $M$ as $p$ is continuous. This  gives the properness of the map $i$ on an open neighbourghood of its image.

\end{document}